\documentclass{amsbook}
\usepackage[utf8]{inputenc}
\usepackage{amsfonts}
\usepackage{amssymb}
\usepackage{amsthm}
\usepackage{mathtools}
\usepackage{enumerate,xcolor,hyperref}
\newtheorem{theorem}{Theorem}[chapter]
\newtheorem{lemma}[theorem]{Lemma}
\newtheorem{remark}[theorem]{Remark}
\newtheorem{definition}[theorem]{Definition}
\newtheorem{corollary}[theorem]{Corollary}
\newtheorem{proposition}[theorem]{Proposition}
\newtheorem{example}[theorem]{Example}
\newtheorem*{non-theorem}{Theorem}
\numberwithin{equation}{chapter}
\numberwithin{table}{chapter}
\numberwithin{figure}{chapter}

\newcommand{\im}{\mbox{\rm im\,}}
\newcommand{\re}{\mbox{\rm re\,}}
\newcommand{\IM}{\mbox{\rm Im\,}}
\newcommand{\RE}{\mbox{\rm Re\,}}
\newcommand{\tr}{\mbox{\rm \,tr\,}}
\newcommand{\diag}{\mbox{\rm diag}}

\newcommand{\rank}{\mbox{\rm rank}}
\newcommand{\Arg}{\mbox{\rm Arg }}

\newcommand{\Real}{\mathbb{R}}
\newcommand{\Comp}{\mathbb{C}}
\newcommand{\eps}{\varepsilon}
\newcommand{\conv}{\mbox{\rm conv}}

\DeclareMathOperator{\ii}{i}

\DeclareMathOperator{\Span}{span}

\newcommand{\eq} [1] {\begin{equation}\label{#1}}
\newcommand{\en} {\end{equation}}
\newcommand {\eqn}  {\begin{eqnarray}}
\newcommand {\enn}  {\end{eqnarray}}
\newcommand {\bstar}    {\begin{eqnarray*}}
\newcommand {\estar}    {\end{eqnarray*}}
\newcommand {\mat}  [1] {\left[\begin{array}{#1}}
\newcommand {\rix}      {\end{array}\right]}
\newcommand{\norm}[1]{\left\| #1 \right\|}
\newcommand{\set}[1]{\left\{ #1 \right\}}

\catcode`@=11     
\font\tenex=cmex10 
\newdimen\p@renwd
\setbox0=\hbox{\tenex B} \p@renwd=\wd0 
\def\bmat#1{\begingroup \m@th
  \setbox\z@\vbox{\def\cr{\crcr\noalign{\kern2\p@\global\let\cr\endline}}%
    \ialign{$##$\hfil\kern2\p@\kern\p@renwd&\thinspace\hfil$##$\hfil
      &&\quad\hfil$##$\hfil\crcr
      \omit\strut\hfil\crcr\noalign{\kern-\baselineskip}%
      #1\crcr\omit\strut\cr}}%
  \setbox\tw@\vbox{\unvcopy\z@\global\setbox\@ne\lastbox}%
  \setbox\tw@\hbox{\unhbox\@ne\unskip\global\setbox\@ne\lastbox}%
  \setbox\tw@\hbox{$\kern\wd\@ne\kern-\p@renwd\left[\kern-\wd\@ne
    \global\setbox\@ne\vbox{\box\@ne\kern2\p@}%
    \vcenter{\kern-\ht\@ne\unvbox\z@\kern-\baselineskip}\,\right]$}%
  \null\;\vbox{\kern\ht\@ne\box\tw@}\endgroup}
\catcode`@=12 

\def\rank{\mathop{\mathrm{rank}}}

\def\Span{\mathop{\mathrm{span}}}
\def\diag{\mathop{\mathrm{diag}}}

\newcommand {\comment}[1]{} 
\author{Oskar Jakub Szymański}

\begin{document}
    \begin{titlepage}
         \begin{center}
        \vspace*{1cm}

        \textbf{STABILITY THEORY FOR MATRIX POLYNOMIALS IN ONE AND SEVERAL VARIABLES WITH EXTENSIONS OF CLASSICAL THEOREMS}
            
        \vspace{1.5cm}

        Oskar Jakub Szymański

        \vspace{0.5cm}
        PhD Dissertation
        \vfill
            
       \raggedleft{Supervisor: \\  dr hab. Michał Wojtylak, prof. UJ}
       \vspace{0.5cm}
        \end{center}
       \vspace{0.8cm}
        \begin{center} 
        \includegraphics[width=100pt]{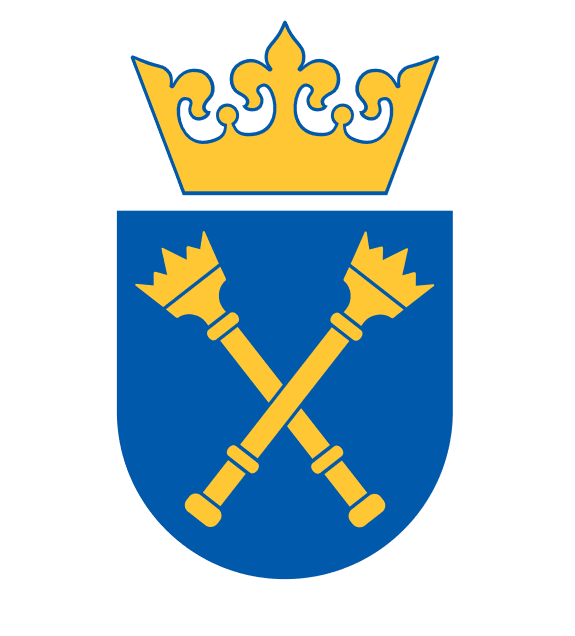}\\
       Jagiellonian University\\
       Faculty of Mathematics and Computer Science\\
       Cracow 2024
        \end{center}
    \end{titlepage}

\tableofcontents

\section*{Preface}
In 
2019 professor Greg Knese from the Washington University in St. Louis published a paper \\

{\small \noindent
G. Knese, Global Bounds on Stable Polynomials,  {\em Complex Analysis and Operator Theory}, 13(4):1895–1915, 2019, no. \cite{Kne19} in the Biblilography.\\}

\noindent The publication contains several generalisations of the classical Sz\'asz inequality
    \begin{equation}\label{ineq}
|p(\lambda)| \leq \exp\bigl(|a_1||\lambda| + 3(|a_1|^2 + |a_2|)|\lambda|^2\bigr)\text{,}\quad \lambda \in \mathbb{C}\text{.}
    \end{equation}
The inequality holds for stable scalar polynomials $p(\lambda) = a_d\lambda^d + \dots + a_1\lambda + 1$ of one complex variable, i.e. polynomials without the roots in the open upper half-plane $H_0$, satisfying the condition $p(0) = 1$. Knese, using a determinantal representation, extended it to scalar mulit-variate stable polynomials.
At the end of 2019, shortly after enrolling in the PhD Program at the Jagiellonian University, I was asked by Professor Łukasz Kosiński to present the aforementioned results at the seminar 'Geometric Function Theory'. This was a starting point of my research on the theory of polynomials stability. Another source of motivation, that came directly afterwards, was the publication \\

{\small \noindent
Julius Borcea and Petter Br\"and\'en, The Lee-Yang and P\'olya-Schur programs. I. Linear operators preserving stability. {\em Inventiones Mathematicae}, 177(3):541, 2009, no. \cite{BorB09} in the Bibliography.} \\

\noindent The authors of this paper present a complete characterization of linear operators that preserve stability of scalar polynomials. They introduce polarisation operators, which play a crucial role in the stability theory.

As noted by my supervisor, Professor Michał Wojtylak, 
the methods of  \cite{Kne19} and \cite{BorB09} could be successfully applied in the area of matrix polynomials, i.e., polynomials of the form
$$
P(\lambda) = \lambda^d A_d + \lambda^{d-1} A_{d-1} + \dots + A_0,
$$
where $A_0, A_1, \dots , A_d$ are complex square matrices and $\lambda$ is a complex indeterminate. 
A reason of this observation was a growing interest in the theory of matrix polynomials (in particular: in localising their spectra) and in their connections with systems of differential equations. The basis for the current research were the publications: \\

{\small \noindent
C. Mehl, V. Mehrmann, and M. Wojtylak, Distance problems for dissipative hamiltonian systems and related matrix polynomials. {\em Linear Algebra Applications}, 623:335–366, 2021, no. \cite{MehMW21} in the Bibliography,} \\

{\small \noindent
C. Mehl, V. Mehrmann, and M. Wojtylak, Matrix pencils with coefficients that have positive semi-definite Hermitian parts. {\em SIAM Journal on Matrix Analysis and Applications}, 43(3):1186–1212, 2022, no. \cite{MehMW22} in the Bibliography.} \\

\noindent Hence, I decided to turn my attention into this direction.
The main subject of the Thesis are regular matrix polynomials with eigenvalues (i.e. zeros of the function $\lambda \mapsto \det P(\lambda)$) localised in a certain set. Usually the half-planes and the unit disc are of main interest. Currently there exist several different methods of localizing eigenvalues of matrix polynomials. Below we concentrate on the one using the numerical range, see \cite{LiR94,Psa03, Psa00}. It appears that the numerical range contains the eigenvalues, which provides a theoretical tool for localisation. However, this method is rather restrictive. In the joint work with my supervisor we have developed several new concepts, which turn out to be more convenient to localise eigenvalues. In particular, we have introduced a notion of {\em hyperstability with respect to $D$} for matrix polynomials: for all $x \in \mathbb{C}^n\setminus\{0\}$ there exists $y \in \mathbb{C}^n\setminus\{0\}$ such that for all $\lambda \in D$ we have $y^*P(\lambda)x \neq 0$, see Definition ~\ref{defla}.
On one hand, if the numerical range lies outiside $D$, then the polynomial is hyperstable with respect to $D$. On the other hand, hyperstability with respect to $D$ implies that there are no eigenvalues in $D$,  cf. Proposition ~\ref{abc}. 

Further, we have discussed matrix version of classical results in complex analysis: Gauss-Lucas theorem, and the aforementioned  Sz\'asz inequality. 
The former one generalises nicely using our hyperstability concept:
\begin{non-theorem}
 Let $D\subseteq\Comp$ be a nonempty open or closed set such that $\Comp\setminus D$ is convex. If a matrix polynomial $P(\lambda)$ is hyperstable with respect to $ D$ and the entries of its derivative $P'(\lambda)$ are linearly independent polynomials over $\mathbb{C}$, then the matrix polynomial $P'(\lambda)$ is also hyperstable with respect to $D$.  
\end{non-theorem}

Generalising the Sz\'asz inequality ~\eqref{ineq} onto matrix polynomials is more difficult and requires a structural analysis. 
In particular, we present examples showing that the dependence on the degree $d$ and dimension $n$ is indispensable, cf. Examples ~\ref{ones} and ~\ref{cmv}.
In connection with this we have showed that there exist several global upper bounds for the Frobenius norm of $p(A)$, where a polynomial $p(\lambda)$ satisfies the assumptions of original Sz\'asz inequality and $A \in \mathbb{C}^{n, n}$. For a comparison of these bounds see Example ~\ref{comp}. \\

Let us discuss the content of the Thesis.
The first Chapter presents a review of the known facts, which form a basis for the Dissertation. 

The second Chapter introduces hyperstability for one- and multi-variable matrix polynomials and deals with its basic properties. Example ~\ref{exa} shows a significant difference between stability and hyperstability. The third Chapter is devoted to an extension of the Gauss-Lucas Theorem on matrix polynomials. It presents the possibilities which hyperstability gives us for improving classical theorems, cf. Theorem ~\ref{GLmat} and Theorem ~\ref{mmgl}. The fourth Chapter offers a variety of Sz\'asz-type inequalities in the matrix case, for example: 
$$
\| P(\lambda) \| \leq 2\exp \left( \lambda_{H}\Bigl[\lambda A_1 -|\lambda|^2 A_2\Bigr] +\frac12 |\lambda|^2 \| A_1\|^2       \right), \quad \lambda \in \Comp,
$$
where $P(\lambda) = \lambda^d A_d + \dots + \lambda A_1 + I_n$ is a matrix polynomial with the numerical range contained in some half plane $H_{\varphi}$.

The fifth Chapter shows how to obtain one-variable hyperstable polynomial from two-variable stable quadratic or cubic polynomial. These 'hyperstability via stability' methods are included in Theorems ~\ref{poly2} (presented below) and ~\ref{poly3}. 
\begin{non-theorem}
Let $P(\lambda) = \lambda^2 A_2 + \lambda A_1 + A_0$ be a  quadratic matrix polynomial and let $D$ be a nonempty open or closed subset of the complex plane $\mathbb{C}$. If at least one of the following conditions holds:
	\begin{enumerate}[\rm (a)]
\item\label{0?D} the multivariate matrix polynomial  $(z_1, z_2) \mapsto z_1^2 A_2 + z_2 A_1 + A_0$ is stable with respect to $D^2$,
\item\label{0notinD1} the  multivariate matrix polynomial $(z_1, z_2) \mapsto z_1z_2 A_2 + z_2 A_1 + A_0$  is stable with respect to $D^2$ and $0 \notin D$,
\item\label{0notinD2} the  multivariate matrix polynomial  $(z_1, z_2) \mapsto z_1^2z_2 A_2 + z_1^2 A_1 + z_2 A_0$ is stable with respect to $D^2$ and $0 \notin D$,
	\end{enumerate}
then the matrix polynomial $\lambda \mapsto P(\lambda)$ is hyperstable with respect to $D$.
\end{non-theorem}
The sixth Chapter describes different operators preserving matrix hyperstability - from basic operators, see Proposition ~\ref{bop}, to polarisation ones, see Theorem ~\ref{Tkappa2}. In the last Section of this Chapter we mention about polarisation operators acting on singular matrix polynomials. The seventh Chapter provides instances of interesting classes of stable and hyperstable matrix polynomials. In this final Chapter, we apply many theorems from previous Chapters and demonstrate how they work, see for example Proposition ~\ref{subadd} or Theorem ~\ref{half-plane}.

All results of Chapters 2 to 7 are original results, most of the content of Chapters 2,3,5,7 was published in \\

{\small \noindent
Oskar Jakub Szymański and Michał Wojtylak, Stability of matrix polynomials in one and several variables. {\em Linear Algebra and its Applications}, 670:42–67, 2023, no. \cite{szymanski2023stability} in the Bibliography.} \\

In the current Thesis the theory is extended by additional examples.
The content of Chapter 4 is announced in \\

{\small \noindent
Piotr Pikul, Oskar Jakub Szymański and Michał Wojtylak, The Sz\'asz inequality for matrix polynomials and functional calculus. {\em arXiv preprint arXiv:2406.08965, 2024}, no. \cite{PSW} in the Bibliography.} \\

The remaining results, i.e., orbits of hyperstability (Theorem ~\ref{orbits}) and a multi-variable version of Gauss-Lucas theorem for hyperstable polynomials (Theorem ~\ref{mmgl}) are in preparation for publication. \\

Hopefully, the newly introduced notion of hyperstability offers many possibilities in analysing regular matrix polynomials. It gives an opportunity to extend classical theorems such as the Gauss-Lucas theorem and the Sz\'asz inequality. My results can be applied in many areas which use localising eigenvalues of matrix polynomials, for example in the control theory of dynamic systems. Therefore, I think that concepts discussed in my Thesis are up-to-date and can be developed much further. One of them are orbits of hyperstability, described extensively in Chapter ~\ref{sHyper}. \\

Finally, I wish to express my gratitude towards my supervisor Prof. Michał Wojtylak for his patience and priceless clues in my research. I also would like to thank a lot Prof. Łukasz Kosiński for guiding me on the hard path of mathematics. Their assistance has proved essential in my development as a research scientist.

\newpage

\section*{Notations}

\begin{enumerate}[]
\item $\mathbb{N} := \{0, 1, 2, \dots\},\quad \mathbb{Z}_+ := \{1, 2, 3, \dots\}$
\item $\Comp^{n,n}[\lambda] := \{P(\lambda) = \sum_{j=0}^d \lambda^j A_j : A_0, A_1, \dots , A_d \in \mathbb{C}^{n, n}\}$
\item $D^{\perp} := \{x \in \mathbb{C}^n : y^*x = 0 \;\text{for all}\; y \in D\}$,\; $D \subseteq \mathbb{C}^n$
\item $\sigma(A)$ - the set of eigenvalues of $A \in \mathbb{C}^{n, n}$
\item $\norm x := \sqrt{x^*x},\; x \in \mathbb{C}^n$
\item $\norm A := \sqrt{\max\sigma(A^*A)}$,\; $A \in \mathbb{C}^{m, n}$
\item $\norm A_F := \sqrt{\tr(A^*A)}$,\; $A \in \mathbb{C}^{m, n}$
\item $\diag(a_1, a_2, \dots , a_n) := [a_{ij}] \in \mathbb{C}^{n, n}$, where $a_{ij} =\delta_{ij}a_{i}$ for $i, j \in \{1,...,n\}$
\item $A^*$ - the conjugate transpose of $A$ 
\item
$
A \oplus B :=
    \begin{bmatrix}
A & 0 \\
0 & B
    \end{bmatrix}
$ (also for non-square matrices $A, B$)
\item $\sigma_{\min}(A) := \sqrt{\min\sigma(A^*A)}$
\item $\lambda_H(X) := \max\sigma\Bigl(\frac{X + X^*}{2}\Bigr)$
\item $A\leq B \iff B - A$ is positive semi-definite
\item $W(A) := \{x^*Ax : x \in \mathbb{C}^n \;\text{and}\; x^*x =1\},\quad A \in \mathbb{C}^{n, n}$
\item $W\bigl(P(\lambda)\bigr) := \big\{\lambda \in \mathbb{C} : x^*P(\lambda)x = 0 \;\;\text{for some}\;\; x \in \mathbb{C}^n \setminus \{0\}\big\},\quad P(\lambda) \in \mathbb{C}^{n, n}[\lambda]$
\item $\re\lambda = (\lambda + \overline{\lambda})/2$,\; $\im\lambda = (\lambda - \overline{\lambda})/(2\ii),\quad \lambda \in \mathbb{C}$
\item $\RE A := (A + A^*)/2$,\; $\IM A := (A - A^*)/(2\ii),\quad A \in \mathbb{C}^{n, n}$
\item $\Arg\lambda := \psi \in (-\pi; \pi]$ such that $\psi \in \arg\lambda$, for $\lambda \in \mathbb{C}\setminus\{0\}$,\; $\Arg 0 := 0$
\item $H_{\varphi} := \{\lambda \in \mathbb{C} : \im(\lambda e^{i\varphi}) > 0\}$, where $\varphi \in [0;2\pi)$
\item $H_0 = \{\lambda \in \mathbb{C} : \im\lambda > 0\}$
\item $H_{\varphi}^\kappa = \bigl\{(z_1, z_2, \dots , z_{\kappa}) \in \mathbb{C}^{\kappa} : \im(z_j e^{\ii\varphi}) > 0 \;\text{for all}\; j \in \{1, 2, \dots , \kappa\}\bigr\}$
\item $A(z^{(j)}) := \{w \in \mathbb{C} : ({z_1}, \dots , {z_{j-1}}, w, {z_{j+1}}, \dots , {z_{\kappa}}) \in A\}$ for $A \subseteq \mathbb{C}^{\kappa}$
\item $\mathcal{U}_n$ - the group of unimodular matrix polynomials with coefficients size $n \times n$
\item $\conv A$ - the smallest (with respect to inclusion) convex set containing $A$
\item $\mathbb{C}_{\kappa}[\lambda]$ - the vector space of complex polynomials with degree less or equal to $\kappa \in \mathbb{Z}_+$
\item 
$$
s_0(z_1, z_2, \dots, z_{\kappa}) := 1,\;\;\; s_j(z_1, z_2, \dots , z_{\kappa}) := \sum_{1 \leq i_1 < i_2 < \dots < i_j \leq \kappa} z_{i_1}z_{i_2} \dots z_{i_j}
$$
\item $T_\kappa : \mathbb{C}^{n, n}_{\kappa}[\lambda] \to \mathbb{C}^{n, n}[z_1, z_2, \dots , z_{\kappa}]$ 
$$
(T_{\kappa}P)(z_1, z_2, \dots , z_{\kappa}):= \sum_{j=0}^{\kappa} \binom{\kappa}{j}^{-1}s_j(z_1, z_2, \dots , z_{\kappa})A_j
$$
    \end{enumerate}
\section*{Definitions}

\begin{enumerate}[]
\item stable polynomial: Definition~\ref{stable1}
\item multi-affine and symmetric polynomial: Definition~\ref{ma_sym}
\item separately convex set: Definition~\ref{scx}
\item singular and regular matrix polynomial: Definition~\ref{rpoly}
\item Smith canonical form: Definition~\ref{Scf}
\item univariate matrix polynomial stable (hyperstable) with respect to $D$: Definition~\ref{defla} 
\item multivariate matrix polynomial stable (hyperstable) with respect to $D$: Definition~\ref{defzi}
\item separately convex set with respect to $j$-th variable: Definition~\ref{scxx}
\end{enumerate}

\chapter{Preliminaries}

\section{Polynomials -- conventions}

We refer to polynomials with complex coefficients as to scalar polynomials and to polynomials with matrix coefficients as to matrix polynomials. By $\mathbb{C}[\lambda]$ we denote the set of one variable complex scalar polynomials. Analogously, $\mathbb{C}[z_1, z_2, \dots , z_{\kappa}]$, where $\kappa \in \mathbb{Z}_+$, stands for the set of $\kappa$-variable complex scalar polynomials. Similar notation is used for matrix polynomials. Let $n \in \mathbb{Z}_+$ denote a size of a square matrix. By $\mathbb{C}^{n, n}[\lambda]$ we mean the set of one variable complex matrix polynomials with coefficients of size $n \times n$ and by $\mathbb{C}^{n, n}[z_1, z_2, \dots , z_{\kappa}]$ we mean the set of $\kappa$-variable complex matrix polynomials with coefficients of size $n \times n$. For scalar coefficients we use lower case letters and for matrix coefficients - upper case letters.

\begin{definition}\label{stable1}\rm
We say that a scalar polynomial $p(\lambda) \in \mathbb{C}[\lambda]$ is \emph{stable with respect to a nonempty set $D \subseteq \mathbb{C}$} if and only if it has no roots in $D$. More generally, we say that a scalar multi-variable polynomial $p(z_1, z_2, \dots , z_{\kappa}) \in \mathbb{C}[z_1, z_2, \dots , z_{\kappa}]$ is stable with respect to a nonempty set $D \subseteq \mathbb{C}^{\kappa}$ (or shortly $D$-stable) if and only if $p(z_1, z_2, \dots , z_{\kappa}) \neq 0$ for all $(z_1, z_2, \dots , z_{\kappa}) \in D$.
\end{definition}
For stability with respect to the open upper half-plane $H_0$, we come with interesting example providing a necessary and sufficient condition at the same time for a quadratic palindromic polynomial to be stable.
\begin{example}\rm 
Consider a complex quadratic palindromic polynomial $p(\lambda) = a\lambda^2 + b\lambda + a$ with $a \neq 0$. Let  $\lambda_1, \lambda_2$ denote the roots of the polynomial $p(\lambda)$. We have the following equivalence: the polynomial $p(\lambda)$ is stable with respect to the open upper half-plane $H_0$ if and only if $\lambda_1, \lambda_2 \in \mathbb{R}$. Only the forward implication requires an argument, hence, assume that the polynomial $p(\lambda)$ is stable with respect to $H_0$. By Vieta's formula for the product of the roots we have $\lambda_2 = 1 / \lambda_1$. Therefore, $\im\lambda_2 = -\im\lambda_1/|\lambda_1|$, which implies $\im\lambda_1\im\lambda_2 = - (\im\lambda_1)^2/|\lambda_1| \leq 0$. On the other hand,  stability implies that $\im\lambda_1\im\lambda_2 \geq 0$. Thus, $\im\lambda_1 = 0$ and $\im\lambda_2 = 0$.

Further, we have the second following equivalence: $\lambda_1, \lambda_2 \in \mathbb{R}$ if and only if there exists a real number $\mu$ such that $|\mu| \geq 2$ and $b = \mu a$. Indeed, assume first that $\lambda_1, \lambda_2 \in \mathbb{R}$ and that $b = \mu a$ with some complex number $\mu$. Then, by the Vieta's formula for the sum of roots, $\mu = -(\lambda_1 + \lambda_2) \in \Real$. Furthermore,  $|\mu| = |\lambda_1 + \lambda_2| = |\lambda_1 + 1/\lambda_1| \geq 2$ from the elementary inequality: $x^2+1 \geq 2x$ valid for all $x \in \mathbb{R}$.
For the converse implication note that 
	\begin{equation*}
p(\lambda) = a\lambda^2 + b\lambda + a = a\lambda^2 + \mu a\lambda + a = a(\lambda^2 + \mu\lambda + 1),
	\end{equation*}
and the discriminant of a polynomial in parentheses from the last equation equals $\mu^2 - 4 $, which is non-negative, since $\mu \in \mathbb{R}$ and $|\mu| \geq 2$. 

Summing up, we combine two above equivalences and obtain that the polynomial $p(\lambda)$ is stable with respect to the open upper half-plane $H_0$ if and only if there exists a real number $\mu$ such that $|\mu| \geq 2$ and $b = \mu a$.
	\end{example}
From other important conventions regarding polynomials, let us recall some definitions for scalar multivariate polynomials. 
    \begin{definition}\label{ma_sym}\rm
We say that a scalar multi-variate polynomial $p(z_1, \dots, z_\kappa)$  is \emph{multi-affine} if and only if its degree with respect to each variable $z_j$ ($j=1,\dots,\kappa$) is less or equal to one. Further, we say that $ p(z_1, \dots, z_\kappa)$ is \emph{symmetric} if and only if any permutation of the variables $z_1, \dots, z_\kappa$ leaves $ p(z_1, \dots, z_\kappa)$ intact.
    \end{definition}


\section{Gauss-Lucas Theorem}

Recall the following classical result.

    \begin{theorem}\label{GLscalar}
{\rm(Gauss-Lucas)} Let $p \in \mathbb{C}[\lambda]$ be a non-constant complex polynomial. Then the following inclusion holds:
        \begin{equation*}
\{\lambda \in \mathbb{C} : p'(\lambda) = 0\}  \subset \conv\{\lambda \in \mathbb{C} : p(\lambda) = 0\},
        \end{equation*}
where the symbol $\conv$ denotes the convex hull of a set in the complex plane $\mathbb{C}$.
    \end{theorem}

There exists an extension of the above result onto several variables, see \cite{kanter}. To review it we first need a notion of the separately convex set. For $z = (z_1, \dots , z_{\kappa})$ and $j \in \{1, 2, \dots , \kappa\}$ let $z^{(j)}$ denotes ordered $(\kappa - 1)$-tuple of coordinates obtained from $z$ by omitting the coordinate $z_j$. For a subset $A \subseteq \mathbb{C}^{\kappa}$ and a fixed $z 
\in \mathbb{C}^{\kappa}$ let us define the section of $A$ determined by $z^{(j)}$ as follows
    \begin{equation}
    A(z^{(j)}) := \{w \in \mathbb{C} : ({z_1}, \dots , {z_{j-1}}, w, {z_{j+1}}, \dots , {z_{\kappa}}) \in A\}.
    \end{equation}

\begin{definition}\label{scx}\rm
We call a set $A \subseteq \mathbb{C}^{\kappa}$ \emph{separately convex} in $\mathbb{C}^{\kappa}$ if the set $A(z^{(j)})$ is a convex subset of the complex plane  for all $j \in \{1, 2, \dots , \kappa\}$ and all $z \in \mathbb{C}^{\kappa}$.
Because an intersection of two separately convex sets in $\mathbb{C}^{\kappa}$ is again a separately convex set in $\mathbb{C}^{\kappa}$, we can define the {\em separately convex hull} of a set $A \subseteq \mathbb{C}^{\kappa}$ as the smallest  separately convex set  containing $A$; we denoted it by $\rm{conv_{\kappa}}(A)$.
\end{definition}

Note that a  convex set is separately convex and hence $\rm{conv}_{\kappa}(A)\subseteq \conv(A)$. With these preparations we have the following (\cite{kanter}).

    \begin{theorem}\label{mgl}
Let $j\in\{1, 2, \dots, \kappa\}$. If $p \in \mathbb{C}[z_1, z_2, \dots , z_{\kappa}]$ is a multivariate polynomial with nonzero $j$-th partial derivative $\partial p/\partial z_j$, then
        \begin{equation}
\Bigl(\frac{\partial p}{\partial z_j}\Bigr)^{-1}(0) \subseteq \conv_{\kappa}\bigl(p^{-1}(0)\bigr)\text{.}
        \end{equation}
    \end{theorem}
In Chapter $3$ we will generalise the above theorems to the context of matrix polynomials. Let us also mention, that the refinements of the Gauss-Lucas Theorem are still a current topic, see, e.g.,  the recent paper \cite{steinerberge}. However, it seems very technical to generalise the results of \cite{steinerberge} onto matrix polynomials.

\section{The Classical Sz\'asz inequality}\label{SzaszClass}

At the end of this Section we present the strong Sz\'asz inequality which will be improved later, in Section \ref{sGL}. The classical inequality discovered by O. Sz\'asz bounds a stable polynomial $p(\lambda) = a_d \lambda^d + a_{d-1} \lambda^{d-1} + \dots + a_1 \lambda + 1 \in \mathbb{C}[\lambda]$ with $p(0) = 1$ in terms of its first few coefficients:
        \begin{equation*}
|p(\lambda)| \leq \exp\bigl(|a_1||z\lambda| + 3(|a_1|^2 + |a_2|)|\lambda|^2\bigr)\text{,}\quad \lambda \in \mathbb{C}\text{.}
        \end{equation*}
Later on the inequality was improved by de Branges \cite[Lemma 5]{deB61} to
        \begin{equation}\label{dB}
|p(\lambda)| \leq \exp\bigl(\re (a_1\lambda)  + \frac12 (|a_1|^2 -2 \re(a_2))|\lambda|^2\bigr),\quad \lambda\in\Comp
        \end{equation}
and Knese showed in \cite[Theorem 1.3]{Kne19} sharpness of the bound. The proof of \eqref{dB} is based on the inequality $\log|1+w|\leq \re w+1/2|w|^2$ applied several times.
    \begin{theorem}
 Let $p(\lambda) = a_d \lambda^d + a_{d-1} \lambda^{d-1} + \dots + a_1 \lambda + 1 \in \mathbb{C}[\lambda]$ be a stable polynomial, where $p(0) = 1$. Then, the following inequality holds:
        \begin{equation*}
|p(\lambda)| \leq \exp\bigl(\re(a_1 \lambda) + \frac{1}{2}(|a_1|^2 - 2\re a_2)|\lambda|^2\bigr)\text{.}
        \end{equation*}
Moreover, the inequality is sharp on the imaginary axis for stable polynomials $p \in \mathbb{R}[\lambda]$ with real coefficients and $p(0) = 1$. Namely, for $c_1, c_2 \in \mathbb{R}$ such that $\gamma := \frac{1}{2}(c_1^2 - 2c_2) > 0$ there exist stable polynomials $p_n(\lambda) = 1 + c_1 \lambda + c_2 \lambda^2 + \dots \in \mathbb{R}[\lambda]$ satisfying the following condition:
        \begin{equation*}
\lim_{n \to \infty}|p_n(iy)| = \exp(\gamma y^2) \;\;\text{for all}\;\; y \in \mathbb{R}\text{.}
        \end{equation*}
    \end{theorem}

In Chapter ~\ref{SzaszM} we present several generalisations of the above inequality to matrix polynomials.

\section{Operators preserving stability of scalar multivariate polynomials}
Let us begin with some background on stability. Two possible generalizations of this notion in the matrix polynomials theory, see Definition \ref{defla} and Definition \ref{defzi}, will be studied further.
\begin{example}\rm
Below there are a few simple examples of $H_0^n$-stable polynomials:
    \begin{enumerate}[(i)]
\item the monomials $(z_1, z_2, \dots, z_{n}) \mapsto \prod_{j=1}^n z_j$,
\item linear polynomials $(z_1, z_2, \dots, z_{n}) \mapsto \sum_{j=1}^n a_jz_j$ with positive coefficients $a_j \in \mathbb{R}_+$,
\item a polynomial $p(z_1, z_2) = c - z_1z_2$, where $c$ is a positive constant, is $H_0^2$-stable (if $z_1, z_2 \in H_0$, then $z_1z_2 \not\in \mathbb{R}_+ \cup \{0\}$).
    \end{enumerate}
Note that a polynomial $q(z_1, z_2) = c + z_1z_2$ is not $H_0^2$-stable, since it has a zero $(z_1, z_2) = (c\ii, \ii)$ belonging to the Cartesian square $H_0^2$ of the open upper half-plane.
\end{example}
We have the following criterion of $H_0^n$-stability, which allows us to verify $H_0^n$-stability of a multivariate polynomial through checking $H_0$-stability of its univariate restriction on some complex lines, cf. Lemma 16.3 in \cite{lecjv}.
    \begin{proposition}
Let $t$ be a complex variable. A polynomial $p \in \mathbb{C}[z_1, z_2, \dots , z_n]$ is $H_0^n$-stable if and only if for every $e \in \mathbb{R}_+^n$ and every $x \in \mathbb{R}^n$ the univariate restriction $t \mapsto p(te + x)$ is a $H_0$-stable polynomial.
    \end{proposition}
There is one fundamental (and highly non-trivial) example of a $H_0^n$-stable polynomial which is, in fact, a determinant of a multi-variable matrix polynomial and which we present here.
\begin{example}\rm
Let $A \in \mathbb{C}^{m, m}$ be a Hermitian matrix and let $B_1, B_2, \dots , B_n \in \mathbb{C}^{m, m}$ be any positive semi-definite matrices. Then the determinantal polynomial
    \begin{equation*}
p(z_1, z_2, \dots , z_n) = \det\Bigl(A + \sum_{j=1}^n z_jB_j\Bigr)
    \end{equation*}
is a $H_0^n$-stable polynomial with real coefficients.
\end{example}
Another remarkable theorem says that in the two-variable case determinantal polynomials are only $H_0^2$-stable ones with real coefficient. However, it is known to be false for more than two variables. Now, let us recall the famous Helton-Vinnikov Theorem (cf. \cite{vinnikov}).
    \begin{theorem}
{\rm(Helton-Vinnikov)} Let $p \in \mathbb{R}[x, y]$ be a bivariate polynomial of total degree $d$. If the polynomial $p$ is $H_0^2$-stable, then there exist real symmetric $d \times d$ positive semi-definite matrices $A, B$ and a real symmetric $d \times d$ matrix $C$ such that
        \begin{equation*}
p(x, y) = \pm\det(xA + yB + C)\text{.}
        \end{equation*}
    \end{theorem}
In the above case we say that the polynomial $p$ has a \emph{determinantal representation}. It turns out that in a general case, all $H_0^2$-stable polynomials of two complex variables have determinantal representations. We acknowledge the following results proven by G. Knese (see \cite{Kne19}, Theorems 3.2 and 4.1)
    \begin{theorem}
Let $p \in \mathbb{C}[z_1, z_2]$ be a bivariate polynomial of degrees $n$ and $m$ with respect to the variables $z_1$ and $z_2$, respectively. If the polynomial $p$ has no zeros in bidisk $\mathbb{D}^2$, then there exist a constant $c$ and an contractive matrix $D$ such that
        \begin{equation}\label{bi}
p(z_1, z_2) = c\det\bigl(I - D\Delta(z_1, z_2)\bigr), 
        \end{equation}
where $\Delta(z_1, z_2) = z_1(I_n \oplus 0_m) + z_2(0_n \oplus I_m)$.
    \end{theorem}
Knese converted the above bidisk formula \eqref{bi} to a formula related to an open upper half-plane $H_0$, see the theorem below. By the total degree of a multi-variable polynomial we mean the maximum degree of all its monomials, i.e., the degree with respect to the system of all its variables.
    \begin{theorem}
If $p \in \mathbb{C}[z_1, z_2]$ is a stable with respect to the product $H_0^2$ polynomial of total degree $d$, then there exist $d \times d$ matrices $A, B_1, B_2$ and a constant $c \in \mathbb{C}$ such that
        \begin{enumerate}[\rm (i)]
\item $\IM A := (A - A^*)/(2i) \geq 0$,
\item $B_1, B_2 \geq 0$,
\item $B_1 + B_2 = I$,
\item $p(z_1, z_2) = c\det(A + z_1 B_1 + z_2 B_2)$.
        \end{enumerate}
    \end{theorem}
Now we claim that $H_0^n$-stability is closed under taking limits. This result is usually proved by induction on $n$, while the case $n = 1$ is a classic complex analysis result, see Lemma 15.3 and Lemma 16.4 in \cite{lecjv} .
    \begin{theorem}\label{limiting}
Let $p_1, p_2, \dots , p_k, \dots \in \mathbb{C}[z_1, z_2, \dots , z_n]$ be a sequence of $H_0^n$-stable polynomials of bounded total degree and assume that $p_k \to p$ coefficient-wise as $k \to \infty$. Then either the limit polynomial $p$  is $H_0^n$- stable or $p \equiv 0$.
    \end{theorem}
    
Let us list some basic operations on multi-variable polynomials which are stable with respect to $n$-th Cartesian power of an open half-plane. Below $H_{\varphi} := \{\lambda \in \mathbb{C} : \im(\lambda e^{i\varphi}) > 0\}$ is an open half-plane with the boundary containing the origin and $H_{\varphi}^n$ is its $n$-th Cartesian power.
    \begin{proposition}\label{op1}
Let $\varphi \in [0; 2\pi)$ be fixed. Then the following linear transformations acting on the complex polynomial space $\mathbb{C}[z_1, z_2, \dots , z_n]$ map every $H_{\varphi}^n$-stable polynomial to another $H_{\varphi}^n$-stable polynomial or to zero:
        \begin{enumerate}[\rm (i)]
\item{Permutation:}
$$
p(z_1, z_2, \dots , z_n) \mapsto p(z_{\sigma(1)}, z_{\sigma(2)}, \dots , z_{\sigma(n)})
$$
for every permutation $\sigma \in S_n$. \\
\item{Scaling:}
$$
p(z_1, z_2, \dots , z_n) \mapsto p(z_1, z_2, \dots , z_{j-1}, az_j, z_{j+1}, \dots , z_n)
$$
for any $a \in \mathbb{R}_+$. \\
\item{Diagonalization:}
$$
p(z_1, z_2, \dots , z_n) \mapsto p(z_j, z_j, \dots , z_j, z_{j+1}, \dots , z_n) \in \mathbb{C}[z_j, z_{j+1}, \dots , z_n]
$$
for any $j \in \{1, 2, \dots , n\}$. \\
\item{Inversion (with rotation):}
$$
p(z_1, z_2, \dots , z_n) \mapsto z_j^dp(z_1, z_2, \dots , z_{j-1}, -e^{-2\ii\varphi}/z_j, z_{j+1}, \dots , z_n), 
$$
where the power $d := \deg_jp$ is the degree of a polynomial $p$ with respect to the variable $z_j$. \\
\item{Specialization:}
$$
p(z) \mapsto p(z_1, z_2, \dots , z_{j-1}, a, z_{j+1}, \dots , z_n) \in \mathbb{C}[z_1, z_2, \dots , z_{j-1}, z_{j+1}, \dots , z_n]
$$
for any $a \in \overline{H}_{\varphi}$, where $z = (z_1, z_2, \dots , z_n)$. \\
\item{Differentiation:}
$$
p(z_1, z_2, \dots , z_n) \mapsto \frac{\partial p}{\partial z_j}(z_1, z_2, \dots , z_n) \in \mathbb{C}[z_1, z_2, \dots , z_n],\; j \in \{1, 2, \dots , n\}.
$$
        \end{enumerate}
    \end{proposition}

As the full proof is hard to find in the literature (cf. Lemma 1.7 in \cite{BorB09}), we present it here for completeness.  

    \begin{proof}
For a proof of (i) it is enough to see that a permutation of variables $z_1, z_2, \dots , z_n$ is a bijection between a region $H_{\varphi}^n$ of the $n$-dimensional complex space $\mathbb{C}^n$ and itself. Similarly, (ii) follows from the fact that scaling by a constant $a \in \mathbb{R}_+$ is a bijection of the half-plane $H_{\varphi}$ on itself. The property (iii) is an immediate consequence of the definition of the stability notion. To prove the property (iv) note that if $\lambda \in H_{\varphi}$, then $f(\lambda) = -e^{-2\ii\varphi}/\lambda \in H_{\varphi}$ as well, as the function $f$ is a superposition of the inversion $\lambda \mapsto -1/\lambda$ and the rotation $\lambda \mapsto \lambda e^{-2\ii\varphi}$ around the origin by the angle $-2\varphi$. Regarding property (v), the case $a \in H_{\varphi}$ follows straight from the definition of the stability notion. For the case $a \in \partial H_{\varphi}$, we should use the limiting argument of Theorem \ref{limiting}. The last property (vi) becomes an easy corollary from Gauss-Lucas Theorem (see Theorem ~\ref{GLscalar}), when we notice that the complement $\mathbb{C} \setminus H_{\varphi}$ of the half-plane $H_{\varphi}$ is a convex set.
    \end{proof}
A full characterisation of operators preserving stability is provided in the seminal paper \cite{BorB09}.

\section{Polarisation operators}

We present here the method of polarisation for scalar polynomials, which will be extended onto matrix polynomials in Section $\ref{sPol}$. Let $\kappa \in \mathbb{Z}_+$, consider the following standard polarization operators $T_{\kappa} : \mathbb{C}_{\kappa}[\lambda] \longrightarrow \mathbb{C}[z_1, z_2, \dots, z_{\kappa}]$ given by the formulas:
    \begin{equation*}
(T_{\kappa}p) (z_1, z_2, \dots , z_{\kappa}):= \sum_{j=0}^{\kappa} \binom{\kappa}{j}^{-1} a_j s_j(z_1, z_2, \dots , z_{\kappa}), 
	\end{equation*}
where the polynomial $p(\lambda) = a_d\lambda^d + a_{d-1}\lambda^{d-1} + \dots + a_0 \in \mathbb{C}_{\kappa}[\lambda]$ has degree $d \leq \kappa$ and the symbol $s_j$ denotes the $j$-th elementary symmetric polynomial
	\begin{equation*}
s_0(z_1, z_2, \dots, z_{\kappa}) := 1,\;\;\; s_j(z_1, z_2, \dots , z_{\kappa}) := \sum_{1 \leq i_1 < i_2 < \dots < i_j \leq \kappa} z_{i_1}z_{i_2} \dots z_{i_j}.
	\end{equation*}
	
It is well known, that these operators preserve stability of scalar polynomials, see Proposition 3.4 in \cite{BorB09} and Theorem \ref{Tkappa2} in case of $1 \times 1$ matrices. We will prove the statement in an extended form for matrix polynomials. The basic tool is the famous result by Grace, Walsh and Szeg\"o (cf. \cite{grace1902, walsh1922, szego1922}), the so called coincidence theorem. For its formulation we need the notion of a circular domain. Namely, an open or closed subset $D$ of the complex plane $\mathbb{C}$ bounded by either a circle or a straight line is called a \emph{circular domain}.
    \begin{theorem}\label{GWS}
\rm{(Grace-Walsh-Szeg\"o)} Let $p \in \mathbb{C}[z_1, z_2, \dots , z_{\kappa}]$ be a symmetric multi-affine polynomial and let $\zeta_1, \zeta_2, \dots , \zeta_{\kappa}$ be points in some circular domain $D$. If the total degree of polynomial $(z_1, z_2, \dots , z_{\kappa}) \mapsto p(z_1, z_2, \dots , z_{\kappa})$ equals $\kappa$ or the domain $D$ is a convex set, then there exists a point $\zeta_0 \in D$ such that
		\begin{equation*}
p(\zeta_1, \zeta_2, \dots , \zeta_{\kappa}) = p(\zeta_0, \zeta_0, \dots , \zeta_0)\text{.}
		\end{equation*}
	\end{theorem}

 \section{Univariate matrix polynomials}
In this Section Smith canonical form of matrix polynomial will be revisited. We need first definition of some basic notions.
    \begin{definition}\rm\label{rpoly}
We call a matrix polynomial $P(\lambda) = \lambda^d A_d + \lambda^{d-1} A_{d-1} + \dots + A_0 \in \mathbb{C}^{n, n}[\lambda]$, where $A_d \neq 0$, a \emph{regular polynomial} if and only if the function $\lambda \mapsto \det P(\lambda)$ is a nonzero scalar polynomial. Otherwise, the polynomial $P(\lambda)$ is called a \emph{singlular polynomial}. In the regular case, we call a vector $x \in \mathbb{C}^n \setminus \{0\}$ an \emph{eigenvector} of the polynomial $P(\lambda)$ if and only if there exists $\lambda_0 \in \mathbb{C}$ such that
		\begin{equation*}
P(\lambda_0)x = 0\text{.}
		\end{equation*}
Then such complex number $\lambda_0$ is called an \emph{eigenvalue} of the polynomial $P(\lambda)$ corresponding to the eigenvector $x$. Note that for a regular matrix polynomial $P(\lambda)$ being an eigenvalue is equivalent to being a root of the scalar polynomial $\det P(\lambda)$. In particular, when $\det P(\lambda) = \text{const.} \neq 0$ we refer to the regular polynomial $P(\lambda)$ as to an \emph{unimodular polynomial}.

    \end{definition}

\begin{example}\rm
The polynomials 
$$    \begin{bmatrix}
1 & 0 \\
0 &\lambda
    \end{bmatrix},\ 
    \begin{bmatrix}
1 & \lambda^2 \\ 
0 & 1
    \end{bmatrix}
$$ 
are regular, the latter being additionally unimodular, while the polynomials 
  $$  \begin{bmatrix}
0 & 0 \\
0 & 0
    \end{bmatrix},\ 
    \begin{bmatrix}
\lambda^2 & \lambda \\
\lambda & 1 
    \end{bmatrix}
$$
are singular.
\end{example}

A general decomposition tool, applicable for  any matrix polynomials is delivered by the following Smith canonical form, see \cite[chapter VI, volume I]{Gan59}, which is a representation of a matrix polynomial under unimodular equivalence.
     \begin{theorem}\label{Scf}
(Smith canonical form) Let $P(\lambda) \in \mathbb{C}^{m, n}[\lambda]$ be a matrix polynomial of a degree $d \in \mathbb{Z}_+ \cup \{0\}$. Then there exist unimodular polynomials $U(\lambda) \in \mathbb{C}^{m, m}[\lambda]$ and $V(\lambda) \in \mathbb{C}^{n, n}[\lambda]$ such that
        \begin{equation*}
U(\lambda)P(\lambda)V(\lambda) = 
    \left[\begin{array}{@{}c|c@{}}
&\\
\underset{}{\diag{\bigl(s_1(\lambda), s_2(\lambda), \dots , s_r(\lambda)\bigr)}} & 0_{r \times (n-r)} \\
\hline
\underset{}{0_{(m-r) \times r}} & 0_{(m-r) \times (n-r)}
    \end{array}\right]
=: S(\lambda)\text{,}
        \end{equation*}
where $0 \leq r \leq \min\{m, n\}$ and the entries $s_1(\lambda), s_2(\lambda), \dots , s_r(\lambda) \in \mathbb{C}[\lambda]$ are uniquely determined monic scalar polynomials such that $s_{j+1}(\lambda)|s_{j}(\lambda)$ for $j \in \{1, 2, \dots , r-1\}$. 
    \end{theorem}

\begin{definition}\label{Scf}\rm
The quasi-diagonal polynomial matrix $S(\lambda) \in \mathbb{C}^{m, n}[\lambda]$ is called the \emph{Smith canonical form} of the polynomial $P(\lambda)$ and the monic polynomials $s_1(\lambda), s_2(\lambda), \dots , s_r(\lambda)$ are called the \emph{invariant factors} of $P(\lambda)$. Their roots are called \emph{finite eigenvalues} of $P(\lambda)$ and in the  case of a regular matrix polynomial  they coincide with eigenvalues as defined in Definition ~\ref{rpoly}. We say that $P(\lambda)$ has an eigenvalue at infinity $\infty$ if and only if its reversal polynomial $R(\lambda) := \lambda^d P(1/\lambda)$ has an eigenvalue at zero. Moreover, the Smith canonical form $S(\lambda)$ of the polynomial $P(\lambda)$ is unique.
\end{definition}

    \begin{theorem}
(characterization of invariant factors)  Let $P(\lambda) \in \mathbb{C}^{m, n}[\lambda]$ be a matrix polynomial and let $S(\lambda)$ be its Smith canonical form. Set $p_0(\lambda) \equiv 1$. For $1 \leq j \leq \min\{m, n\}$ set $p_j(\lambda) \equiv 0$ if all minors of the polynomial $P(\lambda)$ of order $j$ are equal to zero and take $p_j(\lambda)$ as the greatest common divisor of all minors of $P(\lambda)$ otherwise. Then the number $r$ of invariant factors of $P(\lambda)$ is the largest integer such that $p_r(\lambda) \not\equiv 0$, i.e. $r = \rank P(\lambda)$ and the invariant factors $s_1(\lambda), s_2(\lambda), \dots , s_r(\lambda)$ of the polynomial $P(\lambda)$ are ratios of consecutive greatest common divisors $p_r(\lambda), p_{r-1}(\lambda), \dots , p_0(\lambda)$:
        \begin{equation*}
s_j(\lambda) = \frac{p_{r-j+1}(\lambda)}{p_{r-j}(\lambda)} \;\;\;\text{for}\;\;\; j \in \{1, 2, \dots , r\}\text{.}
        \end{equation*}
    \end{theorem}

\section{Matrix pencils}

By a matrix pencil we mean a matrix polynomial of degree $1$. Besides the Smith canonical form there exists a more detailed description, which is the Kronecker canonical form. We will use the following symbols:
   \begin{equation*}
\mathcal{J}_n^{\mu}(\lambda) := 
        \begin{bmatrix}
\;\;\, \lambda - \mu & 1 &  &  &\\
  & \lambda - \mu & \;\ddots &  &\\
  &  & \;\ddots & \; 1 &\\
  &  &  & \!\!\lambda - \mu\!\!\!\!
        \end{bmatrix} \in \mathbb{C}^{n, n}[\lambda]\text{,}
    \end{equation*}

 \begin{equation*}
\mathcal{N}_n(\lambda) :=
        \begin{bmatrix}
\;\;\, 1 & \lambda &  &  &\\
  & 1 & \;\ddots &  &\\
  &  & \;\ddots & \; \lambda &\\
  &  &  & \!\! 1 \!\!\!\!
        \end{bmatrix} \in \mathbb{C}^{n, n}[\lambda]\text{,}
    \end{equation*}
    
    \begin{equation*}
\mathcal{B}_n(\lambda) := \lambda
        \begin{bmatrix}
\; 1 & 0 &  &  &\\
  & \ddots & \ddots &  &\\
  &  & 1 & 0\!\!\!\!\!\!
        \end{bmatrix} 
- 
        \begin{bmatrix}
\; 0 & 1 &  &  &\\
  & \ddots & \ddots &  &\\
  &  & 0 & 1\!\!\!\!\!\!
        \end{bmatrix} \in \mathbb{C}^{n, n+1}[\lambda]\text{.}
    \end{equation*}
  
Next theorem can be found in \cite[chapter XII, volume II]{Gan59}.
    \begin{theorem}
(Kronecker canonical form) Let $E, A \in \mathbb{C}^{m, n}$ be any matrices. Then there exist non-singular matrices $S \in \mathbb{C}^{m, m}$ and $T \in \mathbb{C}^{n, n}$ such that
        \begin{eqnarray*}
S(\lambda E - A)T &=& 
\mathcal{J}^{\lambda_1}_{\tau_1}(\lambda) \oplus \dots \oplus \mathcal{J}^{\lambda_r}_{\tau_r}(\lambda) \\
&\oplus& \mathcal{N}_{\omega_1}(\lambda) \oplus \dots \oplus \mathcal{N}_{\omega_s}(\lambda) \\
&\oplus& \mathcal{B}_{\rho_1}(\lambda) \oplus \dots \oplus \mathcal{B}_{\rho_p}(\lambda) \\
&\oplus& \mathcal{B}_{\sigma_1}^{\top}(\lambda) \oplus \dots \oplus \mathcal{B}_{\sigma_q}^{\top}(\lambda) \\
&\oplus& 0_{M \times N}
        \end{eqnarray*}
where $p, q, r, s \in \mathbb{Z}_+ \cup \{0\}$, $\rho_1, \rho_2, \dots , \rho_p, \sigma_1, \sigma_2, \dots , \sigma_q, \tau_1, \tau_2, \dots , \tau_r, \omega_1, \omega_2, \dots , \omega_s \in \mathbb{Z}_+$, $M, N \in \mathbb{Z}_+$ or $M = N = 0$ and $\lambda_1, \lambda_2, \dots , \lambda_r \in \mathbb{C}$.

    \end{theorem}
    \begin{remark}\rm
The above block-diagonal matrix is unique up to a permutation of the blocks and is called the \emph{Kronecker canonical form} of a matrix pencil $P(\lambda) = \lambda E - A$. The values $\lambda_1, \lambda_2, \dots , \lambda_r$ are the finite eigenvalues of the pencil $\lambda E - A$ as defined in Definition ~\ref{Scf}. Furthermore, $\lambda_0 = \infty$ is an infinite eigenvalue of $\lambda E - A$ if and only if zero is an eigenvalue of $\lambda A - E$. The last statement holds exactly when $s \neq 0$. For $m = n$ ($p = q$ in this case), infinity $\lambda_0 = \infty$ is an eigenvalue of a matrix pencil $\lambda E - A$ if and only if the leading coefficient $E$ is a non-invertible matrix. The sum of all sizes of blocks associated with a fixed eigenvalue $\lambda_0 \in \mathbb{C} \cup \{\infty\}$ of a matrix pencil $\lambda E - A$ is called the \emph{algebraic multiplicity} of $\lambda_0$.
A matrix pencil $P(\lambda) = \lambda E - A \in \mathbb{C}^{m, n}[\lambda]$ is called \emph{regular} if and only if $m = n$ and $\det(\lambda E - A) \neq 0$ for some $\lambda \in \mathbb{C}$; otherwise it is called \emph{singular}. A matrix pencil $\lambda E - A$ is regular if and only if it has no bidiagonal blocks $\mathcal{B}_{\rho_{j}}$ and $\mathcal{B}_{\sigma_{k}}^{\top}$ in its Kronecker canonical form. 
    \end{remark}
Square matrix pencils allow, in general, more factorizations than matrix polynomials. For instance, each of such pencils is unitarily equivalent to an upper-triangular one, what can be seen from the generalised Schur canonical form below, see \cite[Chapter 7]{gol2013}.
    \begin{theorem}\label{Schur}
(generalised Schur canonical form) Let $A, B \in \mathbb{C}^{n, n}$ be any square matrices. Then there exist unitary matrices $Q, Z \in \mathbb{C}^{n, n}$ such that there is a strict equivalence between a matrix pencil $P(\lambda) = \lambda A + B$ and some upper-triangular matrix pencil, namely
        \begin{equation*}
Q^*(\lambda A + B)Z = \lambda S + T,
        \end{equation*}
where matrices $S, T \in \mathbb{C}^{n, n}$ are upper-triangular.
    \end{theorem}
    \begin{remark}\rm
 The upper-triangular matrix pencil $ \lambda S + T$ is called the \emph{generalised Schur canonical form} of the pencil $P(\lambda)$. If $s_{jj} = t_{jj} = 0$ for some $j \in \{1, 2, \dots , n\}$, then $P(\lambda)$ is singular. Otherwise, it is regular and the set of its eigenvalues (including infinity $\infty$) is precisely the set consisting of all quotients $s_{jj}/t_{jj}$ for $j \in \{1, 2, \dots, n\}$. The above canonical form for square matrix pencils is a generalisation of the standard Schur canonical form for square matrices, which allows us to find an upper-triangular matrix unitarily similar to a given one. Indeed, for a matrix $A \in \mathbb{C}^{n, n}$ there exist an upper-triangular matrix $S \in \mathbb{C}^{n, n}$ and a unitary matrix $Q \in \mathbb{C}^{n, n}$ such that $Q^{-1}AQ = S \;(Q^{-1} = Q^*)$. Moreover, diagonal entries of the matrix $S$ are the eigenvalues of the matrix $A$.
    \end{remark}

\section{Numerical range of matrices and matrix polynomials}

Below we recall the definition and some basic properties of the numerical range of a matrix.

    \begin{definition}\rm
Let $A \in \mathbb{C}^{n, n}$ be a square matrix. A subset
        \begin{equation*}
W(A) := \{x^*Ax : x \in \mathbb{C}^n \;\text{and}\; x^*x =1\}
        \end{equation*}
of the complex plane $\mathbb{C}$ is called the \emph{numerical range} of the matrix $A$.
    \end{definition}

Next, we list some properties of the numerical range $W(A)$ according to \cite{HorJ91} (paragraph 1.2).

    \begin{proposition}
Let $A, B \in \mathbb{C}^{n, n}$ be arbitrary square matrices. The following properties hold:
        \begin{enumerate}[\rm (i)]
\item $W(A)$ is a compact subset of $\mathbb{C}$;
\item $W(A)$ is a convex subset of $\mathbb{C}$;
\item $W(A + \alpha I) = W(A) + \alpha$ for all $\alpha \in \mathbb{C}$;
\item $W(\alpha A) = \alpha W(A)$ for all $\alpha \in \mathbb{C}$;
\item $W(\RE A) = \re W(A)$, where $\re W(A) = \{\re z : z \in W(A)\}$;
\item $\sigma(A) \subseteq W(A)$, where the symbol $\sigma(A)$ denotes the spectrum of the matrix $A$;
\item $W(A + B) \subseteq W(A) + W(B)$;
\item $W(U^*AU) = W(A)$ for any unitary matrix $U \in \mathbb{C}^{n, n}$;
\item $W(A) = \conv\bigl(\sigma(A)\bigr)$ for any normal matrix $A \in \mathbb{C}^{n, n}$.
        \end{enumerate}
    \end{proposition}
    
Another important fact, which we use in the proof of Sz\'asz-type inequality for matrix polynomials, see Theorem ~\ref{Szasz} is the inequality between the numerical radius and matrix two-norm. It can be found, e.g., as Problem 23 (g) in \cite{HorJ85}.

    \begin{theorem}
Let $A \in \mathbb{C}^{n, n}$ be an arbitrary square matrix. The following inequalities hold
        \begin{equation*}
\frac{1}{2}\norm{A}_2 \leq w(A) \leq \norm{A}_2\text{,}
        \end{equation*}
where $w(A) := \max\{|z| : z \in W(A)\}$ is so called numerical radius of the matrix $A$.
    \end{theorem}

Besides a notion of the set of the eigenvalues of a matrix polynomial $P(\lambda) \in \mathbb{C}^{n, n}[\lambda]$, we need a notion of a wider set containing all eigenvalues, i.e. a notion of the numerical range of the polynomial $P(\lambda)$. It was introduced in \cite{LiR94}, see also \cite{Psa03,Psa00}.
    \begin{definition}\label{numran}
\rm Let $P(\lambda) \in \mathbb{C}^{n, n}[\lambda]$ be a square matrix polynomial. A subset
        \begin{equation}
W\bigl(P(\lambda)\bigr) := \big\{\lambda \in \mathbb{C} : x^*P(\lambda)x = 0 \;\;\text{for some}\;\; x \in \mathbb{C}^n \setminus \{0\}\big\}
        \end{equation}
of the complex plane $\mathbb{C}$ is called the \emph{numerical range} of the matrix polynomial $P(\lambda)$.
    \end{definition}
\begin{remark}\rm
In the case, when $P(\lambda) = \lambda I_n - A$ is a matrix pencil, where $E = I_n$, the set $W(P)$ reduces to the numerical range of a matrix $A \in \mathbb{C}^{n, n}$. This observation justify the statement that the numerical range of a matrix polynomial is a generalization of the numerical range of a matrix.
\end{remark}
Below we list some fundamental properties of the numerical range of a matrix polynomial, see \cite{LiR94}.
    \begin{proposition}\label{nrp}
Let $P(\lambda) = \lambda^d A_d + \lambda^{d-1} A_{d-1} + \dots + A_0 \in \mathbb{C}^{n, n}[\lambda]$ be a matrix polynomial, where $A_d \neq 0$. Then, the following conditions hold:
        \begin{enumerate}[\rm (a)]
\item\label{zeros} the numerical range $W(P)$ contains all zeros of the function $\lambda \mapsto \det P(\lambda)$,
\item\label{closed} the set $W(P)$ is a closed subset of the complex plane $\mathbb{C}$,
\item if $Q(\lambda) := P(\lambda + \alpha)$, then $W(Q) = W(P) - \alpha$ for all $\alpha \in \mathbb{C}$,
\item for the reversal polynomial $Q(\lambda) := \lambda^d A_0 + \lambda^{d-1} A_1 + \dots + A_d \in \mathbb{C}^{n, n}[\lambda]$ of the polynomial $P(\lambda)$ we have $W(Q)\setminus\{0\} = \{\mu^{-1} \in \mathbb{C} : \mu \in W(P), \;\mu \neq 0\}$,
\item for each $n \times m$ matrix $S$ with $\rank{S} = m \leq n$ we have $W(S^*PS) \subseteq W(P)$.
        \end{enumerate}
    \end{proposition}
As known, the classical numerical range $W(A)$ of a matrix $A \in \mathbb{C}^{n, n}$ is a compact convex subset of the complex plane $\mathbb{C}$. The condition \eqref{closed} of \ref{nrp} says that the set $W(P)$ is always closed. However, the numerical range $W(P)$ of a matrix polynomial $P(\lambda) \in \mathbb{C}^{n, n}[\lambda]$ need not be connected or even bounded, what the following example shows, cf. \cite[example 2.1]{LiR94}.
    \begin{example}\rm
Consider $P(\lambda) = \lambda
        \begin{bmatrix}
1 & 0\\
0 & \!\!\!\! -1\\
        \end{bmatrix} 
- 
        \begin{bmatrix}
1 & 0\\
0 & 1\\
        \end{bmatrix} = 
        \begin{bmatrix}
\lambda - 1 & 0 \\
0 & -\lambda - 1
        \end{bmatrix}$. Then, assuming without loss of generality  that $x^*x=1$,
$$
x^*P(\lambda)x = |x_1|^2(\lambda - 1) + |x_2|^2(-\lambda - 1) = (2|x_1|^2 - 1)\lambda - 1.
$$
Hence,
        \begin{equation*}
W(P) = \{re^{\ii\theta} \in \mathbb{C} : r \geq 1, \;\theta = 0 \;\;\text{or}\;\; \theta = \pi\},
        \end{equation*}
which is neither connected nor bounded.
        \end{example}
The following theorem links the notions of numerical range and singularity. Apparently, the knowledge on this topic is still not complete.   
        
    \begin{proposition}
Let $P(\lambda) = \lambda^d A_d + \lambda^{d-1} A_{d-1} + \dots + A_0 \in \mathbb{C}^{n, n}[\lambda]$ be a matrix polynomial, where $A_d \neq 0$. Consider the following conditions:
        \begin{enumerate}[\rm (a)]
\item\label{sin} the polynomial $P(\lambda)$ is singular,
\item\label{iso} all coefficients $A_0, A_1, \dots, A_d$ have a common isotropic vector $x \in \mathbb{C}^n \setminus \{0\}$, i.e. \\ $x^*A_jx = 0$ for all $j \in \{1, 2, \dots, d\}$,
\item\label{fcp} the numerical range of the polynomial $P(\lambda)$ fully covers the complex plane: $W(P) = \mathbb{C}$.
        \end{enumerate}
Then \eqref{sin} $\Rightarrow$ \eqref{fcp} and \eqref{iso} $\Rightarrow$ \eqref{fcp}. Moreover, for matrix pencils $P(\lambda) = \lambda A_1 + A_0$ also \eqref{sin} $\Rightarrow$ \eqref{iso} holds.
     \end{proposition}
     \begin{proof}
First implication \eqref{sin} $\Rightarrow$ \eqref{fcp} follows from the fact that $\det P(\lambda) = 0$ for all $\lambda \in\mathbb{C}$ and from the condition \eqref{zeros} of Proposition \ref{nrp}. Second implication \eqref{iso} $\Rightarrow$ \eqref{fcp} is a consequence of the equation $x^*P(\lambda)x = x^*(\sum_{j = 0}^d \lambda^j A_j)x = \sum_{j = 0}^d (x^* A_j x)\lambda^j = 0$ and Definition \ref{numran}. For a proof of the additional implication for matrix pencils, see Section 2 of \cite{PP2023matrix}.
    \end{proof}
    \begin{example}\rm
Observe that converse implications do not hold, i.e. for statements \eqref{fcp} $\not\Rightarrow$ \eqref{sin} and \eqref{fcp} $\not\Rightarrow$ \eqref{iso} we have the following counterexample of a matrix pencil, see \cite[Example 2.1]{PP2023matrix}:
        \begin{equation*}
P(\lambda) = \lambda
            \begin{bmatrix}
0 & 0\\
2 & 0
            \end{bmatrix} 
+ 
            \begin{bmatrix}
1 & 0\\
0 & -1
            \end{bmatrix},
        \end{equation*}
which is regular, since $\det P(\lambda) = -1 \neq 0$, and satisfies the condition $W(P) = \mathbb{C}$, but its coefficients have no common isotropic vector $x \in \mathbb{C}^2 \setminus \{0\}$.
    \end{example}
The following example is a contribution of the Author of the Thesis.
    \begin{example}\label{sing}\rm
To show that \eqref{sin} $\not\Rightarrow$ \eqref{iso} consider the singular matrix polynomial
        \begin{equation*}
P(\lambda) = \lambda^2
            \begin{bmatrix}
1 & 0\\
0 & 0
            \end{bmatrix} 
+ \lambda
            \begin{bmatrix}
0 & 1\\
1 & 0
            \end{bmatrix} 
+ 
            \begin{bmatrix}
0 & 0\\
0 & 1
            \end{bmatrix}.
        \end{equation*}
Taking $x = [x_1 \;x_2]^{\top} \in \mathbb{C}^2$ we obtain that $x^*A_2x = 0$ if and only if $x_1 = 0$ and $x^*A_0x = 0$ if and only if $x_2 = 0$. Therefore, all coefficients $A_0, A_1, A_2$ do not have a common isotropic vector $x \in \mathbb{C}^2 \setminus \{0\}$. Also \eqref{iso} $\not\Rightarrow$ \eqref{sin} and it is sufficient to take
    \begin{equation*}
P(\lambda) = \lambda
            \begin{bmatrix}
0 & 1\\
-1 & 0
            \end{bmatrix} 
+ 
            \begin{bmatrix}
0 & -1\\
1 & 0
            \end{bmatrix}.
    \end{equation*}
as a counterexample. This matrix pencil is regular, since $\det P(\lambda) = (\lambda - 1)^2 \not\equiv 0$, and each nonzero real vector $x \in \mathbb{R}^2 \setminus \{0\}$ is an isotropic vector of its coefficients.
    \end{example}

\section{Some results on stability of certain matrix polynomials}

Finally, we close this introductory chapter with citing three results by Mehl, Mehrmann and Wojtylak that initiated the research of the current Thesis. Their extensions can be found in Chapter~\ref{sPos}. We will sometimes use the notation $A \leq B$ to express that the matrix $B - A$ is positive semi-definite. The first result comes from \cite{MehMW22}.
    \begin{theorem}\label{cube}
If $A_0, A_1 \dots, A_d$ are Hermitian positive semi-definite matrices, then the eigenvalues of the matrix polynomial $P(\lambda) = \lambda^d A_d + \lambda^{d-1} A_{d-1} + \dots + A_0 \in \mathbb{C}^{n, n}[\lambda]$ are located in the angle $\{\lambda \in \mathbb{C}: |\Arg\lambda| \geq \pi/d\} \cup \{0\}$.
    \end{theorem}
The second one comes from \cite{MehMW18}.
    \begin{theorem}\label{quad}
If $A_0, A_2, \frac{1}{2}(A_1+A_1^*)$ are Hermitian positive semi-definite, then the eigenvalues of the quadratic matrix polynomial $P(\lambda) = \lambda^2 A_2 + \lambda A_1 + A_0 \in \mathbb{C}^{n, n}[\lambda]$ are located in the closed left half-plane $\overline{H}_{\frac{3}{2}\pi}$. 
    \end{theorem}
The third one comes again from \cite{MehMW22}.
    \begin{theorem}
If $A_0, A_1, A_2, A_3$ are Hermitian positive definite and $A_2 \geq A_3$, $A_1 \geq A_0$, then the eigenvalues of the cubic matrix polynomial $P(\lambda) = \lambda^3 A_3 + \lambda^2 A_2 + \lambda A_1 + A_0 \in \mathbb{C}^{n, n}[\lambda]$ are located in the open left half-plane $H_{\frac{3}{2}\pi}$.
    \end{theorem}
At the end of this Section, we present a partial formulation of Theorem 16 from \cite{MehMW21} which we will need in Chapter ~\ref{sPos}.
    \begin{theorem}\label{kernels}
Let $k \in \mathbb{Z}_+$ and $j \in \{0, \dots , k\}$. If $P(\lambda) = -\lambda^jJ + \sum_{i=0}^k\lambda^iA_i$ with $J, A_0, \dots , A_k \in \mathbb{R}^{n, n}, J^{\top} = -J, A_i^{\top} = A_i \geq 0$ for $i \in \{0, \dots , k\}$, then the matrix polynomial $P(\lambda)$ is singular, i.e. $\det P(\lambda) \equiv 0$, if and only if $\ker J \cap \ker A_0 \cap\dots\cap \ker A_k \neq \{0\}$.
    \end{theorem}

\chapter{Hyperstability of matrix polynomials}\label{sHyper}

We emphasise now two central notions. While the definition of stability is a usual one, the definition of hyperstability is a contribiution of the Dissertation.
	\section{Definition, basic properties}  
	\begin{definition}\label{defla}
\rm Let $D$ denote any nonempty open or closed subset of the complex plane $\mathbb{C}$ and let $P(\lambda) \in \mathbb{C}^{n, n}[\lambda]$ be a regular matrix polynomial. We say that the polynomial $ P(\lambda)$ is \emph{stable with respect to $D$} if and only if a scalar function $\lambda \mapsto \det P(\lambda)$ does not have zeors in $D$.
Further, we say that the polynomial $P(\lambda)$ \emph{is hyperstable with respect to $D$} if and only if for all $x \in \mathbb{C}^n \setminus\{0\}$ there exists $y \in \mathbb{C}^n \setminus \{0\}$ such that 
		\begin{equation}\label{nozeroinD}
y^*P(\mu)x \neq 0 \;\;\text{for all} \;\;\mu \in D\text{.}
		\end{equation}
	\end{definition}
		
Observe that the stability of a matrix polynomial implies its regularity. Further, as mentioned in the Introduction, the notion of hyperstability is situated somewhere between stability and a  numerical range condition. We provide a detailed analysis below.
\begin{proposition}\label{abc}
Let $P(\lambda) \in \mathbb{C}^{n, n}[\lambda]$ be a matrix polynomial and let $D$ be a nonempty open or closed subset of the complex plane $\mathbb{C}$. Consider the following conditions.
\begin{enumerate}[\rm (a)]
\item\label{num}  the numerical range $W(P)$ does not intersect $D$;
\item\label{hyper} the polynomial $P(\lambda)$ is hyperstable with respect to $D$; 
\item\label{stable} the polynomial $ P(\lambda)$ is stable with respect to $D$.
\end{enumerate}
Then \eqref{num}$\;\Rightarrow\;$\eqref{hyper}$\;\Rightarrow\;$\eqref{stable}.
\end{proposition}
		
\begin{proof}
The first implication  follows by setting $y=x$ in \eqref{nozeroinD}. The second implication becomes obvious, as one observes that stability can be reformulated as:
for all $\mu\in D$ and for all $x \in \mathbb{C}^n \setminus\{0\}$ there exists $y \in \mathbb{C}^n \setminus \{0\}$ such that $y^*P(\mu)x \neq 0$.
\end{proof}		

 It appears that in many cases the implication \eqref{num} $\Rightarrow$ \eqref{stable} above is a convenient criterion for stability, see \cite{MehMW22}. Note that the numerical range of a polynomial $P(\lambda)=\lambda I_n - A$ coincides with the numerical range of the matrix $A$. Thus, condition \eqref{num} cannot be, in general, equivalent to stability \eqref{stable}. As for matrix pencils hyperstability \eqref{hyper} is equivalent to 
stability \eqref{stable} (see Theorem~\ref{uppert}\eqref{pq} below), we obtain that in many cases \eqref{num} is not  equivalent to \eqref{hyper}. Further, hyperstability \eqref{hyper} is,  in general, not equivalent to stability \eqref{stable}, cf. the  following  important example.

    \begin{example}\label{exa}\rm	
Let us take
        \begin{equation*}
P(\lambda) := \lambda^2
            \begin{bmatrix}
0 & 0\\
0 & 1
            \end{bmatrix}
\! + \!\lambda
            \begin{bmatrix}
0 & 1\\
1 & 0
            \end{bmatrix}
\! + \!
            \begin{bmatrix}
1 & 0\\
0 & 1
            \end{bmatrix}
 = 
            \begin{bmatrix}
1 & \lambda\\
\lambda & \lambda^2 + 1
            \end{bmatrix}\text{.}
        \end{equation*}
Then $\det P(\lambda)\equiv 1$, hence the matrix polynomial $P(\lambda)$ is stable with respect to any open or closed subset $D$ of the complex plane $\mathbb{C}$.
However, taking $x=[0 \; 1]^{\top}$ and arbitrary $y = [y_1 \; y_2]^{\top}$ yields 
$$
y^*P(\lambda)x = \overline{y}_1 \lambda + \overline{y}_2 (\lambda^2 + 1)\text{.}  
$$	
Note that for any $y \neq 0$ the scalar polynomial on the right side of the above equation has always a root in the closed unit disc $\overline{\mathbb{D}}$, as either $y_2 = 0$ and we have one root $\lambda=0$ or $y_2\neq 0$ and we have two roots (counting multiplicities) with the product equal to $1$, so at least one of them must be in $\overline{\mathbb{D}}$. Therefore, the polynomial $P(\lambda)$ is stable, but not hyperstable on $\overline{\mathbb{D}}$. Replacing $\lambda$ by $\alpha\lambda + \beta$ one obtains an example of a polynomial which is stable, but not hyperstable, with respect to an arbitrary desired half-plane or disc. E.g., the polynomial $P(\lambda-\ii)$  is stable, but not hyperstable, with respect to the upper half-plane $H_0$.

Observe, that the polynomial $P(\alpha\lambda+\beta)$ ($\alpha\neq0$) has infinity as an eigenvalue, i.e., the leading coefficient $\diag(0,\alpha^2)$ is not an invertible matrix. One may wonder if having infinity as an eigenvalue is necessary to construct such an example. Apparently, this is not the case, see Example~\ref{hyper-nsinf} below.
    \end{example}

Let us recall two notions of equivalence for matrix polynomials, cf. \cite{Gan59}.
We say that matrix polynomials $P(\lambda), Q(\lambda) \in \mathbb{C}^{n, n}[\lambda]$ are \emph{equivalent} if there exist  matrix polynomials  $U(\lambda), V(\lambda) \in \mathbb{C}^{n, n}[\lambda]$ with constant nonzero determinants such that $P(\lambda) = U(\lambda)Q(\lambda)V(\lambda)$. We say that the polynomials $ P(\lambda)$ and $Q(\lambda)$ are \emph{strictly equivalent} if, additionally,  the polynomials $ U(\lambda)$ and $ V(\lambda)$ are constant invertible matrices. Clearly, the latter relation preserves hyperstability, in fact a stronger statement holds: 
\begin{lemma}\label{lQ}
Let $D$ be an open or closed subset of the complex plane $\mathbb{C}$, let $Q\in\Comp^{n,n}$ be any matrix and let $S\in\Comp^{n,n}$ be invertible. If $Q^* P(\lambda)S$  is hyperstable with respect to $D$,  then $P(\lambda)$ is hyperstable with respect to $D$.
\end{lemma}

 The proof follows directly from the definition. In a moment we will see that  equivalence does not preserve hyperstability, see Remark~\ref{eqstab}. First, however, we need to show some properties of hyperstability, connected with (block) upper-triangular matrices. 

    \begin{proposition}\label{upperblock}
Let $P(\lambda) \in \mathbb{C}^{n,n}[\lambda]$ be a matrix polynomial and let $D$ be an open or closed subset of the complex plane $\mathbb{C}$. Assume that the polynomial $ P(\lambda)$ is strictly equivalent to a block upper-triangular matrix polynomial 
        \begin{equation}\label{upper_block_matrix}
 \mat{ccc} P_{11}(\lambda) & \cdots & P_{1m}(\lambda)\\ &\ddots&\vdots\\ && P_{mm}(\lambda)  \rix,\quad P_{ij}(\lambda)\in\Comp^{k_i,k_j}[
 \lambda]{, \;\sum_{j=1}^m k_j = n,}
        \end{equation}
with the diagonal entries being hyperstable with respect to $D$. Then the polynomial $P(\lambda)$  is hyperstable with respect to $D$.
    \end{proposition}
	\begin{proof}
As strict equivalence preserves hyperstability, we may assume without loss of generality that the polynomial $P(\lambda)$ is of the form \eqref{upper_block_matrix}. Take $x = \mat{c} x_1^\top  \dots x_m^\top \rix^\top \neq 0$,  with 
$x_j\in \Comp^{k_j}$ and let $r \in \set{1, 2, \dots , m}$ denote the index of the last nonzero  $x_j$. Let $y_r\in\Comp^{k_r}\setminus\set0$ be such that the polynomial $ y_r^* P_{rr}(\lambda)x_r$ is stable with respect to $D$. Taking $y=\mat{c} 0 \dots 0\ y_r^\top 0 \dots 0 \rix^\top  $  we obtain  $y^* P(\lambda)x=y_r^* P_{rr}(\lambda)x_r$, which  is stable with respect to $D$.  
	
	\end{proof}

	\begin{remark}\rm
Proposition~\ref{upperblock} relates somehow to the current work on the (generalised) triangularisation of matrix polynomials, see \cite{ang2021, trian2013, tis2013}. However, note that in these papers one transforms a given matrix polynomials to a  quasi-triangular one by equivalence transformations, and not by strict equivalence. In view of the current paper, strict equivalence with a block upper-triangular matrix polynomial is a rather strong property, which is not easy to obtain. Thus, this will not be the objective below. Instead, we will provide several other ways of showing hyperstability.
	\end{remark}

	The following Theorem presents a  class of matrix polynomials for which the notions of stability and hyperstability coincide. For applications  see Proposition~\ref{MGT} below.

\begin{theorem}\label{uppert}
Let $P(\lambda) \in \mathbb{C}^{n, n}[\lambda]$ be a matrix polynomial and let $D$ be a nonempty open or closed subset of the complex plane $\mathbb{C}$, then the following holds.
\begin{enumerate}[\rm (i)]
\item\label{ut} Assume that a matrix polynomial $P(\lambda)$ is strictly equivalent to an upper-triangular matrix polynomial. Then it is hyperstable with respect to $D$ if and only if it is stable with respect to $D$.
\item\label{pq} If $P(\lambda) = p(\lambda)A + q(\lambda)B$ with some scalar polynomials $p(\lambda), q(\lambda) \in \mathbb{C}[\lambda]$ and $A, B \in \Comp^{n, n}$ (in particular: if $P(\lambda)$ is a matrix pencil), then it is hyperstable with respect to $D$ if and only if it is stable with respect to $D$.
\end{enumerate}
\end{theorem}
	
	\begin{proof}
\eqref{ut} Consider a matrix polynomial $P(\lambda)$ stable with respect to $D$. It follows from the assumption on $P(\lambda)$ that there exist invertible matrices $U, V \in \mathbb{C}^{n, n}$ such that
$$
U^{-1}P(\lambda)V^{-1} = \mat{ccc} p_{11}(\lambda) & \cdots & p_{1n}(\lambda)\\ &\ddots&\vdots\\ && p_{nn}(\lambda)  \rix
$$
with a stable upper-triangular matrix polynomial on the right-hand side of the equation. Hence, the scalar polynomials on the diagonal $p_{11}(\lambda), p_{22}(\lambda), \dots , p_{nn}(\lambda) \in \mathbb{C}[\lambda]$ are stable with respect to $D$. By Proposition~\ref{upperblock}, $P(\lambda)$ is hyperstable with respect to $D$.   
	
\eqref{pq} Using the Kronecker form \cite{Gan59} (or the generalised Schur form, see Definition ~\ref{Schur}) of matrices $A$ and $B$ we obtain the matrix polynomial  $\lambda \mapsto p(\lambda)A + q(\lambda)B$ to be strictly equivalent to an upper-triangular. Hence,  the claim follows from \eqref{ut}.
	\end{proof}

\begin{remark}\label{eqstab}\rm
Note that an arbitrary matrix polynomial $P(\lambda) \in \mathbb{C}^{n, n}[\lambda]$ is equivalent to an upper triangular matrix polynomial (with $V(\lambda) = I_n$), see \cite[Chapter VI]{Gan59}. E.g., the polynomial from Example~\ref{exa} satisfies the following equality:
$$
\mat{cc} 1 & \lambda\\ \lambda & \lambda^2 + 1\rix = \mat{cc} 1 & 0\\ \lambda & 1\rix\mat{cc} 1 & \lambda\\ 0 & 1\rix\text{.}
$$
The upper-triangular matrix polynomial on the right-hand side of the above equation is stable on $\overline{\mathbb{D}}$, since its determinant is identically $1$, and consequently it is hyperstable on $\overline{\mathbb{D}}$ (cf. Proposition \ref{upperblock}). However, the polynomial on the left side is not hyperstable on $\overline{\mathbb{D}}$. Hence, the equivalence of matrix polynomials does not preserve hyperstability.  

Also let us note that by Theorem~\ref{uppert}\eqref{ut} the above-mentioned  polynomial is  not strictly equivalent to an upper-triangular one.

\end{remark}

\section{Orbits of hyperstable polynomials}
As it was shown before, hyperstability is not invariant under the equivalence $E(\lambda)P(\lambda)F(\lambda)$, where $E(\lambda), F(\lambda) \in \mathcal{U}_n$. We study this problem in more detail, by investigating orbits $\{E(\lambda)P(\lambda)F(\lambda) \in \mathbb{C}^{n, n}[\lambda] : E(\lambda), F(\lambda) \in \mathcal{U}_n\}$ of hyperstable matrix polynomials.

    \begin{proposition}
Let $D$ be a nonempty open or closed subset of the complex plane $\mathbb{C}$, $P(\lambda) \in \mathbb{C}^{n, n}[\lambda]$ be an arbitrary matrix polynomial and $S(\lambda) \in \mathbb{C}^{n, n}[\lambda]$ be its Smith canonical form. The following conditions are equivalent:
        \begin{enumerate}[\rm (a)]
\item\label{orb} there exists a hyperstable (with respect to $D$) element of the orbit $\{E(\lambda)P(\lambda)F(\lambda) \in \mathbb{C}^{n, n}[\lambda] : E(\lambda), F(\lambda) \in \mathcal{U}_n\}$ of the polynomial $P(\lambda)$,
\item\label{stab} the polynomial $P(\lambda)$ is stable with respect to $D$,
\item\label{ss} the Smith form $S(\lambda)$ of the polynomial $P(\lambda)$ is stable with respect to $D$,
\item\label{shs} the Smith form $S(\lambda)$ of the polynomial $P(\lambda)$ is hyperstable with respect to $D$.
        \end{enumerate}
    \end{proposition}
    \begin{proof}
\eqref{orb} $\Rightarrow$ \eqref{stab} Since a matrix polynomial $E(\lambda)P(\lambda)F(\lambda)$ is hyperstable with respect to $D$ for some unimodular transformations $E(\lambda), F(\lambda) \in \mathcal{U}_n$, so it is stable with respect to $D$ as well and consequently each its factor is stable with respect to $D$. In particular, the polynomial $P(\lambda)$ is stable with respect to $D$. \\\\
\eqref{stab} $\Rightarrow$ \eqref{ss} We know that the Smith form $S(\lambda)$ is equivalent to the polynomial $P(\lambda)$ which is stable with respect to $D$. Since unimodular equivalence preserves stability, the Smith form $S(\lambda)$ is also stable with respect to $D$. \\\\
\eqref{ss} $\Rightarrow$ \eqref{shs} This implication becomes obvious, when we recall Theorem \ref{uppert}\eqref{ut}. Since the Smith form $S(\lambda)$ is a diagonal (in particular: upper-triangular) matrix polynomial, then it is hyperstable with respect to $D$ if and only if it is stable with respect to $D$. \\\\
\eqref{shs} $\Rightarrow$ \eqref{orb} Because of the fact that Smith form $S(\lambda)$ belongs to the orbit of the polynomial $P(\lambda)$, the last implication follows.
    \end{proof}
The following theorem is a contribution of the Author of the Thesis.
    \begin{theorem}\label{orbits}
Let $D$ be an open or closed disk or open or closed half plane. There is no matrix polynomial $P(\lambda) \in \mathbb{C}^{n, n}[\lambda]$ such that each element of its orbit                      $\{E(\lambda)P(\lambda)F(\lambda) \in \mathbb{C}^{n, n}[\lambda] : E(\lambda), F(\lambda) \in \mathcal{U}_n\}$ is hyperstable with respect to a set $D$.
    \end{theorem}
    \begin{proof}
Let us fix a set $D = \overline{\mathbb{D}}$ and a natural number $n \geq 2$. Consider a matrix polynomial $P(\lambda)$ with square coefficients of the size $n \times n$. If the polynomial $P(\lambda)$ is not hyperstable with respect to $\overline{\mathbb{D}}$, then trivially its orbit contains a polynomial not hyperstable with respect to $\overline{\mathbb{D}}$. Hence, assume that the polynomial $P(\lambda)$ is hyperstable with respect to $\overline{\mathbb{D}}$ and take its Smith canonical form $S(\lambda)$. Since the Smith form $S(\lambda)$ is equvalent to the stable with respect to $\overline{\mathbb{D}}$ polynomial $P(\lambda)$ and unimodular equvalence preserves stability, then the Smith form $S(\lambda)$ is also stable with respect to $\overline{\mathbb{D}}$. Thus, we have:
        \begin{equation*}
S(\lambda) = 
            \begin{bmatrix}
s_1(\lambda) & 0 & \dots & 0\\
0 & s_2(\lambda) & \dots & 0\\
\vdots & \vdots & \ddots & \vdots\\
0 & 0 & \dots & s_n(\lambda)
            \end{bmatrix}\text{,}
        \end{equation*}
where $s_{j+1}(\lambda) | s_j(\lambda)$ for all $j \in \{1, 2, \dots , n-1\}$.
        
From hyperstability of the Smith form $S(\lambda)$ with respect to $\overline{\mathbb{D}}$, we conclude stability of each invariant polynomials $s_1(\lambda), s_2(\lambda), \dots , s_n(\lambda)$ with respect to $\overline{\mathbb{D}}$. Now, let $d := \deg(s_1/s_2) + 2$ and put
        \begin{equation*}
E(\lambda) = 
            \begin{bmatrix}
1 & 0 & \dots & 0\\
0 & 1 & \dots & 0\\
\vdots & \vdots & \ddots & \vdots\\
0 & 0 & \dots & 1
            \end{bmatrix}\;\text{and}\;\;
F(\lambda) = 
            \begin{bmatrix}
1 & \lambda & 0 & \dots & 0\\
\lambda^{d-1} & \lambda^d + 1 & 0 & \dots & 0\\
0 & 0 & 1 & \dots & 0 \\
\vdots & \vdots & \vdots & \ddots & \vdots\\
0 & 0 & 0 & \dots & 1
            \end{bmatrix}\text{.}
        \end{equation*}
From the Laplace'a determinant formula, repeatedly applied to the last row, we obtain immediately
        \begin{equation*}
\det F(\lambda) = 
            \begin{vmatrix}
1 & \lambda\\
\lambda^{d-1} & \lambda^d + 1
            \end{vmatrix}
 = 1\text{.}
        \end{equation*}
\noindent Because $E(\lambda), F(\lambda) \in \mathbb{C}^{n, n}[\lambda]$ and $\det E(\lambda) = \det F(\lambda) = 1$, so we have $E(\lambda), F(\lambda) \in \mathcal{U}_n$. It is easy to check that
        \begin{align*}
E(\lambda)S(\lambda)F(\lambda) &= 
            \begin{bmatrix}
s_1(\lambda) & 0 & 0 & \dots & 0\\
0 & s_2(\lambda) & 0 & \dots & 0\\
0 & 0 & s_3(\lambda) & \dots & 0\\
\vdots & \vdots & \vdots & \ddots & \vdots\\
0 & 0 & 0 & \dots & s_n(\lambda)
            \end{bmatrix}
            \begin{bmatrix}
1 & \lambda & 0 & \dots & 0\\
\lambda^{d-1} & \lambda^d + 1 & 0 & \dots & 0\\
0 & 0 & 1 & \dots & 0 \\
\vdots & \vdots & \vdots & \ddots & \vdots\\
0 & 0 & 0 & \dots & 1
            \end{bmatrix}\\ 
&= 
            \begin{bmatrix}
s_1(\lambda) & \lambda s_1(\lambda) & 0 & \dots & 0\\
\lambda^{d-1}s_2(\lambda) & (\lambda^d + 1)s_2(\lambda) & 0 & \dots & 0\\
0 & 0 & s_3(\lambda) & \dots & 0\\
\vdots & \vdots & \vdots & \ddots & \vdots\\
0 & 0 & 0 & \dots & s_n(\lambda)
            \end{bmatrix}\text{,}
        \end{align*}
which yields for $x = e_2 = [0 \;1 \dots 0]^{\top}$:
        \begin{eqnarray*}
y^*E(\mu)S(\mu)F(\mu)x &=& 
            \begin{bmatrix}
\overline{y}_1 & \overline{y}_2 & \overline{y}_3\dots & \overline{y}_n
            \end{bmatrix}
            \begin{bmatrix}
\mu s_1(\mu)\\
(\mu^d + 1)s_2(\mu)\\
0 \\
\vdots\\
0
            \end{bmatrix}\\ 
&=& \overline{y}_1 \mu s_1(\mu) + \overline{y}_2(\mu^d + 1)s_2(\mu)\text{.}
        \end{eqnarray*}
We will show that for all $y \in \mathbb{C}^n \setminus \{0\}$, there exists $\mu \in \overline{\mathbb{D}}$ such that $\overline{y}_1 \mu s_1(\mu) + \overline{y}_2(\mu^d + 1)s_2(\mu) = 0$. Since $s_2(\lambda) | s_1(\lambda)$ and $s_2(\mu) \neq 0$ for all $\mu \in \overline{\mathbb{D}}$, then the above claim is equivalent to the following one: for all $y \in \mathbb{C}^n \setminus \{0\}$, there exists $\mu \in \overline{\mathbb{D}}$ such that $\overline{y}_1\mu (s_1/s_2)(\mu) + \overline{y}_2(\mu^d + 1) = 0$. If $y_2 = 0$, then it is sufficient to take $\mu = 0$ to satisfy the last condition. If $y_2 \neq 0$, then a scalar polynomial $p(\lambda) = \overline{y}_1\lambda (s_1/s_2)(\lambda) + \overline{y}_2(\lambda^d + 1)$ has degree equal to $d$. Let $\lambda_1, \lambda_2, \dots , \lambda_d$ be all its complex roots. The first and the last coefficient of the polynomial $p(\lambda)$ are both equal to $\overline{y}_2$. Thus, from the Vieta's formula on the product of roots, we have $\lambda_1\lambda_2\dots\lambda_d = (-1)^d$ and consequently $|\lambda_1|\cdot|\lambda_2|\cdot...\cdot|\lambda_d| = 1$. Among complex numbers $\lambda_1, \lambda_2, \dots , \lambda_d$ there is a number with the minimum module, say $\lambda_m$ ($m \in \{1, 2, \dots , d\}$). Because $|\lambda_m|^d \leq |\lambda_1|\cdot|\lambda_2|\cdot...\cdot|\lambda_d| = 1$, so $|\lambda_m| \leq 1$. Hence, the polynomial $p(\lambda)$ has a root $\mu = \lambda_m$ in the disk $\overline{\mathbb{D}}$. In consequence, we have proved that the matrix polynomial $E(\lambda)S(\lambda)F(\lambda)$ is not a hyperstable (with respect to $\overline{\mathbb{D}}$) element of the orbit of the polynomial $P(\lambda)$. 

In order to apply the above argument for any open or closed disk or half-plane, we need to make a linear change of variables and consider a matrix polynomial $P(\alpha\lambda + \beta)$, where $\alpha, \beta \in \mathbb{C}$ and $\alpha \neq 0$, instead of the polynomial $P(\lambda)$ and the proof is completed. 
    \end{proof}

\section{Multivariate hyperstable polynomials}
    \begin{definition}\label{defzi}
\rm Let $\kappa \in \mathbb{Z}_+$ and let $D$ be any nonempty open or closed subset of the complex space $\mathbb{C}^{\kappa}$. We say that a multivariate matrix polynomial $P(z_1, \dots , z_\kappa) \in \mathbb{C}^{n, n}[z_1, \dots , z_{\kappa}]$ \emph{is stable with respect to} $D$ if and only if for all $(w_1, w_2, \dots , w_{\kappa}) \in D$ and for all $x \in \mathbb{C}^n\setminus\{0\}$ there exists $y \in \mathbb{C}^n \setminus \{0\}$ such that
		\begin{equation}
y^*P(w_1, w_2, \dots, w_{\kappa})x \neq 0\text{.}
		\end{equation}
Further, we say that $P(z_1, \dots , z_{\kappa})$ is \emph{hyperstable with respect to} $D$ if and only if for all $x \in \mathbb{C}^n\setminus\{0\}$ there exists $y \in \mathbb{C}^n\setminus\{0\}$ such that 
		\begin{equation}\label{nozeroinDmatrix}
y^*P(w_1, w_2, \dots , w_{\kappa})x \neq 0 \;\;\text{for all} \;\;(w_1, w_2, \dots , w_{\kappa}) \in D\text{.}
		\end{equation}
    \end{definition}

    \begin{example}\rm
Let us examine diagonal multivariate polynomials. Let $n \in \mathbb{Z}_+$ and consider the matrix polynomial
        \begin{equation*}
P(z_1, \dots , z_n) =
            \begin{bmatrix}
z_1 & 0 & \dots & 0 \\
0 & z_2 & \dots & 0 \\
\vdots & \vdots & \ddots & \vdots \\
0 & 0 & \dots & z_n
            \end{bmatrix}\text{.}
        \end{equation*}
We can see, directly from the last definition, that $P(z_1, \dots, z_n)$ is hyperstable with respect to $(\mathbb{C}\setminus\{0\})^n$. Indeed, because $x \in \mathbb{C}^n\setminus\{0\}$ we have $x_j \neq 0$ for some $j \in \{1, 2, \dots , n\}$. For such $x$ we take $y = e_j$ ($j$-th vector of the canonical basis) and then $y^*P(z_1, \dots , z_n)x = x_1 \overline{y}_1 z_1 + x_2 \overline{y}_2 z_2 + \dots + x_n \overline{y}_n z_n = x_j z_j \neq 0$ for all $(z_1, \dots , z_n) \in (\mathbb{C}\setminus\{0\})^n$.
    \end{example}
Next proposition is a multi-variable analogue of Proposition \ref{upperblock}.
    \begin{proposition}\label{multiupperblock}
Let $P(z_1, \dots , z_{\kappa}) \in \mathbb{C}^{n,n}[z_1, \dots , z_{\kappa}]$ be a $\kappa$-variable matrix polynomial and let $D$ be a nonempty open or closed subset of the complex space $\mathbb{C}^{\kappa}$. Assume that the polynomial $ P(z_1, \dots , z_{\kappa})$ is strictly equivalent to a block upper-triangular matrix polynomial, i.e.
        \begin{equation}\label{upper_block_matrix}
 P(\lambda) = S\mat{ccc} P_{11}(z_1, \dots , z_{\kappa}) & \cdots & P_{1m}(z_1, \dots , z_{\kappa})\\ &\ddots & \vdots\\ && P_{mm}(z_1, \dots , z_{\kappa})\rix T\text{,}
        \end{equation}
where $P_{ij}(z_1, \dots , z_{\kappa}) \in \Comp^{k_i,k_j}[z_1, \dots , z_{\kappa}]$, $k_1 + k_2 + \dots + k_m = n$ and $S, T \in \mathbb{C}^{n, n}$ are invertible. If the diagonal entries are hyperstable with respect to $D$, then the polynomial $P(z_1, \dots , z_{\kappa})$ is hyperstable with respect to $D$ as well.
    \end{proposition}
	\begin{proof}
We may assume without loss of generality that $S = T = I$. Take $x = \mat{c} x_1^\top  \dots x_m^\top \rix^\top \neq 0$,  with 
$x_j\in \Comp^{k_j}$ and let $r \in \set{1, 2, \dots , m}$ denote the index of the last nonzero  $x_j$. Let $y_r\in\Comp^{k_r}\setminus\set0$ be such that the polynomial $ y_r^* P_{rr}(z_1, \dots , z_{\kappa})x_r$ is stable with respect to $D$. Taking $y=\mat{c} 0 \dots 0\ y_r^\top 0 \dots 0 \rix^\top  $  we obtain  $y^* P(z_1, \dots , z_{\kappa})x = y_r^* P_{rr}(z_1, \dots , z_{\kappa})x_r$, which is stable with respect to $D$.
	\end{proof}
 The proposition stated above allows us to give another example of multivariate hyperstable matrix polynomial which is upper-triangular.
    \begin{example}\rm
The following polynomial, which has $z_1$ on the main diagonal, $z_2$ on the super-diagonal, etc. is hyperstable with respect to the set $(\mathbb{C}\setminus\{0\}) \times \mathbb{C}^{n-1}$:
        \begin{equation*}
P(z_1, \dots , z_n) =
            \begin{bmatrix}
z_1 & z_2 & z_3 & \dots & z_n \\
0 & z_1 & z_2 & \dots & z_{n-1} \\
0 & 0 & z_1 & \dots & z_{n-2} \\
\vdots & \vdots & \vdots & \ddots & \vdots \\
0 & 0 & 0 & \dots & z_1
            \end{bmatrix}\text{.}
        \end{equation*}
It is due to Preposition \ref{multiupperblock}, because the diagonal entry $p(z_1, \dots , z_{\kappa}) = z_1$ is a stable polynomial with respect to $(\mathbb{C}\setminus\{0\}) \times \mathbb{C}^{n-1}$ (respectively hyperstable treated as a one-entry matrix).
    \end{example}

\chapter[Gauss-Lucas Theorem for matrix polynomials]{Gauss-Lucas Theorem for matrix polynomials}\label{sGL}

\section{Gauss-Lucas Theorem for univariate matrix polynomials}

We present here an important extension of a major result for scalar polynomials onto matrix polynomials. 
The Gauss-Lucas Theorem was reviewed in the first chapter (Theorem~\ref{GLscalar}). Using the current terminology, one can reformulate it as follows.
\begin{theorem}\label{GLreform}
If $D\subseteq\Comp$ is such that $\Comp\setminus D$ is convex, and $p(\lambda)$ is non-constant and stable with respect to $D$, then $p'(\lambda)$ is stable with respect to $D$. 
\end{theorem}

An analogous  similar statement is not true for matrix polynomials,  we present a very general example which will serve in future constructions.
\begin{example}\label{nonGL}\rm
Consider a polynomial 
$$
P(\lambda)=\mat{cc} \lambda p(\lambda)+q(\lambda) & p(\lambda)\\ \lambda & 1\rix,
$$
where $q(\lambda)$ has its roots outside $D$ and $p'(\lambda)$ is non-constant and has some of its roots inside $D$. Then
$$
\det P(\lambda)=q(\lambda),\quad \det P'(\lambda) = -p'(\lambda),
$$
i.e., the polynomial $P(\lambda)$ is stable with respect to $D$ but its derivative $P'(\lambda)$ is not.
\end{example}

The example above says that, in general, there is no relation between the location of the eigenvalues of a matrix polynomial and of its derivative.
However, a hyperstable version of the Gauss-Lucas Theorem holds. Below, by saying 'linear independent polynomials', we mean linear independence in the complex vector space $\mathbb{C}[\lambda]$. 

\begin{theorem}\label{GLmat}
 Let $D\subseteq\Comp$ be a nonempty open or closed set such that $\Comp\setminus D$ is convex. If a matrix polynomial $P(\lambda)$ is hyperstable with respect to $ D$ and the entries of its derivative $P'(\lambda)$ are linearly independent polynomials, then the matrix polynomial $P'(\lambda)$ is also hyperstable with respect to $D$.  
\end{theorem}

\begin{proof}
Fix a nonzero $x\in\Comp^n\setminus\{0\}$. As $P(\lambda)$ is hyperstable with respect to $D$, there exists $y\in\Comp^n\setminus\{0\}$ such that a scalar polynomial $p(\lambda)=y^*P(\lambda)x$ has all roots in $\Comp\setminus D$. Note that the polynomial $p'(\lambda)=y^*P'(\lambda)x$ is a non-trivial linear combination of the entries of $P'(\lambda)$, hence $p'(\lambda)$ is a nonzero polynomial.
As $\Comp\setminus D$ is convex, by Theorem~\ref{GLreform}, $p'(\lambda)$ has all roots in $\Comp\setminus D$.
\end{proof}

Firstly,  note that the assumption on the derivative $P'(\lambda)$ having independent entries cannot be dropped:

\begin{example}\rm
The polynomial $P(\lambda)=\diag(\lambda, 1)$ is hyperstable with respect to the outside of the unit disc $\mathbb{D}$, but its derivative $P'(\lambda)$ is singular, hence it is not hyperstable with respect to any nonempty set.  
\end{example}

Secondly, note that the assumption on $P(\lambda)$ being hyperstable cannot be relaxed to stability, even if we keep the assumption on the entries of $P'(\lambda)$ being linearly independent and relax the claim to $P'(\lambda)$ being stable with respect to $D$.

\begin{example}\label{hyper-nsinf}\rm
We specify and modify the polynomial $P(\lambda)$ from Example~\ref{nonGL}. Let $D=\{z\in\Comp:|z|>1\}$, $q(\lambda)=\lambda^2$, $p(\lambda)=\lambda^3-4\lambda$. Consider a perturbed polynomial
$$
P_\eps(\lambda):= 
\mat{cc} \lambda p(\lambda)+q(\lambda) & p(\lambda)\\ \lambda & 1+\eps \lambda^4\rix= \mat{cc} \lambda^4-3\lambda^2 & \lambda^3-4\lambda\\ \lambda & 1+\eps \lambda^4\rix.
$$
Fix  $\eps>0$ sufficiently small, so that  $P_\eps(\lambda)$ has still its eigenvalues inside the unit disc $\mathbb{D}$ and $P'_\eps(\lambda)$ has still its eigenvalues outside the closed unit disc $\overline{\mathbb{D}}$. This choice of $\eps$ is possible due to the continuity of zeros of a polynomial with respect to its coefficients. Thus, the polynomial $P'_\eps(\lambda)$ is not hyperstable on $D$. Further, note that
$$
P'_\eps(\lambda) =  \mat{cc} 4\lambda^3-6\lambda & 3\lambda^2-4\\ 1 & 4\eps \lambda^3\rix,
$$
which clearly has linearly independent entries. By 
Theorem~\ref{GLmat} the polynomial $P_\eps(\lambda)$ is not hyperstable with respect to $D$.
Further, note that the leading coefficient of this polynomial is the invertible matrix $\diag(1, \eps)$, which gives the desired example of a stable, but not hyperstable polynomial, with no eigenvalues at infinity. 
\end{example}

Let us notice that the notions of linear independence of entries of a matrix polynomial and its regularity are quite different. The following two simple examples show that neither linear independence  implies regularity nor regularity implies linear independence.
\begin{example}\label{stab1}\rm
   Consider two polynomials
   \begin{equation}
P_1(\lambda) =
        \begin{bmatrix}
1 & \lambda \\
\lambda^2 & \lambda^3
        \end{bmatrix},\text{ and }
P_2(\lambda) =
        \begin{bmatrix}
\lambda & 0 \\
0 & 1
        \end{bmatrix}\text{.}
    \end{equation}
    Then the entries of $P_1(\lambda)$ are linearly independent polynomials (over the field $\Comp$), but $P_1(\lambda)$ is singular. The polynomial $P_2(\lambda)$ is clearly regular, and its entries are linearly dependent due to the zeros on the off-diagonal positions. 
\end{example}
Note also a following obvious corollary of Theorem~\ref{GLmat}.
\begin{corollary}
If a matrix polynomial $P(\lambda)$ is hyperstable with respect to some nonempty open or closed set $ D$  such that $\Comp\setminus D$ is convex and the entries of its derivative $P'(\lambda)$ are linearly independent polynomials, then the matrix polynomial $P'(\lambda)$ is regular.
\end{corollary}

Note that the assumption of hyperstability is crucial here, if $P(\lambda)$ is regular and $P'(\lambda)$ has independent entries than $P'(\lambda)$ can be still singular, as shown by the following example.

\begin{example}\rm 
Let 
$$
P_0(\lambda) =
        \begin{bmatrix}
\lambda  & \frac12\lambda^2 \\
\frac 13\lambda^3 &\frac14 \lambda^4
        \end{bmatrix}.
$$
Then $P_0'(\lambda)=P_1(\lambda)$ (see Example~\ref{stab1}), which is clearly a singular polynomial, though its entries are linearly independent polynomials. In particular, $P_0(\lambda)$ cannot be hyperstable with respect to any nonempty open or closed set $ D$  such that $\Comp\setminus D$ is convex.
\end{example}


    \begin{remark}  \rm
The famous Sendov conjecture states the following.\medskip
    
{\em If all roots $\lambda_1, \lambda_2, \dots , \lambda_d$ of a polynomial $ p \in \mathbb{C}[\lambda] $  of degree $d\geq 2$  lie in the closed unit disk $ \overline{\mathbb{D}}$, then for every   $j \in \{1, 2, \dots , d\}$ there is a root $\xi$  of $p'$  such that $|\xi - \lambda_j| \leq 1$.}
\medskip

The conjecture has been verified for polynomials of degree at most $8$, cf. \cite{brown1999proof}. Terence Tao proved it for large enough $d$, see \cite{tao2020sendovs}. 
It is tempting to consider also a generalisation of the Sendov conjecture onto matrix polynomials using the hyperstability notion.  However, the main problem here would be the assumption that the degree of $p(\lambda)$ is larger or equal than 2. When presenting a similar construction as above we would have to guarantee that  $y^*P(\lambda)x$ keeps this assumption, which would lead to severe technical difficulties.  
\end{remark}

\section{Gauss-Lucas Theorem for multivariate matrix polynomials}

In this Section we show the most general version of the theorem. Let us recall the notation
    \begin{equation*}
A(z^{(j)}) := \{w \in \mathbb{C} : ({z_1}, \dots , {z_{j-1}}, w, {z_{j+1}}, \dots , {z_{\kappa}}) \in A\}\text{.}
    \end{equation*}
    \begin{definition}\label{scxx}\rm
Let $\kappa \in \mathbb{Z}_+ \setminus \{1\}$ and $j \in \{1, 2, \dots , \kappa\}$. We say that a subset $A$ of the complex space $\mathbb{C}^{\kappa}$ is \emph{separately convex with respect to the $j$-th variable in $\mathbb{C}^{\kappa}$} if and only if the set $A(z^{(j)})$ is a convex subset of the complex plane $\mathbb{C}$ for all $z \in \mathbb{C}^{\kappa}$. 
    \end{definition}
    \begin{example}
\rm Let
$$
D = \Comp^2\setminus \{(z_1,z_2):|z_1|\leq 1,\ |z_2|\leq 1\}.
$$
We have $(\mathbb{C}^2 \setminus D)(z^{(1)}) = \{z_1 \in \mathbb{C} : (z_1, {z_2}) \in \mathbb{C}^2 \setminus D\}$. When $|{z_2}| \leq 1$, then $(\mathbb{C}^2 \setminus D)(z^{(1)}) = \{z_1 \in \mathbb{C} : |z_1| \leq 1\}$. When $|{z_2}| > 1$, then $(\mathbb{C}^2 \setminus D)(z^{(1)}) = \emptyset$. In both cases, $(\mathbb{C}^2 \setminus D)(z^{(1)})$ is a convex set, therefore $\mathbb{C}^2 \setminus D$ is separately convex with respect to the first variable. 

\end{example}
The following theorem is a contribution of the Author of the Thesis.
    \begin{theorem}\label{mmgl}
Let $\kappa \in \mathbb{Z}_+ \setminus \{1\}$, $j \in \{1, 2, \dots , \kappa\}$  and let $D \subseteq \mathbb{C}^{\kappa}$ be a nonempty open or closed set such that its complement  $\mathbb{C}^{\kappa} \setminus D$ is separately convex with respect to the $j$-th variable in $\mathbb{C}^{\kappa}$.
If a matrix polynomial $P(z_1,  \dots , z_{\kappa}) \in \mathbb{C}^{n, n}[z_1,  \dots , z_{\kappa}]$ is hyperstable with respect to $D$ and the entries of its partial derivative $\frac{\partial P}{\partial z_j}(z_1, \dots , z_{\kappa})$  are linearly independent polynomials, then  $\frac{\partial P}{\partial z_j}(z_1, \dots , z_{\kappa})$ is  hyperstable with respect to $D$ as well.
    \end{theorem}
    \begin{proof}
Fix $x \in \mathbb{C}^n \setminus \{0\}$. As $P(z_1, \dots , z_{\kappa})$ is hyperstable with respect to $D$, there exists $y \in \mathbb{C}^n \setminus \{0\}$ such that a scalar multivariate polynomial 
$$
p(z_1, \dots , z_{\kappa}) = y^*P(z_1, \dots , z_{\kappa})x
$$
is stable with respect to $D$. Note that 
$$
\frac{\partial p}{\partial z_j}(z_1,  \dots , z_{\kappa}) = y^*\frac{\partial P}{\partial z_j}(z_1,  \dots , z_{\kappa})x
$$
is a non-trivial linear combination of the entries of $\frac{\partial P}{\partial z_j}(z_1, \dots , z_{\kappa})$, hence the polynomial $\frac{\partial p}{\partial z_j}(z_1,  \dots , z_{\kappa})$ is nonzero. To finish the proof it is enough to show that it is stable with respect to $D$. Consider the  univariate polynomial  
$$
q_{z^{(j)}}(w) := p({z_1},  \dots , {z}_{j-1}, w , {z}_{j+1}, \dots , {z_{\kappa}}),
$$
where  $w\in D(z^{(j)})$ 
 and $z^{(j)}$ is arbitrary, fixed.  
By definition, if $w \in D(z^{(j)})$, then the point $({z_1}, , \dots , {z}_{j-1}, w, {z}_{j+1}, \dots , {z_{\kappa}})$ belongs to $D$. 
 As $p(z_1, \dots , z_{\kappa})$ is stable with respect to $D$, we have that $q_{z^{(j)}}(w)$ is stable with respect to $D(z^{(j)})$. Since $\frac{\partial p}{\partial z_j}(z_1,  \dots , z_{\kappa})$ is nonzero, $q_{z^{(j)}}(w)$ is non-constant.  Further, note that
$$
(\mathbb{C}^{\kappa} \setminus D)(z^{(j)}) = \mathbb{C} \setminus D(z^{(j)})\text{.}
$$
The set on the left hand side is convex, which follows from the assumption of $\mathbb{C}^{\kappa} \setminus D$ being separately convex with respect to $j$-th variable. Thus we may apply the scalar Gauss-Lucas theorem (Theorem~\ref{GLscalar}) and obtain the stability of
$
q_{z^{(j)}}'(w)$ 
with respect to $D(z^{(j)})$ for any $z^{(j)} \in \mathbb{C}^{\kappa - 1}$. 
Hence, the scalar polynomial $\frac{\partial p}{\partial z_j}(z_1, \dots , z_{\kappa})$ is  stable with respect to $D$. 
 Indeed, if $\frac{\partial p}{\partial z_j}(z_1, \dots , z_{\kappa})$ had a zero $(\zeta_1, \dots , \zeta_{\kappa}) \in D$, we would have 
 $q_{z^{(j)}}'(\zeta_j) = 0$ for $\zeta_j \in D(\zeta^{(j)})$, which is impossible.  This finishes the proof.

    \end{proof}

  \begin{remark}
\rm Let us notice we have proved a result stronger than the one of Kanter \cite{kanter}, see Theorem \ref{mgl} above. Indeed, let us take a scalar $\kappa$-variable polynomial $p$. As $n = 1$, by definition (see Definition \ref{defzi}) hyperstability is equivalent to stability and linear independence of entries is equivalent to $\partial p/\partial z_j$ being nonzero. Regarding a set $D$, we take $D = \mathbb{C}^{\kappa} \setminus \conv_{\kappa}(p^{-1}(0))$. Then $\mathbb{C}^{\kappa} \setminus D = \conv_{\kappa}(p^{-1}(0))$, which is separately convex (in particular, separately convex with respect to $j$-th variable). Since a polynomial $p$ is stable with respect to $\mathbb{C}^{\kappa} \setminus \conv_{\kappa}(p^{-1}(0))$, we receive from the above Theorem~\ref{mmgl} stability of its partial derivative $\partial p/\partial z_j$ with respect to $\mathbb{C}^{\kappa} \setminus \conv_{\kappa}(p^{-1}(0))$, which means that desired inclusion in Theorem \ref{mgl} holds.
    \end{remark}

We show now that the assumption of linear independence of the entries in Theorem~\ref{mmgl} is needed.

    \begin{example}\rm 
Consider the set $\mathbb{C}^2 \setminus H_0^2$, which is separately convex with respect to the first variable. It is clear, because $(\mathbb{C}^2 \setminus H_0^2)(z^{(1)}) = \{z_1 \in \mathbb{C} : (z_1, {z_2}) \in \mathbb{C}^2 \setminus H_0^2\}$ is equal to $\mathbb{C}$, when $\im{z_2} \leq 0$ and equal to $\mathbb{C} \setminus H_0$, when $\im{z_2} > 0$. For necessity of the assumption of linear independence of the entries in Theorem \ref{mmgl}, consider the following example of two-variable polynomial being hyperstable with respect to the Cartesian square of the upper half-plane $H_0^2$:
        \begin{equation*}
P(z_1, z_2) = 
            \begin{bmatrix}
z_1^2 & 0 \\
0 & z_2^2
            \end{bmatrix}.
        \end{equation*}
Its hyperstability follows directly from the definition (cf. Definition \ref{defzi}). Indeed, we have $y^*P(z_1, z_2)x = x_1 \overline{y}_1 z_1^2 + x_2 \overline{y}_2 z_2^2$. Therefore, when $x_1 \neq 0$, we take $y = e_1 = [1\; 0]^{\top}$ to get $y^*P(z_1, z_2)x \neq 0$ for all $(z_1, z_2) \in H_0^2$ and when $x_2 \neq 0$, we take $y = e_2 = [0\; 1]^{\top}$ to satisfy the condition $y^*P(z_1, z_2)x \neq 0$ for all $(z_1, z_2) \in H_0^2$. Now, let us calculate first order partial derivative of the polynomial $P(z_1, z_2)$ with respect to the variable $z_1$:
        \begin{equation*}
\frac{\partial P}{\partial z_1}(z_1, z_2) =
            \begin{bmatrix}
2z_1 & 0 \\
0 & 0
            \end{bmatrix}\text{.}
        \end{equation*}
The last polynomial is not hyperstable with respect to $H_0^2$. As above, we have $y^*(\partial P/\partial z_1)(z_1, z_2)x = 2x_1 \overline{y}_1 z_1$. Clearly, for $x = [0\; 1]^{\top}$ there is no $y \in \mathbb{C}^2\setminus\{0\}$ such that $y^*(\partial P/\partial z_1)(z_1, z_2)x \neq 0$. A reason why it happens is the fact that the entries of $(\partial P/\partial z_1)(z_1, z_2)$ are linearly dependent and the assumption of their linear independence in Theorem \ref{mgl} has been violated.
    \end{example}


\chapter{Sz\'asz-type inequality  for stable matrix polynomials}\label{SzaszM}
The corresponding results for scalar polynomials were reviewed in Section \ref{SzaszClass}, we concentrate now on the matrix case.  Recall that  the norm $||\cdot||$ is the matrix two-norm induced by the vector two-norm and $||\cdot||_F$ is the Frobenius norm. Further, 
$$
H_\varphi := \{\lambda \in \mathbb{C} : \im(\lambda e^{i\varphi}) > 0\},\quad\varphi \in [0;2\pi),
$$
in particular $H_0$ is the open upper half-plane. 

\section{Employing  numerical range}

 Let $X$ be any matrix. The symbol $\lambda_{H}(X)$ below denotes the largest (possibly negative) eigenvalue of the Hermitian matrix $\frac{X+X^*}2$. Note that 
 \begin{equation}\label{lambdaH}
\lambda_{H}(X)=\max_{\norm x=1} x^*\left(  \frac{X+X^*}2 \right)x = \max_{\norm x=1} (  \re  x^* X x    ).
 \end{equation}
      
    \begin{theorem}\label{Szasz}
Consider a matrix polynomial $P(\lambda) = \lambda^d A_d + \dots + \lambda A_1 + I_n \in \mathbb{C}^{n, n}[\lambda]$. If the numerical range $W(P)$ is contained in some half-plane $H_\varphi$, $\varphi \in [0;2\pi)$, then
$$
\| P(\lambda) \| \leq 2\exp \left( \lambda_{H}\Bigl[\lambda A_1 -|\lambda|^2 A_2\Bigr] +\frac12 |\lambda|^2 \| A_1\|^2       \right), \quad \lambda \in \Comp.
$$
    \end{theorem}

\begin{proof}
Without loss of generality, we may assume that the numerical range of $P(\lambda)$ is contained in the open upper half-plane $H_0$. By the well-known fact that the operator two-norm of a matrix is less or equal to twice its numerical radius, see, e.g., \cite{HorJ91},  we obtain for $\lambda\in\Comp$
\begin{eqnarray*}
\|P(\lambda)\| &\leq & 2\max_{\|x\| = 1} |x^*P(\lambda)x| \\
& = &  2\max_{\|x\| = 1} |\lambda^d x^*A_dx + \dots + \lambda x^*A_1x + 1| \\
& \leq & 2\max_{\|x\|=1} \exp \left(\re(x^*A_1x \lambda) - \re(x^* A_2 x) |\lambda|^2  + \frac12 |x^*A_1x|^2|\lambda|^2  \right)=\\
& =&  2\max_{\|x\|=1} \exp \left(\re(  x^* ( \lambda  A_1  -  |\lambda|^2 A_2 )x )    + \frac12 \norm{A_1}^2|\lambda|^2  \right),
\end{eqnarray*}
where the second inequality follows from \eqref{dB} applied to $x^*P(\lambda)x$.
Taking the maximum under the exponent and applying \eqref{lambdaH} we see that the assertion follows.
\end{proof}

 As the assumption that the numerical range lies in a half-plane is rather stong, we have asked in \cite{szymanski2023stability} whether hyperstability implies any Sz\'asz-type inequality. 
 Apparently, this is not the case, as the following example shows.

 \begin{example}\rm
Consider a matrix polynomial of the form
$$
P(\lambda)=\begin{bmatrix}
    1 & \lambda^d \\ 0 & 1
\end{bmatrix},
$$
with $d \geq 3$. Then $P(0)=I_2$, and $P(\lambda)$ is hyperstable with respect to any set $D\subseteq \Comp$ (cf. Theorem~\ref{upperblock}). However,
$$
\norm{P(\lambda)} \geq \frac{1}{\sqrt 2} \norm{P(\lambda)}_F= \sqrt{1+\frac{|p(\lambda)|^2}{2}}.
$$
Hence, as $A_1 = A_2 = 0$, there can not be any global bound on  $\norm{P(\lambda)}$ of Sz\'asz type. 

The last thing we can do here is to compute the numerical range of $P(\lambda)$. Let $x \in \mathbb{C}^2\setminus\{0\}$. We have $x^*P(\lambda)x = |x_1|^2 + |x_2|^2 + \lambda^d\overline{x}_1 x_2$. In particular, for $x_1 = x_2 = 1$, we see that the numerical range of $P(\lambda)$ is not contained in any half-plane. 
 \end{example}

All this motivates us to provide several other inequalities of Sz\'asz type for matrix polynomials.

\section{Employing factorizations}

In this paragraph, matrix-type Sz\'asz inequalities are considered. First inequality of such type was presented in Proposition \ref{Szasz}. Below we sometimes use the Frobenius norm of a matrix
$$
\norm{A}_F := \sqrt{\tr(A^*A)}
$$
instead of the usual operator norm $\norm{A}=\sup_{\norm x=1}\norm{Ax}$.
We begin with two lemmas. 

    \begin{lemma}\label{mlog}
Let $A \in \mathbb{C}^{n,n} \setminus \{-n^{-1/2}I\}$. Then
            \begin{equation*}
\log\norm{\frac{1}{\sqrt{n}}I + A}_F \leq \frac{1}{\sqrt{n}}\tr\RE A + \frac{1}{2}||A||_F^2\text{.}
            \end{equation*}           
    \end{lemma}
        \begin{proof}
Using the inequality $\log(1 + x) \leq x$, $x \in (-1; \infty)$, one obtains
            \begin{eqnarray*}
\log\Big|\Big|\frac{1}{\sqrt{n}}I + A\Big|\Big|_F &=& \frac{1}{2}\log\Big|\Big|\frac{1}{\sqrt{n}}I + A\Big|\Big|_F^2 \\
&=& \frac{1}{2}\log\tr\Bigl[\Bigl(\frac{1}{\sqrt{n}}I + A\Bigr)^*\Bigl(\frac{1}{\sqrt{n}}I + A\Bigr)\Bigr] \\
& = &\frac{1}{2}\log\tr(\frac{1}{n}I + \frac{1}{\sqrt{n}}A + \frac{1}{\sqrt{n}}A^* + A^*A)\\
&=& \frac{1}{2}\log(1 + \frac{2}{\sqrt{n}}\tr\RE A + ||A||_F^2)  \\
& \leq & \frac{1}{2}(\frac{2}{\sqrt{n}}\tr\RE A + ||A||_F^2) \\
&=& \frac{1}{\sqrt{n}}\tr\RE A + \frac{1}{2}||A||_F^2\text{.}
            \end{eqnarray*}
        \end{proof}
The following lemma is a contribution of the Author of the Thesis.
    \begin{lemma}\label{sums}
Let $B_1, B_2, \dots , B_d$ be $n \times n$ complex matrices such that $\IM B_j \leq 0$ for $j \in \{1, 2, \dots , d\}$ and let $||\cdot||_F$ denotes the Frobenius norm of a matrix. Then:
        \begin{equation*}
\sum_{j=1}^d ||B_j||_F^2 \leq \big|\big|\sum_{j=1}^d B_j\big|\big|_F^2 - 2\tr\RE\sum_{1 \leq j < k \leq d}B_jB_k\text{.}
        \end{equation*}
    \end{lemma}
    \begin{proof}
Observe that
        \begin{eqnarray*}
\norm{\sum_{j=1}^d B_j}_F^2 &=& \tr\Bigl(\sum_{j=1}^d B_j^*\sum_{k=1}^d B_k\Bigr)\\
&=& \tr\Bigl(\sum_{j=1}^d B_j^*B_j + \sum_{1 \leq j \neq k \leq d}B_j^*B_k\Bigr)  \\
& =& \tr\sum_{j=1}^d B_j^*B_j + \tr\sum_{1 \leq j < k \leq d}(B_j^*B_k + B_k^*B_j) \\
&=& \sum_{j=1}^d\tr B_j^*B_j + \tr\sum_{1 \leq j < k \leq d}2\RE(B_j^*B_k) \\
& =& \sum_{j=1}^d ||B_j||_F^2 + 2\tr\RE\sum_{1 \leq j < k \leq d}B_j^*B_k\text{,}
        \end{eqnarray*}
therefore it is sufficient to prove the following statement
        \begin{equation*}
2\tr\RE\sum_{1 \leq j < k \leq d}(B_j^*B_k - B_jB_k) \geq 0\text{.}
        \end{equation*}
Now, note that
        \begin{eqnarray*}
2\tr\!&\RE&\!\sum_{1 \leq j < k \leq d}(B_j^*B_k - B_jB_k) \\
&=& 
\sum_{1 \leq j < k \leq d}\bigl(\tr(B_j^*B_k) - \tr(B_jB_k) + \tr(B_k^*B_j) - \tr(B_k^*B_j^*)\bigr) \\
& =& \sum_{1 \leq j < k \leq d}\bigl(\tr(B_j^*B_k) - \tr(B_jB_k) + \tr(B_jB_k^*) - \tr(B_j^*B_k^*)\bigr) \\
&=& \sum_{1 \leq j < k \leq d}\tr\bigl((B_j - B_j^*)(B_k^* - B_k)\bigr) \\
& =& 4\sum_{1 \leq j < k \leq d}\tr(\IM B_j\IM B_k)\text{,}
        \end{eqnarray*}
which is indeed a non-negative number since
        \begin{equation*}
\tr(\IM B_j\IM B_k) = \tr\bigl(\sqrt{-\IM B_j}(-\IM B_k)(\sqrt{-\IM B_j})^*\bigr) \geq 0
        \end{equation*} 
as the trace of $\sqrt{-\IM B_j}(-\IM B_k)(\sqrt{-\IM B_j})^* \geq 0$, which is positive semi-definite due to the assumption that $\IM B_j$ and $\IM B_k$ are both negative semi-definite matrices.
    \end{proof}

Now we able to show a matrix version of the Sz\'asz inequality without the assumption of the numerical range lying in a half-plane. Instead we assume a factorisation of the matrix polynomial. We refer the reader to \cite{LiR94}, Section 3, for a relation between factorisation and the location of the  numerical range.

\begin{proposition}\label{Frob}
Let $P(\lambda) = I + \sum_{j=1}^d \lambda^j A_j \in \mathbb{C}^{n, n}[\lambda]$ be a matrix polynomial with the following factorization $P(\lambda) = \prod_{j=1}^d(I + \lambda B_j)$, where $\IM B_j \leq 0$ for $j \in \{1, 2, \dots , d\}$. Then the Frobenius norm of the polynomial $P(\lambda)$ can be estimated from above:
        \begin{equation*}
\norm{P(\lambda)}_F \leq n^{d/2}\exp\Bigl(\frac{1}{n}\tr\RE(\lambda A_1) + \frac{1}{2n}(||A_1||_F^2 - 2\tr\RE A_2)|\lambda|^2\Bigr), \;\lambda \in \mathbb{C}\text{.}
        \end{equation*}
    \end{proposition}
    \begin{proof}
To obtain the desired inequality, we will apply Lemma \ref{mlog} and Lemma \ref{sums} consecutively. First, we use Lemma \ref{mlog}  $d$ times:
        \begin{eqnarray}\label{nowenowe}
\sum_{j=1}^d\log||I + \lambda B_j||_F & \!=\! & d\log\sqrt{n} + \sum_{j=1}^d\log||\frac{1}{\sqrt{n}}I + \lambda \tilde{B}_j||_F\\ &\!\leq\!& \log\sqrt{n^d} + \frac{1}{\sqrt{n}}\sum_{j=1}^d \tr\RE(\lambda \tilde{B}_j) + \frac{1}{2}\sum_{j=1}^d ||\lambda \tilde{B}_j||_F^2\text{,}
        \end{eqnarray}
where $\tilde{B}_j := n^{-1/2}B_j$ for $j \in \{1,2, \dots , d\}$.
Since $I + \sum_{j=1}^d \lambda^j A_j = \prod_{j=1}^d(I + \lambda B_j)$, we have $\sum_{j=1}^d \tilde{B}_j = n^{-1/2}A_1$ and $\sum_{1 \leq j < k \leq d}\tilde{B}_j\tilde{B}_k = n^{-1}A_2$. Then, the right side of the last inequality can be rewritten and estimated by Lemma~\ref{sums}  as follows:
        \begin{eqnarray*}
\log\sqrt{n^d} + \frac{1}{n}\tr\RE(\lambda A_1) + \frac{1}{2}|\lambda|^2\sum_{j=1}^d||\tilde{B}_j||_F^2 \\
\leq \log\sqrt{n^d} + \frac{1}{n}\tr\RE(\lambda A_1) + \frac{1}{2n}(||A_1||_F^2 - 2\tr\RE A_2)|\lambda|^2.
        \end{eqnarray*}
Therefore, we obtain
        \begin{equation*}
\log||P(\lambda)||_F \leq \log\sqrt{n^d} + \frac{1}{n}\tr\RE(\lambda A_1) + \frac{1}{2n}(||A_1||_F^2 - 2\tr\RE A_2)|\lambda|^2.
        \end{equation*}
Taking the exponent of both sides of the above inequality ends the proof.
    \end{proof}
Two next examples are contributions of the Author of the Thesis.

The following example shows that the Frobenius norm of a matrix polynomial cannot be estimated from above independent of both: the degree $d$ of the polynomial and the size $n$ of the matrix coefficients.
    \begin{example}\label{ones}\rm
Let us take
        $$
B_j = B =
        \begin{bmatrix}
1 & 1 & \dots & 1 \\
1 & 1 & \dots & 1 \\
\vdots & \vdots & \ddots & \vdots \\
1 & 1 & \dots & 1
        \end{bmatrix},\;\;\;
j \in \{1, \dots , d\},
$$
and consider a polynomial $P(\lambda) = \prod_{j=1}^d(I + \lambda B_j)$. Then due to the fact that $B^j = n^{j-1}B$, what can be easily proven by induction on $j$, we have
        \begin{equation*}
P(\lambda) = (I + \lambda B)^d = I + \sum_{j=1}^d \binom{d}{j}\lambda^{j}B^{j} = I + \sum_{j=1}^d \binom{d}{j}n^{j-1}\lambda^{j}B\text{.}
        \end{equation*}
Next, we have
        \begin{equation*}
\norm{P(\lambda)}_F^2 = n\Bigl|1 + \sum_{j=1}^d \binom{d}{j}n^{j-1}\lambda^j\Bigr|^2 + (n^2 - n)\Bigl|\sum_{j=1}^d \binom{d}{j}n^{j-1}\lambda^j\Bigr|^2\text{.}
        \end{equation*}
Taking $\lambda > 0$ we obtain
        \begin{eqnarray*}
||P(\lambda)||_F^2 &=& n + 2n\sum_{j=1}^d \binom{d}{j}n^{j-1}\lambda^j + n^2\Bigl(\sum_{j=1}^d \binom{d}{j}n^{j-1}\lambda^j\Bigr)^2 \\
&=& (n\lambda + 1)^{2d} + n - 1 \geq n^{2d}\lambda^{2d}\text{.}
        \end{eqnarray*}
        This lower bound shows that  fixing $d>0$ and assuming factorisation 
        $$
        P(\lambda)=\prod_{j=1}^d (I + \lambda B_j),\quad \IM B_j\leq 0
        $$
        is still too less for a global Sz\'asz-type bound independent from $n$. Similarly,  for $n>1$ fixed there cannot be any global Sz\'asz-type bound that is independent from $d$.
        Note that the factorisation above is a strong assumption, implying stability with respect to the upper half-plane $H_0$.
        
Let us compare this lower bound with the upper bound in Proposition \ref{Frob}. There, also for $\lambda>0$, we receive
        \begin{eqnarray*}
\norm{P(\lambda)}_F &\leq& n^{d/2}\exp\Bigl(\frac{1}{n}\tr\RE(\lambda A_1) + \frac{1}{2n}(\norm{A_1}_F^2 - 2\tr\RE A_2)|\lambda|^2\Bigr) \\
&=& n^{d/2}\exp\biggl[\frac{1}{n}\tr\RE(\lambda d B) + \frac{1}{2n}\biggl(\norm{dB}_F^2 - 2\tr\RE\Bigl(\frac{d(d-1)}{2}nB\Bigr)\biggr)|\lambda|^2\biggr] \\
&=&  n^{d/2}\exp\Bigl(d\lambda + \frac{1}{2n}\bigl(n^2d^2 - n^2d(d-1)\bigr)|\lambda|^2\Bigr)\\
&=& n^{d/2}\exp\Bigl(d\lambda + \frac{nd}{2}|\lambda|^2\Bigr)\text{.}
        \end{eqnarray*}
Finally, let us compute the numerical range of $P(\lambda)$. Let $x \in \mathbb{C}^n \setminus \{0\}$. We have
$$
x^*x = \sum_{k=1}^n |x_k|^2,\;\;\; x^*Bx = \Bigl|\sum_{k=1}^n x_k\Bigr|^2\text{.}
$$
Using these equalities, we can write
$$
x^*P(\lambda)x = x^*x + \sum_{j=1}^d \binom{d}{j}n^{j-1}\lambda^j x^*Bx = \sum_{k=1}^n |x_k|^2 + \Bigl|\sum_{k=1}^n x_k\Bigr|^2\sum_{j=1}^d \binom{d}{j}n^{j-1}\lambda^j\text{.}
$$
By the binomial formula, we obtain that
$$
x^*P(\lambda)x = \sum_{k=1}^n |x_k|^2 + \Bigl|\sum_{k=1}^n x_k\Bigr|^2\frac{(n\lambda + 1)^d - 1}{n}\text{.}
$$
Solving $x^*P(\lambda)x = 0$ with respect to $\lambda$, for $x$ with $\sum_{k=1}^n x_k \neq 0$,
we obtain
    \begin{equation}\label{roots}
(n\lambda + 1)^d = -n\sum_{k=1}^n |x_k|^2\Big/\Bigl|\sum_{k=1}^n x_k\Bigr|^2 + 1.
    \end{equation}
Hence, the solutions of $x^*P(\lambda)x = 0$ form a regular polygon with $d$ vertices, centred at the point $-1/n$. Now, we can choose $x$ such that an expression on the right-hand side of \eqref{roots} has an arbitrary large negative value, say $-R < 0$. Then distances between the centre of the polygon and its vertices are equal $\sqrt[d]{R}/n$, i.e. they can be arbitrary large. Therefore, if $d \geq 3$, then the numerical range $W(P)$ cannot be contained in any half-plane $H_{\varphi}$, where $\varphi \in [0; 2\pi)$.
    \end{example}

Below we present an example of a sequence of matrix polynomials $P_k(\lambda)$ of the same size $n\times n$ but increasing degree   (cf. Proof of Theorem 1.3 in \cite{Kne19} for the scalar version). 
We  will be able to compute the exact values of $\lim_{k\to\infty}\norm{P_k(\lambda)}_F$. 
For a fixed $n$ this limit is finite, while Theorem~\ref{Frob} gives us only a sequence of bounds on $\norm{P_k(\lambda)}_F$ increasing to infinity. 
    
  \begin{example}\label{cmv}\rm
Let $c_1, c_2$ be any real numbers such that $\gamma := (c_1^2 - 2c_2)/2 > 0$. Then choose $k_0 \in \mathbb{Z}_+$ such that for $k \geq k_0$ we have $d_k := \gamma - c_1^2/(2k) \geq 0$. Now, consider a polynomial sequence
        \begin{equation*}
p_k(\lambda) = \biggl(1 + \frac{c_1\lambda}{k}\biggr)^k\biggl(1 + \frac{\sqrt{d_k}\lambda}{\sqrt{k}}\biggr)^k\biggl(1 - \frac{\sqrt{d_k}\lambda}{\sqrt{k}}\biggr)^k\text{.}
        \end{equation*}
This is the sequence of stable with respect to the open upper half-plane $H_0$ polynomials with real coefficients, for which the following condition of convergence holds, cf. \cite{Kne19}:
        \begin{equation}\label{lim}
\lim_{k \to \infty}|p_k(\ii y)| = \exp(\gamma y^2)\; \text{for all}\; y \in \mathbb{R}\text{.}
        \end{equation}
Let us construct a matrix version of this scalar example. For this aim  we take two matrices of size $n \times n$:
        \begin{equation*}
C_1 =
            \begin{bmatrix}
1 & 1 & \dots & 1 \\
1 & 1 & \dots & 1 \\
\vdots & \vdots & \ddots & \vdots \\
1 & 1 & \dots & 1
            \end{bmatrix}
\;\;\;\text{and}\;\;\;
C_2 =
            \begin{bmatrix}
-1 & 0 & \dots & 0 \\
0 & -1 & \dots & 0 \\
\vdots & \vdots & \ddots & \vdots \\
0 & 0 & \dots & -1
            \end{bmatrix}\text{.}
        \end{equation*}
Then, keeping the analogy with ~\cite{Kne19}, we have:
        \begin{equation*}
\Gamma := (C_1^2 - 2C_2)/2 = I + \frac{n}{2}C_1
        \end{equation*}
and
        \begin{eqnarray*}
D_k := \Gamma - C_1^2/(2k) = I + \frac{n(k-1)}{2k}C_1\text{.}
        \end{eqnarray*}
Now, let us define
$$
P_k(\lambda) := \biggl(I + \frac{C_1\lambda}{k}\biggr)^k\biggl(I + \frac{\sqrt{D_k}\lambda}{\sqrt{k}}\biggr)^k\biggl(I - \frac{\sqrt{D_k}\lambda}{\sqrt{k}}\biggr)^k
$$
and observe it has the form as in Theorem \ref{Frob}.
Further, observe that 
            \begin{eqnarray*}
P_k(\ii y)
&=&\biggl(I + \frac{\ii y C_1}{k}\biggr)^k\biggl(I + \frac{y^2 D_k}{k}\biggr)^k \\
&=& \Bigl(\frac{y^2 + k}{k}\Bigr)^k\biggl(I + \sum_{j=1}^k \binom{k}{j}n^{j-1}\frac{(\ii y)^{j}}{k^j}C_1\biggr) \\
&&\cdot\ \biggl(I + \sum_{j=1}^k \binom{k}{j}n^{2j-1}\frac{y^{2j}(k-1)^j}{(2k)^j(y^2 + k)^j}C_1\biggr) \\
&=& \Bigl(\frac{y^2 + k}{k}\Bigr)^k\Biggl[I + \frac{1}{n}\biggl[\Bigl(\frac{ny\ii}{k} + 1\Bigr)^k - 1\biggr]C_1\Biggr] \\
&&\cdot\ \Biggl[I + \frac{1}{n}\biggl[\Bigl(\frac{n^2y^2(k-1)}{2k(y^2 + k)} + 1\Bigr)^k - 1\biggr]C_1\Biggr] = \\
&=& \Bigl(\frac{y^2 + k}{k}\Bigr)^k\Biggl[I + \frac{1}{n}\biggl[\Bigl(\frac{n^2y^2(k-1)}{2k(y^2 + k)} + 1\Bigr)^k - 1\biggr]C_1 \\
&& + \ \frac{1}{n}\biggl[\Bigl(\frac{ny\ii}{k} + 1\Bigr)^k - 1\biggr]C_1 + \frac{1}{n^2}\biggl[\Bigl(\frac{ny\ii}{k} + 1\Bigr)^k - 1\biggr] \\ &&\cdot\ \biggl[\Bigl(\frac{n^2y^2(k-1)}{2k(y^2 + k)} + 1\Bigr)^k - 1\biggr]C_1^2\Biggr]\text{.}
            \end{eqnarray*}
Note that $C_1^2 = nC_1$. Therefore, we have
        \begin{equation*}
P_k(\ii y) = \Bigl(\frac{y^2 + k}{k}\Bigr)^k\Biggl[I + \frac{1}{n}\biggl[\Bigl(\frac{ny\ii}{k} + 1\Bigr)^k\Bigl(\frac{n^2y^2(k-1)}{2k(y^2 + k)} + 1\Bigr)^k - 1\biggr]C_1\Biggr]\text{.}
        \end{equation*}  
Finally, we can calculate the Frobenius norm of $P_k(\ii y)$:
        \begin{eqnarray*}
\norm{P_k(\ii y)}_F^2 
&=& n\Bigl(\frac{y^2 + k}{k}\Bigr)^{2k}\ \biggl|1 + \frac{1}{n}\biggl[\Bigl(\frac{ny\ii}{k} + 1\Bigr)^k\Bigl(\frac{n^2y^2(k-1)}{2k(y^2 + k)} + 1\Bigr)^k - 1\biggr]\biggr|^2 \\
&& + \ \frac{n - 1}{n}\Bigl(\frac{y^2 + k}{k}\Bigr)^{2k}\biggl|\Bigl(\frac{ny\ii}{k} + 1\Bigr)^k\Bigl(\frac{n^2y^2(k-1)}{2k(y^2 + k)} + 1\Bigr)^k - 1\biggr|^2\text{.}
        \end{eqnarray*}
The last step is taking limit of $\norm{P_k(\ii y)}_F$ with $k \to \infty$. To do this, we make use of a classical property of exponent:
        \begin{equation*}
\lim_{n \to \infty}\bigl(1 + f(n)\bigr)^{g(n)} = \exp\bigl[\ \lim_{n \to \infty} f(n)g(n)\ \bigr]\text{,}
        \end{equation*}
which holds whenever $f(n) \to 0$ and $g(n) \to \infty$ with $n \to \infty$. Thus, we obtain
        \begin{eqnarray*}
\lim_{k \to \infty}\norm{P_k(\ii y)}_F
&=& \Bigl[ne^{2y^2}\Bigl|1 + \frac{1}{n}(e^{ny\ii}e^{n^2y^2/2} - 1)\Bigr|^2 + \frac{n - 1}{n}e^{2y^2}|e^{ny\ii}e^{n^2y^2/2} - 1|^2\Bigr]^{1/2}\text{.}
        \end{eqnarray*}
To finish these calculations, we need to find both squared moduli appearing in the last limit. We have
       \begin{eqnarray*}
\Bigl|1 + \frac{1}{n}(e^{ny\ii}e^{n^2y^2/2} - 1)\Bigr|^2 
&=& \Bigl[1 + \frac{1}{n}(e^{ny\ii}e^{n^2y^2/2} - 1)\Bigr]\Bigl[1 + \frac{1}{n}(e^{-ny\ii}e^{n^2y^2/2} - 1)\Bigr] \\
&=& 1 + \frac{1}{n}(e^{-ny\ii}e^{n^2y^2/2} - 1) + \frac{1}{n}(e^{ny\ii}e^{n^2y^2/2} - 1) \\
&&+\ \frac{1}{n}(e^{ny\ii}e^{n^2y^2/2} - 1)\frac{1}{n}(e^{-ny\ii}e^{n^2y^2/2} - 1) \\
&=& 1 + \frac{1}{n}(e^{-ny\ii}e^{n^2y^2/2} - 1) + \frac{1}{n}(e^{ny\ii}e^{n^2y^2/2} - 1) \\
&&+\ \frac{1}{n^2}e^{n^2y^2} - \frac{1}         {n^2}e^{ny\ii}e^{n^2y^2/2} - \frac{1}{n^2}e^{-ny\ii}e^{n^2y^2/2} + \frac{1}{n^2} \\
&=& \frac{1}{n^2}e^{n^2y^2} + \frac{n-1}{n^2}e^{ny\ii}e^{n^2y^2/2} + \frac{n-1}{n^2}e^{-ny\ii}e^{n^2y^2/2} \\
&&+\ \frac{n^2 - 2n + 1}{n^2}
        \end{eqnarray*}
and
        \begin{eqnarray*}
|e^{ny\ii}e^{n^2y^2/2} - 1|^2
&=& (e^{ny\ii}e^{n^2y^2/2} - 1)(e^{-ny\ii}e^{n^2y^2/2} - 1) \\
&=& e^{n^2y^2} - e^{ny\ii}e^{n^2y^2/2} - e^{-ny\ii}e^{n^2y^2/2} + 1\text{.}
        \end{eqnarray*}
Therefore, we see that
        \begin{equation}\label{mlim}
\lim_{k \to \infty}\norm{P_k(\ii y)}_F = e^{y^2}(e^{n^2y^2} + n - 1)^{1/2}\text{.}
        \end{equation}
Since for $n=1$ we have $C_1 = c_1 = 1,\; C_2 = c_2 = -1,\; \norm{P_k(\ii y)}_F =|p_k(\ii y)| $, we conclude from the equation \eqref{mlim} that $\lim_{k \to \infty} |p_k(\ii y)| = e^{3/2y^2}$, which is exactly the same value as we would obtain from the equation \eqref{lim} with $\gamma = 3/2$.
    \end{example}

    \begin{remark}\label{alt}\rm
 We present here another possible way of obtaining the Sz\'asz inequality, unfortunately providing essentially worse result than 
Theorem ~\ref{Frob}. For this aim, we need a slightly different version of Lemma ~\ref{mlog}.

Let $A \in \mathbb{C}^{n,n} \setminus \{-I\}$. Then
        \begin{equation*}
\log||I + A||_F \leq \tr\RE A + \frac{1}{2}||A||_F^2 + \frac{n-1}{2}\text{.}
        \end{equation*}

In the following estimations we use a standard inequality $\log(n + x) \leq x + n - 1$,   $x \in (-n; \infty)$, obtaining 
        \begin{eqnarray*}
\log||I + A||_F &=& \frac{1}{2}\log||I + A||_F^2 \\
&=& \frac{1}{2}\log\tr\bigl((I + A)^*(I + A)\bigr)  \\
& =& \frac{1}{2}\log\tr(I + A + A^* + A^*A)\\
&=& \frac{1}{2}\log(n + 2\tr\RE A + ||A||_F^2) \\
& \leq& \frac{1}{2}(2\tr\RE A + ||A||_F^2 + n - 1)\\
&=& \tr\RE A + \frac{1}{2}||A||_F^2 + \frac{n-1}{2}\text{.}
        \end{eqnarray*}
Now, we can formulate another Sz\'asz inequality mentioned above.

Let $P(\lambda) = I + \sum_{j=1}^d \lambda^j A_j \in \mathbb{C}^{n, n}[\lambda]$ be a matrix polynomial with the following factorization $P(\lambda) = \prod_{j=1}^d(I + \lambda B_j)$, where $\IM B_j \leq 0$ for $j \in \{1, 2, \dots , d\}$. Then the Frobenius norm of the polynomial $P(\lambda)$ can be estimated from above:
        \begin{equation*}
||P(\lambda)||_F \leq \exp\Bigl(\tr\RE(\lambda A_1) + \frac{1}{2}(||A_1||_F^2 - 2\tr\RE A_2)|\lambda|^2 + \frac{d(n-1)}{2}\Bigr), \;\lambda \in \mathbb{C}\text{.}
        \end{equation*}
The proof is similar as before
        \begin{eqnarray*}
\log||P(\lambda)||_F &\leq& \log\prod_{j=1}^d ||I + \lambda B_j||_F\\
&=& \sum_{j=1}^d\log||I + \lambda B_j||_F\\
&\leq& \sum_{j=1}^d \tr\RE(\lambda B_j) + \frac{1}{2}\sum_{j=1}^d ||\lambda B_j||_F^2 + \frac{d(n-1)}{2}.
        \end{eqnarray*}
Since 
$$
I + \sum_{j=1}^d \lambda^j A_j = \prod_{j=1}^d(I + \lambda B_j),
$$
we have $\sum_{j=1}^d B_j = A_1$ and $\sum_{1 \leq j < k \leq d}B_jB_k = A_2$. Continuing the estimation of $\log||P(\lambda)||_F$ we obtain by Lemma~\ref{sums} the following
      \begin{eqnarray*}
\log||P(\lambda)||_F &\leq& 
\tr\RE(\lambda A_1) + \frac{1}{2}|\lambda|^2\sum_{j=1}^d||B_j||_F^2 + \frac{d(n-1)}{2} \\
&\leq& \tr\RE(\lambda A_1) + \frac{1}{2}(||A_1||_F^2 - 2\tr\RE A_2)|\lambda|^2 + \frac{d(n-1)}{2}.
        \end{eqnarray*}
Therefore, 
        \begin{equation*}
\log||P(\lambda)||_F \leq \tr\RE(\lambda A_1) + \frac{1}{2}(||A_1||_F^2 - 2\tr\RE A_2)|\lambda|^2 + \frac{d(n-1)}{2},
        \end{equation*}
which is equivalent to the desired inequality.
    \end{remark}

\begin{remark} \rm Observe that both inequalities in Theorem \ref{Frob} and Remark \ref{alt} are based on the same method: estimation of $
\alpha\norm{ \alpha^{-1} I + \lambda (\alpha^{-1} B) }_F$, see e.g. \eqref{nowenowe}. In Lemma ~\ref{mlog} the parameter $\alpha$ was chosen as $\sqrt{n}$ while in Remark \ref{alt} $\alpha = 1$. One may seek for other values of $
\alpha$, that could give different bounds. Examples \ref{ones} and \ref{cmv} are, however, discouraging, in this matter.
\end{remark}

\section{Sz\'asz type inequalities for functional calculus}

In this Section we present one inequality for $\norm{p(A)}_F$ and two inequalities for $\norm{p(A)}$. In Proposition ~\ref{pA1} the estimation depends on both the size of an matrix $n$ and the degree $d$ of a polynomial $p(\lambda)$. In the next Proposition ~\ref{pA2} the estimation is only dependent on the parameter $d$, while the last Theorem ~\ref{SvN} contains the estimation independent on both parameters $n$ and $d$.
   
    \begin{proposition}\label{pA1}
Let $p(\lambda) = 1 + \sum_{j=1}^d a_j\lambda^j \in \mathbb{C}[\lambda]$ be a scalar polynomial, which is stable with respect to the upper half-plane $H_0$. Then the Frobenius norm of the polynomial $p(A)$ can be estimated from above as follows:
        \begin{equation*}
\norm{p(A)}_F \leq \sqrt{n^d}\exp\Bigl(\frac{1}{n}\tr\RE(a_1 A) + \frac{1}{2n}(|a_1|^2 - 2\re a_2)||A||_F^2\Bigr), \; A \in \mathbb{C}^{n,n}.
        \end{equation*}
    \end{proposition}
    \begin{proof}
Because of stability of polynomial $p(\lambda)$ we can write $p(\lambda) = \prod_{j=1}^d (1 + \alpha_j\lambda)$, where $\im\alpha_j \leq 0$ for $j \in \{1, 2, \dots , d\}$. Then obviously, we have $p(A) = \prod_{j=1}^d (I + \alpha_j A)$ and $\sum_{j=1}^d \alpha_j = a_1$, and $\sum_{1 \leq j<k \leq d} \alpha_j\alpha_k = a_2$. We will use both Lemma \ref{mlog} (exactly $d$ times) and Lemma \ref{sums} once again, but Lemma \ref{sums} will be used in a scalar version this time, i.e. for $n = 1$.
        \begin{eqnarray*}
\log||p(A)||_F 
&\leq& \log\prod_{j=1}^d ||I + \alpha_j A||_F \\
&=& \sum_{j=1}^d \log||I + \alpha_j A||_F \\
&=& d\log\sqrt{n} + \sum_{j=1}^d \log||\frac{1}{\sqrt{n}}I + \tilde{\alpha}_j A||_F \\
&\leq& \log\sqrt{n^d} + \frac{1}{\sqrt{n}}\sum_{j=1}^d \tr\RE(\tilde{\alpha}_j A) + \frac{1}{2}\sum_{j=1}^d ||\tilde{\alpha}_j A||_F^2\text{,}
        \end{eqnarray*}
where $\tilde{\alpha}_j := n^{-1/2}\alpha_j$ for $j \in \{1, 2, \dots, d\}$. Now, we obviously have $\sum_{j=1}^d \tilde{\alpha}_j = n^{-1/2}a_1$ and $\sum_{1 \leq j<k \leq d} \tilde{\alpha}_j\tilde{\alpha}_k = n^{-1}a_2$, so we can rewrite the right side of the last inequality and then, using the scalar version of Lemma \ref{sums} ($\im\tilde{\alpha}_j \leq 0$ for $j \in \{1, 2, \dots, d\}$), estimate it further as follows:
        \begin{align*}
&\log\sqrt{n^d} + \frac{1}{n}\tr\RE(a_1 A) + \frac{1}{2}||A||_F^2\sum_{j=1}^d |\tilde{\alpha}_j|^2 \\
&\leq \log\sqrt{n^d} + \frac{1}{n}\tr\RE(a_1 A) + \frac{1}{2n}(|a_1|^2 - 2\re a_2)||A||_F^2
        \end{align*}
Therefore, we obtain
        \begin{equation*}
\log||p(A)||_F \leq \log\sqrt{n^d} + \frac{1}{n}\tr\RE(a_1 A) + \frac{1}{2n}(|a_1|^2 - 2\re a_2)||A||_F^2
        \end{equation*}
Taking the exponent on the both side of the above inequality ends the proof.
    \end{proof}
In the next results, we move on to estimating the operator two-norm of $p(A)$. The proof of  Proposition ~\ref{pA2} is a modification of the  proof of Proposition ~\ref{pA1} above, but first we need a simple lemma.
    \begin{lemma}\label{imm}
Let $A \in \mathbb{C}^{n, n} \setminus \{-I\}$. Then
        \begin{equation*}
\log||I + A||\leq ||A||\text{.}
        \end{equation*}
    \end{lemma}
        \begin{proof}
The proof follows from the standard scalar inequality:
$\log\norm{I+A}\leq \log (1+\norm A)\leq \norm A$.
        \end{proof}
        \begin{proposition}\label{pA2}
Let $p(\lambda) = 1 + \sum_{j=1}^d a_j\lambda^j \in \mathbb{C}[\lambda]$ be a scalar polynomial, which is stable with respect to the open upper half-plane $H_0$. Then the operator two-norm of the polynomial $p(A)$ can be estimated from above as follows:
            \begin{equation*}
\norm{p(A)} \leq \exp\Bigl(||A||\sqrt{d(|a_1|^2 - 2\re a_2)}\Bigr), \; A \in \mathbb{C}^{n,n}.
            \end{equation*}
        \end{proposition}
        \begin{proof}
As before, because of stability of polynomial $p(\lambda)$ we can write $p(\lambda) = \prod_{j=1}^d (1 + \alpha_j\lambda)$, where $\im\alpha_j \leq 0$ for $j \in \{1, 2, \dots , d\}$. Then we have $p(A) = \prod_{j=1}^d (I + \alpha_j A)$ and $\sum_{j=1}^d \alpha_j = a_1$, and $\sum_{1 \leq j<k \leq d} \alpha_j\alpha_k = a_2$. Now, we use Lemma \ref{imm} (exactly $d$ times) and the scalar version of Lemma \ref{sums} consecutively.
            \begin{eqnarray*}
\log||p(A)||
&\leq& \log\prod_{j = 1}^d ||I + \alpha_j A|| \\
&=& \sum_{j = 1}^d\log ||I + \alpha_j A|| \\
&\leq& \sum_{j = 1}^d ||\alpha_j A|| \\
&=& ||A||\sum_{j = 1}^d |\alpha_j| \\
&\leq& ||A||\Bigl(d\sum_{j = 1}^d |\alpha_j|^2\Bigr)^{1/2} \\
&\leq& ||A||\sqrt{d(|a_1|^2 - 2\re a_2)}\text{,}
            \end{eqnarray*}
where the but last one inequality comes from the inequality of arithmetic and square means. Therefore, we have an inequality
            \begin{equation*}
\log||p(A)|| \leq ||A||\sqrt{d(|a_1|^2 - 2\re a_2)}\text{,}
            \end{equation*}
which is equivalent to the desired one after taking the exponent of its both sides.
        \end{proof}
        \begin{remark}\rm
From the proof of the previous proposition, we can see that for a stable with respect to the open upper half-plane polynomial $p(\lambda)$ with $p(0)=1$ the expression $|a_1|^2 - 2\re a_2 = \sum_{j = 1}^d |\alpha_j|^2 \geq 0$ has always a non-negative value.
        \end{remark}
As one may observe, in Proposition~\ref{pA1} the estimate contained $d$ and $n$ while in Proposition~\ref{pA2} the estimate contained $d$. We present now an estimate that is independent on both $d$ and $n$. Its proof is based on the von Neumann inequality
$$
\norm{p(A)}\leq \sup_{|z|\leq 1} |p(z)|,
$$
where $p(z)$ is a polynomial and $\norm{A}\leq 1$. For our purposes we will need the following form
\begin{equation}\label{vN}
 \norm{p(A)}\leq \sup_{|z|\leq \norm A} |p(z)|, 
\end{equation}
where $A$ is an arbitrary square matrix.  To see this let us take any $A\in \mathbb{C}^{n, n}$ and apply the classical inequality to the polynomial $\tilde p(z)=p(z\cdot\norm A)$ and the matrix $A/||A||$:
$$
\norm{p(A)} = \norm{\tilde p(A/\norm A)} \leq \sup_{|z| \leq 1}|\tilde p(z)| = \sup_{|z|\leq \norm A} |p(z)|\text{.}
$$
    \begin{theorem}\label{SvN}
Let $p(\lambda) \in \mathbb{C}[\lambda]$ be a stable with respect to the open upper half-plane $H_0$ polynomial with $p(0) = 1$ and let $A \in \mathbb{C}^{n, n}$ be a complex square matrix. Then the following inequality holds:
        \begin{equation*}
\norm{p(A)} \leq \exp\bigl(|a_1|\norm A  + \frac{1}{2} (|a_1|^2 -2 \re(a_2))\norm A^2\bigr)\text{.}
        \end{equation*}
    \end{theorem}
    \begin{proof}
As mentioned before, we use the general von Neumann inequality \eqref{vN} first. The next step is to use a scalar Sz\'asz inequality (cf.\,(\ref{dB}) in preliminaries). The last inequality is simply obtained by taking the supremum under the exponent.
\begin{eqnarray*}
    \norm{p(A)}&\leq& \sup_{|z|\leq \norm A} |p(z)| \\
    &\leq & \sup_{|z|\leq \norm A} \exp\bigl(\re (a_1z)  + \frac12 (|a_1|^2 -2 \re(a_2))|z|^2\bigr) \\
    &\leq & \exp\bigl(|a_1|\norm A  + \frac12 (|a_1|^2 -2 \re(a_2))\norm A^2\bigr)\text{.}
\end{eqnarray*}
This ends the proof.
    \end{proof}
Let us recall that for any matrix $A \in \mathbb{C}^{n, n}$ the relation $||A||_F \leq \sqrt{n}||A||$ between the operator two-norm and the Frobenius norm holds. Below, we make a comparison of global bounds presented in this Section.

The following example is a contribution of the Author of the Thesis.
    \begin{example}\label{comp}\rm
Using the relation between norms mentioned above, we can get estimations of $\norm{p(A)}_F$ from Proposition ~\ref{pA2} and Theorem ~\ref{SvN}. All inequalities, including one from Proposition ~\ref{pA1}, are listed here:
        \begin{eqnarray}
\label{tri}\hspace{1cm}\norm{p(A)}_F &\leq& \sqrt{n^d}\exp\Bigl(\frac{1}{n}\tr\RE(a_1 A) + \frac{1}{2n}(|a_1|^2 - 2\re a_2)\norm{A}_F^2\Bigr), \\
\label{sqrti}\hspace{1cm}\norm{p(A)}_F &\leq& \sqrt{n}\exp\Bigl(\norm{A}\sqrt{d(|a_1|^2 - 2\re a_2)}\Bigr), \\
\label{ui}\hspace{1cm}\norm{p(A)}_F &\leq& \sqrt{n}\exp\bigl(|a_1|\norm{A} + \frac{1}{2}(|a_1|^2 -2 \re a_2)\norm{A}^2\bigr)\text{.}
        \end{eqnarray}
Let us compare right-hand sides of these inequalities. Set $n=2$ and $d=3$. Our goal is to find three different pairs of values for $p(\lambda)$ and $A$ such that each of these inequalities gives the best upper-bound for $\norm{p(A)}_F$ in exactly one of the cases. Then, take $p(\lambda) = -(\lambda - 1)^3 = -\lambda^3 + 3\lambda^2 - 3\lambda + 1$, which is a stable polynomial with respect to the open upper half-plane $H_0$ with $p(0) = 1$. We have $a_1 = -3$ and $a_2 = 3$. For 
        \begin{equation*}
A =
            \begin{bmatrix}
1 & 1\\
1 & 1
            \end{bmatrix},\;\;\;
\norm{A}_F = 2,\;\;\; \norm{A} = 2
        \end{equation*}
right-hand sides of \eqref{tri}, \eqref{sqrti}, \eqref{ui} are equal $2\sqrt{2},\; e^6\sqrt{2},\; e^{12}\sqrt{2}$, respectively. Obviously, in this case the inequality \eqref{tri} gives the best estimation of $\norm{p(A)}_F$. 

Next, for
        \begin{equation*}
A =
            \begin{bmatrix}
-1 & 1\\
1 & -1
            \end{bmatrix},\;\;\;
\norm{A}_F = 2,\;\;\; \norm{A} = 2
        \end{equation*}
r.h.s. of \eqref{tri} - \eqref{ui} are equal $2e^6\sqrt{2},\; e^6\sqrt{2},\; e^{12}\sqrt{2}$, respectively. This time, second inequality \eqref{sqrti} does the best.

Now, put $p(\lambda) = (\lambda - 1)^2(\lambda + 1) = \lambda^3 - \lambda^2 - \lambda + 1$. Again, $p(\lambda)$ is stable with respect to the open upper half-plane $H_0$ and $p(0) = 1$. We see that $a_1 = a_2 = -1$. 

For 
        \begin{equation*}
A =
            \begin{bmatrix}
-1 & 0\\
0 & -1
            \end{bmatrix},\;\;\;
\norm{A}_F = \sqrt{2},\;\;\; \norm{A} = 1
        \end{equation*}
r.h.s. of \eqref{tri} - \eqref{ui} are equal $2e^{5/2}\sqrt{2},\; e^3\sqrt{2},\; e^{5/2}\sqrt{2}$, respectively. Clearly, in this case third inequality \eqref{ui} gives the best upper bound for $\norm{p(A)}_F$.
    \end{example}

\chapter[Hyperstability via stability of bivariate matrix polynomials]{Hyperstability of matrix polynomials in one variable via stability of bivariate matrix polynomials}\label{s5}

In this Chapter we will proceed in the following way. Given a one-variable matrix polynomial $P(\lambda)$ of a degree $2$ or $3$, we construct different bivariate matrix polynomials. We show that if one of those new polynomials is stable with respect to some $D^2$, then the original polynomial $P(\lambda)$ is hyperstable with respect to $D$.  
	
\section{Quadratic matrix polynomials}

We present now one of the major outcomes of the current Thesis. It provides several sufficient conditions for the hyperstability of a given univariate quadratic polynomials. The conditions are formulated using only the classical notion of the stability (not hyperstability!)  of multivariate matrix polynomials.
\begin{theorem}\label{poly2}
Let $P(\lambda) = \lambda^2 A_2 + \lambda A_1 + A_0$ be a  quadratic matrix polynomial and let $D$ be a nonempty open or closed subset of the complex plane $\mathbb{C}$. If at least one of the following conditions holds:
	\begin{enumerate}[\rm (a)]
\item\label{0?D} the multivariate matrix polynomial  $(z_1, z_2) \mapsto z_1^2 A_2 + z_2 A_1 + A_0$ is stable with respect to $D^2$,
\item\label{0notinD1} the  multivariate matrix polynomial $(z_1, z_2) \mapsto z_1z_2 A_2 + z_2 A_1 + A_0$  is stable with respect to $D^2$ and $0 \notin D$,
\item\label{0notinD2} the  multivariate matrix polynomial  $(z_1, z_2) \mapsto z_1^2z_2 A_2 + z_1^2 A_1 + z_2 A_0$ is stable with respect to $D^2$ and $0 \notin D$,
	\end{enumerate}
then the matrix polynomial $\lambda \mapsto P(\lambda)$ is hyperstable with respect to $D$.
\end{theorem}

\begin{proof}
First observe that setting $z_1 = z_2 = \lambda$ implies, in each case, that the polynomial $\lambda \mapsto P(\lambda)$ is stable with respect to $D$. In particular $\det P(\lambda) \neq 0$ on $D$, hence it is regular. In the case \eqref{0notinD2} the stability of a polynomial $\lambda \mapsto \lambda P(\lambda)$ is equivalent to the stability of the polynomial $\lambda \mapsto P(\lambda)$. To show that the polynomial $\lambda \mapsto P(\lambda)$ is hyperstable with respect to $D$, fix $x \in \Comp^n\setminus\{0\}$. Below by $\perp$ we denote orthogonality with respect to the standard complex inner-product. The proof in each case is similar, however, it requires certain adaptations and thus we present all the details.

%

Assume \eqref{0?D} holds. If there exists a vector $y \in \mathbb{C}^n\setminus\{0\}$ such that $y \perp A_2x, \; y \perp A_1x$ and $y \not\perp A_0x$, then 
	\begin{equation*}
y^* P(\lambda)x = y^*A_0x
	\end{equation*}
and \eqref{nozeroinD} is simply satisfied.

Hence, further we assume that such a vector $y$ does not exist, i.e.,
$
\{A_2x, A_1x\}^{\perp} \subseteq \{A_0x\}^{\perp}$,  equivalently $A_0x\in \Span\set{A_2x, A_1x}$.
Thus, we have 
\begin{equation}\label{A_0x}
A_0x = \alpha_0 A_2x + \beta_0 A_1x
\end{equation}
for some $\alpha_0, \beta_0 \in \Comp$. 
Next, we can write
	\begin{align*}
P(\lambda)x = \lambda^2 A_2x + \lambda A_1x + (\alpha_0 A_2x + \beta_0 A_1x) = (\lambda^2 + \alpha_0)A_2x + (\lambda + \beta_0)A_1x\text{.}
	\end{align*}
	
We show now that at least one of the scalar polynomials $\lambda \mapsto \lambda^2 + \alpha_0$ or $\lambda \mapsto \lambda + \beta_0$ is stable with respect to $D$.	
By the assumption that the polynomial $(z_1, z_2) \mapsto z_1^2 A_2 + z_2 A_1 + A_0$ is stable with respect to $D^2$, we have that the equation 
$$
(z_1^2 A_2 + z_2 A_1 + A_0)x = 0
$$
has no solutions $(z_1, z_2) \in D^2$. Substituting  \eqref{A_0x} we obtain that the equation
$$
(z_1^2 + \alpha_0)A_2x + (z_2 + \beta_0)A_1x = 0
$$
has no solutions $(z_1, z_2) \in D^2$. This implies that at least one of the scalar polynomials 
$\lambda \mapsto \lambda^2 + \alpha_0$ or $\lambda \mapsto \lambda + \beta_0$ has no roots in $D$ and the claim follows.
%

If the polynomial $\lambda \mapsto \lambda^2 + \alpha_0$ is stable with respect to $D$, then similarly as before we seek for a vector $y$ such that $y \not\perp A_2x$ and $y \perp A_1x$. If such a $y$ exists, we have
	\begin{equation*}
y^*P(\lambda)x = (y^*A_2x)(\lambda^2 + \alpha_0)
	\end{equation*}
and condition \eqref{nozeroinD} is satisfied.
If such a vector $y$ does not exist, we have $\{A_1x\}^{\perp} \subseteq \{A_2x\}^{\perp}$ and consequently $A_2x = c_1 A_1x$ for some constant $c_1 \in \mathbb{C}$. Then
	\begin{equation*}
P(\lambda)x = (c_1\lambda^2 + \lambda + c_1\alpha_0 + \beta_0)A_1x
	\end{equation*}
and since the polynomial $\lambda \mapsto P(\lambda)$ is regular and has no eigenvalues in $D$ we have that $A_1x \neq 0$ and  $\lambda \mapsto c_1\lambda^2 + \lambda + c_1\alpha_0 + \beta_0$ has no roots in $D$. Setting $y = A_1x$  we have
	\begin{equation*}
y^*P(\lambda)x = (A_1x)^*(c_1\lambda^2 + \lambda + c_1\alpha_0 + \beta_0)A_1x = \|A_1x\|^2(c_1\lambda^2 + \lambda + c_1\alpha_0 + \beta_0)
	\end{equation*}
and \eqref{nozeroinD} is again satisfied.

Similarly, if the polynomial $\lambda + \beta_0$ is stable with respect to $D$, then we seek for a vector $y$ such that $y \perp A_2x$ and $y \not\perp A_1x$. If such a $y$ exists, then
	\begin{equation*}
y^*P(\lambda)x = (y^*A_1x)(\lambda + \beta_0)
	\end{equation*}
and \eqref{nozeroinD} is satisfied.
If such the vector $y$ does not exist, we have $\{A_2x\}^{\perp} \subseteq \{A_1x\}^{\perp}$ and consequently $A_1x= c_2 A_2x$ for some constant $c_2 \in \mathbb{C}$ and 
	\begin{equation*}
P(\lambda)x =(\lambda^2 + c_2\lambda + \alpha_0 + c_2\beta_0)A_2x\text{.}
	\end{equation*}
As before, the scalar polynomial $\lambda \mapsto \lambda^2 + c_2\lambda + \alpha_0 + c_2\beta_0$ is stable and $A_2x \neq 0$. Taking $y = A_2x$ we have \eqref{nozeroinD} satisfied.

 Assume that \eqref{0notinD1} holds.
If there exists a vector $y \in \mathbb{C}^n\setminus\{0\}$ such that $y \perp A_2x, \; y \not\perp A_1x$ and $y \perp A_0x$, then 
	\begin{equation*}
y^* P(\lambda)x = (y^* A_1x)\lambda
	\end{equation*}
and \eqref{nozeroinD} is satisfied, due to the assumption that $0 \notin D$.

Hence, in the rest of  the proof of part (b) we assume that such a vector $y$ does not exist, i.e.,
$
\{A_2x, A_0x\}^{\perp} \subseteq \{A_1x\}^{\perp}$,  equivalently $A_1x\in \Span\set{A_2x, A_0x}$.
Thus, we have 
\begin{equation}\label{A_1x}
A_1x = \alpha_0 A_2x + \beta_0 A_0x
\end{equation}
for some $\alpha_0, \beta_0 \in \Comp$. 
Next, we can write
	\begin{align*}
P(\lambda)x = \lambda^2 A_2x + \lambda(\alpha_0 A_2x + \beta_0 A_0x) + A_0x = \lambda(\lambda + \alpha_0)A_2x + (\beta_0\lambda + 1)A_0x\text{.}
	\end{align*}
	
We show now that at least one of the scalar polynomials $\lambda \mapsto \lambda(\lambda + \alpha_0)$ or $\lambda \mapsto \beta_0\lambda + 1$ is stable with respect to $D$.	
By the assumption that the polynomial $(z_1, z_2) \mapsto z_1z_2 A_2 + z_2  A_1 + A_0$ is stable with respect to $D^2$ we have that the equation 
$$
(z_1z_2 A_2 + z_2  A_1 + A_0)x = 0
$$
has no solutions $(z_1, z_2) \in D^2$. Substituting \eqref{A_1x} we obtain that the equation
$$
z_2\left(z_1+\alpha_0\right) A_2 x+ \left( \beta_0 z_2 +1 \right)A_0x = 0
$$
has no solutions $(z_1, z_2) \in D^2$. This implies that at least one of the scalar polynomials 
$\lambda \mapsto \lambda(\lambda + \alpha_0)$ or $\lambda \mapsto \beta_0\lambda + 1$ has no roots in $D$ and the claim follows.
%

If the polynomial $\lambda \mapsto \lambda(\lambda + \alpha_0)$ is stable with respect to $D$, then similarly as before we seek for a vector $y$ such that $y \not\perp A_2x$ and $y \perp A_0x$. If such a $y$ exists, we have
	\begin{equation*}
y^*P(\lambda)x = (y^*A_2x)\lambda(\lambda + \alpha_0)
	\end{equation*}
and condition \eqref{nozeroinD} is satisfied.
If such the vector $y$ does not exist, we have $A_2x = c_1 A_0x$ for some constant $c_1 \in \mathbb{C}$. Then
	\begin{equation*}
P(\lambda)x = \Big[c_1 \lambda^2 + (c_1\alpha_0 + \beta_0)\lambda + 1\Big]A_0x
	\end{equation*}
and since the polynomial $\lambda \mapsto P(\lambda)$ is regular and has no eigenvalues in $D$ we have that $A_0x \neq 0$ and $c_1\lambda^2 + (c_1\alpha_0 + \beta_0)\lambda + 1$ has no roots in $D$. Setting $y = A_0x$  we have
	\begin{equation*}
y^*P(\lambda)x = (A_0x)^*\Big[c_2 \lambda^2 + (c_2\alpha_0 + \beta_0)\lambda + 1\Big]A_0x = \|A_0x\|^2\Big[c_2 \lambda^2 + (c_2\alpha_0 + \beta_0)\lambda + 1\Big] 
	\end{equation*}
and \eqref{nozeroinD} is again satisfied.

Similarly, if the polynomial $\lambda \mapsto \beta_0\lambda + 1$ is stable with respect to $D$, then we seek for a vector $y$ such that $y \perp A_2x$ and $y \not\perp A_0x$. If such a $y$ exists then
	\begin{equation*}
y^*P(\lambda)x = (y^*A_0x)(\beta_0\lambda + 1)
	\end{equation*}
and \eqref{nozeroinD} is satisfied.
If such a vector $y$ does not exist, we have $A_0x = c_2 A_2x$ for some constant $c_2 \in \mathbb{C}$ and 
	\begin{equation*}
P(\lambda)x = \Bigl[\lambda^2 + (\alpha_0 + c_2\beta_0 )\lambda + c_2\Bigr]A_2x\text{.}
	\end{equation*}
As before, the scalar polynomial $\lambda \mapsto \lambda^2 + (\alpha_0 + c_2\beta_0)\lambda + c_2$ is stable and $A_2x \neq 0$. Taking $y = A_2x$ we have \eqref{nozeroinD} satisfied.

Finally, assume that \eqref{0notinD2} holds. If there exists a vector $y \in \mathbb{C}^n\setminus\{0\}$ such that $y \not\perp A_2x, \; y \perp A_1x$ and $y \perp A_0x$, then 
	\begin{equation*}
y^* P(\lambda)x = (y^*A_2x)\lambda^2
	\end{equation*}
and \eqref{nozeroinD} is satisfied, due to the assumption that $0\notin D$.

Hence, further we assume that such a vector $y$ does not exist, i.e.,
$
\{A_1x, A_0x\}^{\perp} \subseteq \{A_2x\}^{\perp}$, equivalently $A_2x \in \Span\set{A_1x, A_0x}$.
Thus, we have
\begin{equation}\label{A_2x}
A_2x = \alpha_0 A_1x + \beta_0 A_0x
\end{equation}
for some $\alpha_0, \beta_0 \in \Comp$. 
Next, we can write
	\begin{align*}
P(\lambda)x = (\alpha_0\lambda^2 + \lambda)A_1x + (\beta_0\lambda^2 + 1)A_0x = \lambda(\alpha_0\lambda + 1)A_1x + (\beta_0\lambda^2 + 1)A_0x\text{.}
	\end{align*}
	
We show now that at least one of the scalar polynomials $\lambda \mapsto \lambda(\alpha_0\lambda + 1)$ or $\lambda \mapsto \beta_0\lambda^2 + 1$ is stable with respect to $D$.	
By the assumption that the polynomial $(z_1, z_2) \mapsto z_1^2z_2 A_2 + z_1^2 A_1 + z_2 A_0$ is stable with respect to $D^2$, we have that the equation 
$$
(z_1^2z_2 A_2 + z_1^2 A_1 + z_2 A_0)x = 0
$$
has no solutions $(z_1, z_2) \in D^2$. Substituting \eqref{A_2x} we obtain that the equation
$$
\label{noso}z_1^2(\alpha_0z_2 + 1)A_1x + z_2(\beta_0z_1^2 + 1)A_0x = 0
$$
has no solutions $(z_1, z_2) \in D^2$. This implies that at least one of the scalar polynomials 
$\lambda \mapsto \lambda(\alpha_0\lambda + 1)$ or $\lambda \mapsto \beta_0\lambda^2 + 1$ has no roots in $D$. Otherwise, there would exist $\lambda_1, \lambda_2 \in D$ such that $\lambda_1(\alpha_0\lambda_1 + 1) = \beta_0\lambda_2^2 + 1 = 0$. Since $\lambda_1 \neq 0$ we would have $\alpha_0\lambda_1 + 1 = \beta_0\lambda_2^2 + 1 = 0$ and taking $(z_1, z_2) = (\lambda_2, \lambda_1) \in D^2$, which is a solution of the equation \eqref{noso},  a contradiction. Therefore, the claim follows.
%

If the polynomial $\lambda \mapsto \lambda(\alpha_0\lambda + 1)$ is stable with respect to $D$, then similarly as before we seek for a vector $y$ such that $y \not\perp A_1x$ and $y \perp A_0x$. If such a $y$ exists, we have
	\begin{equation*}
y^*P(\lambda)x = (y^*A_1x)\lambda(\alpha_0\lambda + 1)
	\end{equation*}
and condition \eqref{nozeroinD} is satisfied.
If such vector $y$ does not exist, we have $\{A_0x\}^{\perp} \subseteq \{A_1x\}^{\perp}$ and consequently $A_1x = c_1 A_0x$ for some constant $c_1 \in \mathbb{C}$. Then:
	\begin{equation*}
P(\lambda)x = \Bigl[\lambda(\alpha_0\lambda + 1)c_1 + \beta_0\lambda^2 + 1\Bigr]A_0x = \Bigl[(\alpha_0c_1 + \beta_0)\lambda^2 + c_1\lambda + 1\Bigr]A_0x
	\end{equation*}
and since the polynomial $\lambda \mapsto P(\lambda)$ is regular and has no eigenvalues in $D$ we have that $A_0x \neq 0$ and the polynomial $\lambda \mapsto (\alpha_0c_1 + \beta_0)\lambda^2 + c_1\lambda + 1$ has no roots in $D$. Setting $y = A_0x$  we have
	\begin{equation*}
y^*P(\lambda)x = (A_0x)^*\Bigl[(\alpha_0c_1 + \beta_0)\lambda^2 + c_1\lambda + 1\Bigr]A_0x = \|A_0x\|^2\Bigl[(\alpha_0c_1 + \beta_0)\lambda^2 + c_1\lambda + 1\Bigr]
	\end{equation*}
and \eqref{nozeroinD} is satisfied again.

Similarly, if the polynomial $\lambda \mapsto \beta_0\lambda^2 + 1$ is stable with respect to $D$, then we seek for a vector $y$ such that $y \perp A_1x$ and $y \not\perp A_0x$. If such a $y$ exists, then
	\begin{equation*}
y^*P(\lambda)x = (y^*A_0x)(\beta_0\lambda^2 + 1)
	\end{equation*}
and \eqref{nozeroinD} is satisfied.
If such the vector $y$ does not exist, we have $\{A_1x\}^{\perp} \subseteq \{A_0x\}^{\perp}$ and consequently $A_0x= c_2 A_1x$ for some constant $c_2 \in \mathbb{C}$ and 
	\begin{equation*}
P(\lambda)x = \Bigl[\lambda(\alpha_0\lambda + 1) + (\beta_0\lambda^2 + 1)c_2\Bigr]A_1x = \Bigl[(\alpha_0 + \beta_0c_2)\lambda^2 + \lambda + c_2\Bigr]A_1x\text{.}
	\end{equation*}
As before, the scalar polynomial $\lambda \mapsto (\alpha_0 + \beta_0c_2)\lambda^2 + \lambda + c_2$ is stable with respect to $D$ and $A_1x \neq 0$. Taking $y = A_1x$ we have \eqref{nozeroinD} satisfied.
	\end{proof}

\section{Cubic matrix polynomials}

From Theorem~\ref{uppert}(ii) we know that the  hyperstability of palindromic matrix polynomials of degree three is equivalent to their stability. In addition to this, we consider now polynomials of the form $P(\lambda) = \lambda^3 A_0 + \lambda^2 A_2 + \lambda A_1 + A_0$. The following Theorem delivers us cubic hyperstable polynomials of this form.
    \begin{theorem}\label{poly3}
Let $P(\lambda) = \lambda^3 A_0 + \lambda^2 A_2 + \lambda A_1 + A_0$ be a cubic matrix polynomial and let $D$ be a nonempty open or closed subset of the complex plane $\mathbb{C}$. If at least one of the following conditions holds:
        \begin{enumerate}[\rm (a)]
\item\label{forA_0}
the multivariate matrix polynomial 
$$
(z_1, z_2) \mapsto (z_1^3 z_2^3 + z_1^3 + z_2^3)A_0 + (z_1^2 z_2^3 +z_1^2)A_2 + (z_1^3z_2 + z_2) A_1 + A_0
$$
is stable with respect to $D^2$ and $-1, \frac{1}{2} - \frac{\sqrt{3}}{2}i, \frac{1}{2} + \frac{\sqrt{3}}{2}i \not\in D$,
\item\label{forA_1}
the multivariate matrix polynomial 
$$
(z_1, z_2) \mapsto z_2^3 A_0 + z_1 z_2 A_2 + z_2 A_1 + A_0
$$ 
is stable with respect to $D^2$ and $0 \not\in D$,
\item\label{forA_2}
the multivariate matrix polynomial 
$$
(z_1, z_2) \mapsto z_1 z_2^3 A_0 + z_1 z_2^2 A_2 + z_2^2 A_1 + z_1 A_0
$$
is stable with respect to $D^2$ and $0 \not\in D$,
        \end{enumerate}
then the matrix polynomial $\lambda \mapsto P(\lambda)$ is hyperstable with respect to $D$.
    \end{theorem}
    \begin{proof}
In this proof, we proceed as in the proof of Theorem 5.2. Setting $z_1 = z_2 = \lambda$ leads in each case to the stability of the polynomial $\lambda \mapsto P(\lambda)$ with respect to $D$. In particular, $P(\lambda) $ is regular.  In the case \eqref{forA_0} the stability of a polynomial $\lambda \mapsto (\lambda^3 + 1)P(\lambda)$ is equivalent to the stability of the polynomial $\lambda \mapsto P(\lambda)$. For the  hyperstability of the polynomial $\lambda \mapsto P(\lambda)$ with respect to $D$, fix $x \in \mathbb{C}^n\setminus\{0\}$.

Assume that  condition \eqref{forA_0} holds. If there exists a vector $y \in \mathbb{C}^n\setminus\{0\}$ such that $y \not\perp A_0x, y \perp A_2x, y \perp A_1x$, then
    \begin{equation}
y^*P(\lambda)x = (y^*A_0x)(\lambda^3 + 1)
    \end{equation}
and the condition \eqref{nozeroinD} is satisfied because of the fact that $-1, \frac{1}{2} - \frac{\sqrt{3}}{2}i, \frac{1}{2} + \frac{\sqrt{3}}{2}i \not\in D$. Otherwise, we have $\{A_2x, A_1x\}^{\perp} \subseteq \{A_0x\}^{\perp}$, which equivalently means that 
$$
A_0x = \alpha_0 A_2x + \beta_0 A_1x
$$
for some $\alpha_0, \beta_0 \in \mathbb{C}$. Thus we can write
    \begin{equation*}
P(\lambda)x = (\alpha_0\lambda^3 + \lambda^2 + \alpha_0)A_2x + (\beta_0\lambda^3 + \lambda + \beta_0)A_1x\text{.}
    \end{equation*}
Due to the stability of the multivariate polynomial in \eqref{forA_0} we conclude that the equation
    \begin{equation*}
\Bigl[(z_1^3 z_2^3 + z_1^3 + z_2^3)A_0 + (z_1^2 z_2^3 +z_1^2)A_2 + (z_1^3z_2 + z_2) A_1 + A_0\Bigr]x = 0
    \end{equation*}
does not have solutions $(z_1, z_2) \in D^2$. Then the equation
    \begin{eqnarray*}
&&(\alpha_0 z_1^3 z_2^3 + z_1^2 z_2^3 + \alpha_0 z_1^3 + \alpha_0 z_2^3 + z_1^2 + \alpha_0)A_2x \\
&+& (\beta_0 z_1^3 z_2^3 + z_1^3 z_2 + \beta_0 z_1^3 + \beta_0 z_2^3 + z_2 + \beta_0)A_1x \\
&=& 0
    \end{eqnarray*}
does not have such solutions as well. Therefore, at least one of the polynomials $\lambda \mapsto \alpha_0\lambda^3 + \lambda^2 + \alpha_0$ or $\lambda \mapsto \beta_0\lambda^3 + \lambda + \beta_0$ is stable with respect to $D$. If the polynomial $\lambda \mapsto \alpha_0\lambda^3 + \lambda^2 + \alpha_0$ is stable with respect to $D$, then we seek for a vector $y \not\perp A_2x$ and $y \perp A_1x$. If such a vector $y$ exists, then
    \begin{equation*}
y^*P(\lambda)x = (y^*A_2x)(\alpha_0\lambda^3 + \lambda^2 + \alpha_0)
    \end{equation*}
and the condition \eqref{nozeroinD} is also satisfied. Otherwise, we have $\{A_1x\}^{\perp} \subseteq \{A_2x\}^{\perp}$, which equivalently means that $A_2x = c_1A_1x$ for some $c_1 \in \mathbb{C}$. Finally, we can write
    \begin{equation*}
P(\lambda)x = \Bigl[(c_1\alpha_0 + \beta_0)\lambda^3 + c_1\lambda^2 + \lambda + (c_1\alpha_0 + \beta_0)\Bigr]A_1x,
    \end{equation*}
where a polynomial in the square brackets is stable because of the stability of the polynomial $\lambda \mapsto P(\lambda)$. Now, we just take $y = A_1x$ to satisfy the condition \eqref{nozeroinD}:
    \begin{equation*}
y^*P(\mu)x = ||A_1x||^2\Bigl[(c_1\alpha_0 + \beta_0)\mu^3 + c_1\mu^2 + \mu + (c_1\alpha_0 + \beta_0)\Bigr] \neq 0
    \end{equation*}
for $\mu \in D$. If the polynomial $\lambda \mapsto \beta_0\lambda^3 + \lambda + \beta_0$ is stable with respect to $D$, then we seek in turn for a vector $y \perp A_2x$ and $y \not\perp A_1x$. If such a vector $y$ exists, then
    \begin{equation*}
y^*P(\lambda)x = (y^*A_1x)(\beta_0\lambda^3 + \lambda + \beta_0)
    \end{equation*}
and the condition \eqref{nozeroinD} is satisfied again. Otherwise, we have $\{A_2x\}^{\perp} \subseteq \{A_1x\}^{\perp}$, which equivalently means that $A_1x = c_2A_2x$ for some $c_2 \in \mathbb{C}$. As before, we can write:
\begin{equation*}
P(\lambda)x = \Bigl[(\alpha_0 + c_2\beta_0)\lambda^3 + \lambda^2 + c_2\lambda + (\alpha_0 + c_2\beta_0)\Bigr]A_2x,
    \end{equation*}
where a polynomial in the square brackets is stable because of the stability of the polynomial $\lambda \mapsto P(\lambda)$. As previously, we take $y = A_2x$ to satisfy the condition \eqref{nozeroinD}:
    \begin{equation*}
y^*P(\mu)x = ||A_2x||^2\Bigl[(\alpha_0 + c_2\beta_0)\lambda^3 + \lambda^2 + c_2\lambda + (\alpha_0 + c_2\beta_0)\Bigr] \neq 0
    \end{equation*}
for $\mu \in D$, which ends the proof in this case.

Assume that the condition \eqref{forA_1} holds. If there exists a vector $y \in \mathbb{C}^n\setminus\{0\}$ such that $y \perp A_0x, y \perp A_2x, y \not\perp A_1x$, then
    \begin{equation}
y^*P(\lambda)x = (y^*A_1x)\lambda
    \end{equation}
and the condition\eqref{nozeroinD} is satisfied because of the fact that $0 \not\in D$. Otherwise, we have $\{A_0x, A_2x\}^{\perp} \subseteq \{A_1x\}^{\perp}$, which equivalently means that 
$$
A_1x = \alpha_0 A_2x + \beta_0 A_0x
$$ 
for some $\alpha_0, \beta_0 \in \mathbb{C}$. Thus we can write
    \begin{equation*}
P(\lambda)x = (\lambda^2 + \alpha_0\lambda)A_2x + (\lambda^3 +\beta_0\lambda +1)A_0x\text{.}
    \end{equation*}
Due to the stability of the multivariate polynomial in \eqref{forA_1} we conclude that the equation
    \begin{equation*}
\Bigl[z_2^3 A_0 + z_1 z_2 A_2 + z_2 A_1 + A_0\Bigr]x = 0
    \end{equation*}
does not have solutions $(z_1, z_2) \in D^2$. Then the equation
    \begin{equation*}
(z_1 z_2 + \alpha_0 z_2)A_2x + (z_2^3 + \beta_0z_2 + 1)A_0x = 0
    \end{equation*}
does not have such solutions as well. Therefore, at least one of the polynomials $\lambda \mapsto \lambda^2 + \alpha_0\lambda$ or $\lambda \mapsto \lambda^3 + \beta_0\lambda + 1$ is stable with respect to $D$. If the polynomial $\lambda \mapsto \lambda^2 + \alpha_0\lambda$ is stable with respect to $D$, then we seek for a vector $y \not\perp A_2x$ and $y \perp A_0x$. If such vector $y$ exists, then
    \begin{equation*}
y^*P(\lambda)x = (y^*A_2x)(\lambda^2 + \alpha_0\lambda)
    \end{equation*}
and the condition \eqref{nozeroinD} is also satisfied. Otherwise, we have $\{A_0x\}^{\perp} \subseteq \{A_2x\}^{\perp}$, which equivalently means that $A_2x = c_1A_0x$ for some $c_1 \in \mathbb{C}$. Finally, we can write
    \begin{equation*}
P(\lambda)x = \Bigl[\lambda^3 + c_1\lambda^2 + (c_1\alpha_0 + \beta_0)\lambda + 1\Bigr]A_0x,
    \end{equation*}
where a polynomial in the square brackets is stable because of stability of the polynomial $\lambda \mapsto P(\lambda)$. Now, we just take $y = A_0x$ to satisfy the condition \eqref{nozeroinD}:
    \begin{equation*}
y^*P(\mu)x = ||A_0x||^2\Bigl[\mu^3 + c_1\mu^2 + (c_1\alpha_0 + \beta_0)\mu + 1\Bigr] \neq 0
    \end{equation*}
for $\mu \in D$. If the polynomial $\lambda \mapsto \lambda^3 + \beta_0\lambda + 1$ is stable with respect $D$, then we seek in turn for a vector $y \perp A_2x$ and $y \not\perp A_0x$. If such vector $y$ exists, then
    \begin{equation*}
y^*P(\lambda)x = (y^*A_0x)(\lambda^3 + \beta_0\lambda + 1)
    \end{equation*}
and the condition \eqref{nozeroinD} is satisfied again. Otherwise, we have $\{A_2x\}^{\perp} \subseteq \{A_0x\}^{\perp}$, which equivalently means that $A_0x = c_2A_2x$ for some $c_2 \in \mathbb{C}$. As before, we can write:
\begin{equation*}
P(\lambda)x = \Bigl[c_2\lambda^3 + \lambda^2 + (\alpha_0 + c_2\beta_0)\lambda + c_2\Bigr]A_2x,
    \end{equation*}
where the polynomial in the square brackets is stable because of the stability of the polynomial $\lambda \mapsto P(\lambda)$. As previously, we take $y = A_2x$ to satisfy the condition \eqref{nozeroinD}:
    \begin{equation*}
y^*P(\mu)x = ||A_2x||^2\Bigl[c_2\mu^3 + \mu^2 + (\alpha_0 + c_2\beta_0)\mu + c_2\Bigr] \neq 0
    \end{equation*}
for $\mu \in D$, which ends the proof in this case.

Assume that the condition \eqref{forA_2} holds. If there exists a vector $y \in \mathbb{C}^n\setminus\{0\}$ such that $y \perp A_0x, y \not\perp A_2x, y \perp A_1x$, then
    \begin{equation}
y^*P(\lambda)x = (y^*A_2x)\lambda^2
    \end{equation}
and the condition\eqref{nozeroinD} is satisfied because of the fact that $0 \not\in D$. Otherwise, we have $\{A_1x, A_0x\}^{\perp} \subseteq \{A_2x\}^{\perp}$, which equivalently means that 
$$
A_2x = \alpha_0 A_1x + \beta_0 A_0x
$$ 
for some $\alpha_0, \beta_0 \in \mathbb{C}$. Thus we can write
    \begin{equation*}
P(\lambda)x = (\alpha_0\lambda^2 + \lambda)A_1x + (\lambda^3 + \beta_0\lambda^2 + 1)A_0x\text{.}
    \end{equation*}
Due to the stability of the multivariate polynomial in \eqref{forA_2} we conclude that the equation
    \begin{equation*}
\Bigl[z_1 z_2^3 A_0 + z_1 z_2^2 A_2 + z_2^2 A_1 + z_1 A_0\Bigr]x = 0
    \end{equation*}
does not have solutions $(z_1, z_2) \in D^2$. Then the equation
    \begin{equation*}
(\alpha_0 z_1 z_2^2 + z_2^2)A_1x + (z_1 z_2^3 + \beta_0 z_1 z_2^2 + z_1)A_0x = 0
    \end{equation*}
does not have such solutions as well. Therefore, at least one of the polynomials $\lambda \mapsto \alpha_0\lambda^2 + \lambda$ or $\lambda \mapsto \lambda^3 + \beta_0\lambda^2 + 1$ is stable with respect to $D$. If the polynomial $\lambda \mapsto \alpha_0\lambda^2 + \lambda$ is stable with respect to $D$, then we seek for a vector $y \not\perp A_1x$ and $y \perp A_0x$. If such a vector $y$ exists, then
    \begin{equation*}
y^*P(\lambda)x = (y^*A_1x)(\alpha_0\lambda^2 + \lambda)
    \end{equation*}
and the condition \eqref{nozeroinD} is also satisfied. Otherwise, we have $\{A_0x\}^{\perp} \subseteq \{A_1x\}^{\perp}$, which equivalently means that $A_1x = c_1A_0x$ for some $c_1 \in \mathbb{C}$. Finally, we can write
    \begin{equation*}
P(\lambda)x = \Bigl[\lambda^3 + (c_1\alpha_0 + \beta_0)\lambda^2 + c_1\lambda + 1\Bigr]A_0x,
    \end{equation*}
where the polynomial in the square brackets is stable because of the stability of the polynomial $\lambda \mapsto P(\lambda)$. Now, we just take $y = A_0x$ to satisfy the condition \eqref{nozeroinD}:
    \begin{equation*}
y^*P(\mu)x = ||A_0x||^2\Bigl[\mu^3 + (c_1\alpha_0 + \beta_0)\mu^2 + c_1\mu + 1\Bigr] \neq 0
    \end{equation*}
for $\mu \in D$. If the polynomial $\lambda \mapsto \lambda^3 + \beta_0\lambda^2 + 1$ is stable with respect to $D$, then we seek in turn for a vector $y \perp A_1x$ and $y \not\perp A_0x$. If such a vector $y$ exists, then
    \begin{equation*}
y^*P(\lambda)x = (y^*A_0x)(\lambda^3 + \beta_0\lambda^2 + 1)
    \end{equation*}
and the condition \eqref{nozeroinD} is satisfied again. Otherwise, we have $\{A_1x\}^{\perp} \subseteq \{A_0x\}^{\perp}$, which equivalently means that $A_0x = c_2A_1x$ for some $c_2 \in \mathbb{C}$. As before, we can write:
\begin{equation*}
P(\lambda)x = \Bigl[c_2\lambda^3 + (\alpha_0 + c_2\beta_0)\lambda^2 + \lambda + c_2\Bigr]A_1x,
    \end{equation*}
where a polynomial in the square brackets is stable because of stability of the polynomial $\lambda \mapsto P(\lambda)$. As previously, we take $y = A_1x$ to satisfy the condition \eqref{nozeroinD}:
    \begin{equation*}
y^*P(\mu)x = ||A_1x||^2\Bigl[c_2\mu^3 + (\alpha_0 + c_2\beta_0)\mu^2 + \mu + c_2\Bigr] \neq 0
    \end{equation*}
for $\mu \in D$, which ends the proof in this case and the whole proof is completed.
\end{proof}

\chapter{Operators preserving matrix hyperstability}\label{sPol}

In this Chapter we develop the theory  of hyperstable matrix polynomials of several variables. 
This, besides an independent interest,  will later on  serve as a tool for showing the hyperstability of univariate matrix polynomials. 

\section{Basic operations}

Let us list some basic operations on multi-variable matrix polynomials which are hyperstable with respect to $\kappa$-th Cartesian power of an open half-plane with the boundary containing the origin. Below 
$$
H_{\varphi}^{\kappa} = \bigl\{(z_1, z_2, \dots , z_{\kappa}) \in \mathbb{C}^{\kappa} : \im(z_j e^{\ii\varphi}) > 0 \;\text{for all}\; j \in \{1, 2, \dots , \kappa\}\bigr\}
$$ 
be the $\kappa$-th Cartesian power of the open half-plane $H_{\varphi}$.
   \begin{proposition}\label{bop}
Let $\varphi \in [0; 2\pi)$ be fixed. Then the following linear transformations acting on the complex matrix polynomial space $\mathbb{C}^{n, n}[z_1, z_2, \dots , z_{\kappa}]$ map every hyperstable with respect to a product $H_{\varphi}^{\kappa}$ matrix polynomial to hyperstable with respect to $H_{\varphi}^{\kappa}$ matrix polynomial:
        \begin{enumerate}[\rm (i)]
\item{Permutation:}
$$
P(z_1, z_2, \dots , z_{\kappa}) \mapsto P(z_{\sigma(1)}, z_{\sigma(2)}, \dots , z_{\sigma(\kappa)})
$$
for every permutation $\sigma \in S_n$. \\
\item{Scaling:}
$$
P(z_1, z_2, \dots , z_{\kappa}) \mapsto P(z_1, z_2, \dots , z_{j-1}, az_j, z_{j+1}, \dots , z_{\kappa})
$$
for any $a \in \mathbb{R}_+$. \\
\item{Diagonalization:}
$$
P(z_1, z_2, \dots , z_{\kappa}) \mapsto P(z_j, z_j, \dots , z_j, z_{j+1}, \dots , z_{\kappa}) \in \mathbb{C}^{n, n}[z_j, z_{j+1}, \dots , z_{\kappa}]
$$
for every $j \in \{1, 2, \dots , \kappa\}$. \\
\item{Inversion (with rotation):}
$$
P(z_1, z_2, \dots , z_{\kappa}) \mapsto z_j^d P(z_1, z_2, \dots , z_{j-1}, -e^{-2\ii\varphi}/z_j, z_{j+1}, \dots , z_{\kappa}), 
$$
where the power $d := \deg_j P$ is the degree of a polynomial $P$ with respect to the variable $z_j$. \\
\item{Specialization:}
$$
P(z) \mapsto P(z_1, z_2, \dots , z_{j-1}, a, z_{j+1}, \dots , z_n) \in \mathbb{C}[z_1, z_2, \dots , z_{j-1}, z_{j+1}, \dots , z_n]
$$
for any $a \in {H}_{\varphi}$, where $z = (z_1, z_2, \dots , z_n)$.
     \end{enumerate}
Furthermore, differentiation 
$$
P(z_1, z_2, \dots , z_{\kappa}) \mapsto \frac{\partial P}{\partial z_j}(z_1, z_2, \dots , z_{\kappa}) \in \mathbb{C}^{n, n}[z_1, z_2, \dots , z_{\kappa}],\; j \in \{1, 2, \dots , \kappa\}
$$
also preserves matrix hyperstability, provided that the entries of the partial derivative $(\partial P)/(\partial z_j)(z_1, z_2, \dots , z_{\kappa})$ are linearly independent polynomials.
    \end{proposition}
    \begin{proof}
For the proof of first property \rm{(1)} we use an analogous property for scalar polynomials (cf. Proposition~\ref{op1}, \rm{(1)}). Let us fix $x \in \mathbb{C}^n \setminus \{0\}$ and note that hyperstability (with respect to $H_{\varphi}^{\kappa}$) of the matrix polynomial $P(z_1, z_2, \dots , z_{\kappa})$ means that  there exists $y \in \mathbb{C}^n $ such that the scalar polynomial 
$$
p(z_1, z_2, \dots , z_{\kappa}) = y^*P(z_1, z_2, \dots , z_{\kappa})x
$$ 
is stable with respect to $H_{\varphi}^{\kappa}$. 
It follows from Proposition~\ref{op1}, \rm{(1)} that the  polynomial $p(z_{\sigma(1)}, z_{\sigma(2)}, \dots , z_{\sigma(\kappa)})$ is  stable with respect to $H_{\varphi}^{\kappa}$. Since 
$$
p(z_{\sigma(1)}, z_{\sigma(2)}, \dots , z_{\sigma(\kappa)}) = y^*P(z_{\sigma(1)}, z_{\sigma(2)}, \dots , z_{\sigma(\kappa)})x\text{,}
$$
we see that for the matrix polynomial $P(z_{\sigma(1)}, z_{\sigma(2)}, \dots , z_{\sigma(\kappa)})$ the  vector $y$ satisfies the definition of hyperstability, i.e., $P(z_{\sigma(1)}, z_{\sigma(2)}, \dots , z_{\sigma(\kappa)})$ is also hyperstable with respect to $H_{\varphi}^{\kappa}$. 

Similarly, properties \rm{(2)-(5)} follow from the corresponding properties for scalar polynomials (cf. ~\ref{op1}, \rm{(2)-(5)}) via the equations:
        \begin{eqnarray*}
&&p(z_1, z_2, \dots , z_{j-1}, az_j, z_{j+1}, \dots , z_{\kappa}) = y^*P(z_1, z_2, \dots , z_{j-1}, az_j, z_{j+1}, \dots , z_{\kappa})x\text{,} \\
&&p(z_j, z_j, \dots , z_j, z_{j+1}, \dots , z_{\kappa}) = y^*P(z_j, z_j, \dots , z_j, z_{j+1}, \dots , z_{\kappa})x\text{,} \\
&&z_j^d p(z_1, z_2, \dots , -e^{-2\ii\varphi}/z_j, \dots , z_{\kappa}) = y^*\bigl[z_j^d P(z_1, z_2, \dots , -e^{-2\ii\varphi}/z_j, \dots , z_{\kappa})\bigr]x\text{,} \\
&&p(z_1, z_2, \dots , z_{j-1}, a, z_{j+1}, \dots , z_{\kappa}) = y^*P(z_1, z_2, \dots , z_{j-1}, a, z_{j+1}, \dots , z_{\kappa})x\text{,}
        \end{eqnarray*}
respectively. For a proof of the last property \rm{(6)}, it is sufficient to notice that the set $\mathbb{C}^{\kappa} \setminus H_{\varphi}^{\kappa}$ is separately convex with respect to each variable. Then the assertion follows from our matrix version of the Gauss-Lucas Theorem, see Theorem~\ref{mmgl}. This ends the  proof.
    \end{proof}
    \begin{example}\rm
Let us notice that, unlike in the scalar case,  specialization
$$
P(z_1, z_2, \dots , z_{\kappa}) \mapsto P(z_1, z_2, \dots , z_{j-1}, a, z_{j+1}, \dots , z_{\kappa})
$$
may not preserve hyperstability  for $a \in \partial {H}_{\varphi}$. Indeed, take a bivariate polynomial, hyperstable with respect to the product $H_0^2$
        \begin{equation*}
P(z_1, z_2) = 
            \begin{bmatrix}
z_1 & 0 \\
0 & z_2
            \end{bmatrix}
        \end{equation*}
and insert $z_2 = a = 0$. The obtained matrix polynomial
        \begin{equation*}
P(z_1, 0) = 
            \begin{bmatrix}
z_1 & 0 \\
0 & 0
            \end{bmatrix}
        \end{equation*}
is not hyperstable on any set, since for $x = [0\; 1]^{\top}$ we have $y^*P(z_1, 0)x = 0$ for all $y \in \mathbb{C}^2$.
    \end{example}
Similarly, an analogue of Theorem ~\ref{limiting} for matrix polynomials is not true, as the following example shows.
    \begin{example}\rm
Let us consider a sequence of matrix polynomials
        \begin{equation*}
P_k(z_1, z_2) = 
            \begin{bmatrix}
z_1 & 0 \\
0 & z_2/k
            \end{bmatrix} = z_1
            \begin{bmatrix}
1 & 0 \\
0 & 0
            \end{bmatrix} + z_2
            \begin{bmatrix}
0 & 0 \\
0 & 1/k
            \end{bmatrix}\text{.}
        \end{equation*}
All of these polynomials have total degree equal one, they are hyperstable with respect to the product $H_0^2$ and the sequence $P_k(z_1, z_2)$ converges coefficient-wise to the matrix polynomial
        \begin{equation*}
P(z_1, z_2) =
            \begin{bmatrix}
z_1 & 0 \\
0 & 0
            \end{bmatrix}\text{,}
        \end{equation*}
which is neither hyperstable with respect to $H_0^2$ nor identically equal zero.
    \end{example}

\section{Polarisation preserves hyperstability}


 We will apply Theorem \ref{GWS} to the case when $\Omega=D^\kappa$, where $D$ is an open or closed disc or open or closed half-plane. Note that such an $\Omega$ is convex and the assumptions on the total degree of $p(z_1, \dots , z_{\kappa})$ is not  needed.
 
Let $\kappa \in \mathbb{Z}_+$ and let $P(\lambda)=\lambda^d A_d + \lambda^{d-1} A_{d-1} + \dots + A_0$ be a matrix polynomial of the degree $d \leq \kappa$, we put $A_{d+1} = A_{d+2} = \dots = A_{\kappa} = 0$. We define the \emph{polarisation operator} $T_{\kappa}$ of the polynomial $P(\lambda)$ by 
		\begin{equation}\label{Tdef1}
(T_{\kappa}P)(z_1, z_2, \dots , z_{\kappa}):= \sum_{j=0}^{\kappa} \binom{\kappa}{j}^{-1}s_j(z_1, z_2, \dots , z_{\kappa})A_j, 
		\end{equation}
where the symbol $s_j$ denotes the $j$-th elementary symmetric polynomial
		\begin{equation}\label{Tdef2}
s_0(z_1, z_2, \dots, z_{\kappa}) := 1,\;\;\; s_j(z_1, z_2, \dots , z_{\kappa}) := \sum_{1 \leq i_1 < i_2 < \dots < i_j \leq \kappa} z_{i_1}z_{i_2} \dots z_{i_j}.
		\end{equation}
The polarisation operator $T_\kappa$ defined above  is a well known object in the theory of  scalar multivariate polynomials, see, e.g., \cite{BorB09}. Above we have extended its action onto matrix polynomials. Note that the operator $T_\kappa$ is injective with its image being the set consisting of all symmetric multi-affine polynomials (with matrix coefficients). We have the following result. 	



	\begin{theorem}\label{Tkappa2}
Let  $\kappa \in \mathbb{Z}_+$ and let the operator $T_\kappa$ be defined by \eqref{Tdef1} above with $D$ being an open or closed disc or open or closed half-plane. Then, a one variable matrix polynomial $P(\lambda)$ is hyperstable with respect to  $D$ if and only if a $\kappa$-variable matrix polynomial $(T_\kappa P)(z_1, z_2, \dots , z_{\kappa})$ is hyperstable with respect to $D^{\kappa}$.
	\end{theorem}
	\begin{proof}
Consider first the case $n=1$. For that purpose, take a scalar polynomial $p(\lambda)$ and suppose that there exists a point $(\zeta_1, \zeta_2, \dots , \zeta_{\kappa}) \in D^{\kappa}$ such that $(T_{\kappa}p)(\zeta_1, \zeta_2, \dots , \zeta_{\kappa}) = 0$. Since the polynomial $(T_\kappa p)(z_1, z_2, \dots z_{\kappa})$ is a symmetric multi-affine polynomial, then from Theorem~\ref{GWS} we conclude that there exists  a point $\zeta_0 \in D$ such that $(T_\kappa p)(\zeta_1, \zeta_2, \dots , \zeta_{\kappa}) = (T_\kappa p)(\zeta_0, \zeta_0, \dots , \zeta_0)$. However, note that
$ (T_\kappa p)(\zeta_0, \zeta_0, \dots , \zeta_0)=p(\zeta_0)$, which shows the forward implication.
The converse implication is obvious, as $(T_{\kappa}p)(\lambda, \lambda, \dots , \lambda) = p(\lambda)$.	
	
Assume now that $P(\lambda)$ is a matrix polynomial, hyperstable with respect to $D$. Take an arbitrary  vector $x \in \mathbb{C}^n \setminus\set0$. By definition of hyperstability, there exists $y\in \mathbb{C}^n \setminus\set0$ such that the scalar polynomial $p(\lambda)=y^*P(\lambda)x$ is stable with respect to $D$. By the first part of the proof,  the polynomial $(T_\kappa p)(z_1,\dots,z_\kappa)$ is stable with respect to $D^\kappa$. Observe that 
$$
(T_\kappa p)(z_1,\dots,z_\kappa)= y^* \left( (T_\kappa P   )(z_1,\dots,z_\kappa)\right)x,
$$
 which shows that  $(T_\kappa P   )(z_1,\dots,z_\kappa)$ is hyperstable with respect to $D^\kappa$.  
 The converse implication is again obvious.
\end{proof}

One may wonder if the following version of Theorem~\ref{Tkappa2} is true: if $P(\lambda)$ is stable with respect to $D$, then $(T_\kappa P)(z_1,\dots,z_\kappa)$ is stable with respect to $D^{\kappa}$. However, this is not the case.

\begin{example}\label{nonstab}\rm  We continue with the matrix polynomial $P(\lambda)$ from Example \ref{exa}. Again, the polynomial $P(\lambda)$ is stable with respect to any $D$ but 
$$
\det((T_2 P)(z_1,z_2))=\det\mat{cc} 1 & \frac{z_1+z_2}2 \\ \frac{z_1+z_2}2 & z_1z_2+1\rix=1-\left( \frac{z_1-z_2}2\right)^2. 
$$
Hence, every ordered pair $(\mu_1, \mu_2)$ with $\mu_1 - \mu_2 = \pm 2$ is an eigenvalue, in particular a matrix polynomial $(T_2P)(z_1,z_2)$ is not stable with respect to $D^2 = H_0^2$. 
\end{example}

We provide now a  general tool for creating hyperstable univariate matrix polynomials, based on  the operator $T_\kappa$. Its applications will be given in the next Section. 

\begin{theorem}\label{increasedegree}
Let $D$ be an open or closed disc or open or closed half-plane. If a matrix polynomial $P(\lambda)$ of degree $d $ is hyperstable with respect to $D$, then for any scalar polynomials $p_1(\lambda),$ $ \dots,$ $p_\kappa(\lambda)$, $\kappa\geq d$, the matrix polynomial $Q(\lambda) := (T_{\kappa}P)(p_1(\lambda),  \dots , p_\kappa(\lambda))$ is hyperstable with respect to
$E := p_1^{-1}(D) \cap \dots \cap p_\kappa^{-1}(D)\subseteq\Comp$.
\end{theorem}
	
\begin{proof}
Fix a nonzero vector $x\in\Comp^n\setminus\{0\}$. By Theorem~\ref{Tkappa2}, the multivariate polynomial $(T_{\kappa}P)(z_1, \dots , z_\kappa)$ is hyperstable with respect to $D^\kappa$, i.e., there exists $y\in\Comp^n\setminus\{0\}$ such that $y^*(T_{\kappa}P)(z_1,  \dots , z_\kappa)x \neq 0$ for all $z_1,  \dots, z_\kappa \in D$. 
In particular, for any $\lambda\in\Comp$ such that $p_j(\lambda)\in D$ ($j=1,\dots,\kappa$) one has 
	$y^*(T_{\kappa}P)(p_1(\lambda),\dots , p_\kappa(\lambda))x \neq 0$. Hence, $Q(\lambda)$ is hyperstable with respect to $E$.
\end{proof}

\section{Polarisation and singularity}

As it was shown above hyperstability of univariate matrix polynomials implies its regularity.
The notion of singular multivariate matrix polynomials is not yet defined, and we restrict ourselves from defining it in this Dissertation. 
The reasons of our hesitation can be found below.

    \begin{proposition}
Let $\kappa \in \mathbb{Z}_+$ and $A, B \in \mathbb{C}^{n, n}$. If $\det(\lambda A + B) \equiv c \in \mathbb{C}$, then $\det T_{\kappa}(\lambda A + B) \equiv c$. In particular, if a matrix pencil $\lambda A + B$ is singular, then the determinant of a $\kappa$-variable matrix polynomial $T_{\kappa}(\lambda A + B)$ is identically zero.
    \end{proposition}
    \begin{proof}
Since the determinant of a pencil $\lambda A + B$ is constant, it does not depend on $\lambda$. In particular, we have 
$$
 \det\left ( \frac{\mu_1+\dots+\mu_\kappa}\kappa A + B\right)  = c 
 $$
 for any $\mu_1,\dots\mu_\kappa\in\Comp $. However, 
 $$
 T_{\kappa}(\lambda A + B)(z_1,\dots,z_\kappa) = (T_\kappa(\lambda) A + B) = \frac{z_1+\cdots+z_\kappa}\kappa A + B,
 $$
 which clearly has constant in $z_1,\dots,z_\kappa$ determinant. 
    \end{proof}
    The proposition above extends only in some sense to polynomials of an arbitrary degree, as the following Proposition and Example show.
    \begin{proposition}
If $\kappa \in \mathbb{Z}_+$ and $P(\lambda) \in \mathbb{C}_{\kappa}^{n, n}[\lambda]$, then a matrix polynomial $P(\lambda)$ is singular if and only if $\det T_{\kappa}(P)(\zeta, \zeta, \dots , \zeta) = 0$ for all $\zeta \in \mathbb{C}$.
    \end{proposition}
    \begin{proof}
By the definition of singularity, the polynomial $P(\lambda)$ is singular if and only if $\det P(\zeta) = 0$ for all $\zeta \in \mathbb{C}$. Note that taking $z_1 = z_2 = \dots = z_{\kappa} = \zeta$ we have $T_{\kappa}(P)(\zeta, \zeta, \dots , \zeta) = P(\zeta)$, then the assertion follows.
    \end{proof}
    \begin{example}
\rm Let us return to the singular quadratic polynomial $P(\lambda)$ from the Example \ref{sing}. Now, we calculate the value of the polarization operator $T_2$ for $P(\lambda)$ 
        \begin{equation*}
T_2(P)(z_1, z_2) = T_2
            \begin{bmatrix}
\lambda^2 & \lambda\\
\lambda & 1
            \end{bmatrix}
 = 
            \begin{bmatrix}
z_1z_2 & \frac{z_1 + z_2}{2}\\
\frac{z_1 + z_2}{2} & 1
            \end{bmatrix}, 
\end{equation*}
and its determinant
\begin{equation*}
\det
        \begin{bmatrix}
z_1z_2 & \frac{z_1 + z_2}{2}\\
\frac{z_1 + z_2}{2} & 1
        \end{bmatrix}
 = -\frac{(z_1 - z_2)^2}{4} \not\equiv 0\text{.}
        \end{equation*}
We see that the determinant of a 2-variable matrix polynomial $T_2(P)(z_1, z_2)$ is not identically zero, although $\det P(\lambda)$ is.
    \end{example}

\chapter{Stability and hyperstability of some classes of matrix polynomials}\label{sPos}

\section{Basic methods for hyperstability}

In this chapter, we will apply a few theorems concerning stability and hyperstability of matrix poylnomials developed in previous chapters of the Thesis. These results provide examples of basic methods for proving hyperstability of certain polynomials. First, let us deal with a relatively simple matrix polynomial, appearing in \cite{KalN19,MehMW22}.
	\begin{proposition}\label{MGT} 
Let $P(\lambda)=\lambda^3 I_n +a I_n \lambda^2 +\lambda b R + cR$ with  $R\in\Comp^{n,n}$ positive definite. If $a>1$ and $b>c$ then $P(\lambda)$ is hyperstable with respect to the open right half-plane $H_{\pi/2}$.
    \end{proposition}
    \begin{proof}
It was showed in  \cite{MehMW22}  that $P(\lambda)$ is stable with respect to the open right half-plane. Since $P(\lambda) = (\lambda^3 + a\lambda^2)I_n + (\lambda b + c)R$, we can directly use Theorem~\ref{uppert}\eqref{pq} to show that it is also hyperstable with respect to $H_{\pi/2}$. 
    \end{proof}
Chapter ~\ref{sHyper} can be illustrated with more complicated examples. Besides Theorem ~\ref{uppert}, we can apply Proposition ~\ref{upperblock} as well. The following example presents how we use them altogether in order to show hyperstability of a matrix polynomial with a block upper-triangular structure. 
    \begin{example}\rm
Let us consider the following $4 \times 4$ block upper-triangular matrix polynomial
        \begin{equation*}
P(\lambda) = 
\left[\begin{array}{@{}c|c@{}}
  \begin{matrix}
\lambda & 1 \\
1 & \lambda
  \end{matrix}
  & 
  \begin{matrix}
\lambda^4 - 1 & \lambda^2 - 1 \\
\lambda^2 + 1 & 1
  \end{matrix} \\
\hline
  \begin{matrix}
0 & 0 \\
0 & 0
  \end{matrix} &
  \begin{matrix}
\lambda - 1 & 0 \\
0 & \lambda + 1
  \end{matrix}
    \end{array}\right]\text{.}
        \end{equation*}
Note that on the main diagonal of $P(\lambda)$ we have two matrix pencils which are stable, and consequently hyperstable with respect to the open upper half-plane $H_0$ (cf. Theorem ~\ref{uppert}(ii) for matrix pencils). According to Proposition ~\ref{upperblock} matrix polynomial $P(\lambda)$ is also hyperstable with respect to $H_0$, regardless the fact that the polynomial in the upper right corner of $P(\lambda)$ is singular.
    \end{example}

\section{Hyperstability via bivariate stability}

In this section we will apply two main theorems from Chapter~\ref{s5}. Recall that it requires considering a bivariate polynomial $Q(z_1, z_2)$ of a special form, satisfying $Q(\lambda,\lambda) = P(\lambda)$.  
For this polynomial $Q(z_1, z_2)$ we only need to show it is \emph{stable} with respect to some $D^2$. Then the linear algebra machinery from the proofs of Theorems~\ref{poly2} and~\ref{poly3} shows that a matrix polynomial $P(\lambda)$ is hyperstable with respect to $D$, which is a stronger notion than stability, cf. Proposition ~\ref{abc}. 
Note that the stability of the polynomial in question is in some instances obvious or was shown before. However, by proving hyperstablity we gain some extra techniques, developed in Chapter ~\ref{sHyper}, that allow us to manipulate with the polynomial and receive further results.

Stability of $P(\lambda)$ below is almost obvious (follows from the subadditivity of the norm), but showing hyperstability requires Theorem~\ref{poly2}.

    \begin{proposition}\label{subadd}
Let $A_0, A_1, A_2 \in \mathbb{C}^{n, n}$ and let $D = \set{\lambda \in \Comp: |\lambda| < r}$, $r > 0$ be such that $r\norm{A_1} + r^2\norm{A_2} < \sigma_{\min}(A_0)$. Then the bivariate polynomial $Q(z_1, z_2) = z_1^2 A_2 + z_2 A_1 + A_0$ is stable with respect to $D^2$ and consequently the regular matrix polynomial $P(\lambda) = \lambda^2 A_2 + \lambda A_1 + A_0$ is hyperstable with respect to $D$.
    \end{proposition}
    \begin{proof}
Note that by definition, stability of $Q(z_1, z_2)$ implies stability of $P(\lambda) = Q(\lambda, \lambda)$. In particular, $P(\lambda)$is regular. We have
$$
\norm{z_1^2 A_2 + z_2 A_1} < r\norm {A_1} + r^2\norm{A_2} < \sigma_{\min}(A_0)
$$
for $(z_1, z_2) \in D^2$, which implies the stability of the bivariate polynomial $Q(z_1, z_2) = z_1^2 A_2 + z_2 A_1 + A_0$ with respect to $D^2$. Indeed, if for some $(\mu_1, \mu_2) \in D^2$ and $x_0 \in \mathbb{C}^2\setminus\{0\}$ it was 
$$
(\mu_1^2 A_2 + \mu_2 A_1 + A_0)x_0 = 0\text{,}
$$
we would have
$$
\norm{(\mu_1^2 A_2 + \mu_2 A_1)x_0} = \norm{A_0x_0}\text{,}
$$
which would lead us to 
$$
\norm{A_0x_0} \leq \norm{\mu_1^2 A_2 + \mu_2 A_1}\norm{x_0} < \sigma_{\min}(A_0)\norm{x_0}\text{.}
$$
The inequality 
$$
\norm{A_0x_0} < \sigma_{\min}(A_0)\norm{x_0}
$$ 
is contradictory, because in fact the opposite inequality holds:
$$
\norm{A_0x} \geq \sigma_{\min}(A_0)\norm{x}
$$
for all $x \in \mathbb{C}^n$. Application of Theorem~\ref{poly2}\eqref{0?D} finishes the proof.
    \end{proof}
We want to present also second version of Proposition ~\ref{subadd}, which gives us the same result, but from different inequality in the assumption.
    \begin{proposition}\label{subadd2}
Let $A_0, A_1, A_2 \in \mathbb{C}^{n, n}$ and let $D = \set{\lambda \in \Comp: |\lambda| \geq r}$, $r > 0$ be such that $r\norm{A_1}+\norm{A_0}< r^2\sigma_{\min}(A_2)$. Then the bivariate polynomial $Q(z_1, z_2) = z_1^2 A_2 + z_2 A_1 + A_0$ is stable with respect to $D^2$ and consequently the regular matrix polynomial $P(\lambda) = \lambda^2 A_2 + \lambda A_1 + A_0$ is hyperstable with respect to $D$.
    \end{proposition}
    \begin{proof}
As before, we have
$$
\norm{z_2 A_1 + A_0} < r\norm {A_1} + \norm{A_0} < r^2\sigma_{\min}(A_2),
$$
which implies the stability of the bivariate polynomial $Q(z_1, z_2) = z_1^2 A_2 + z_2 A_1 + A_0$ with respect to $D^2$.
Again, application of Theorem~\ref{poly2}\eqref{0?D} ends the proof. 
    \end{proof}


Next, we present some stronger result for quadratic polynomials than one originating form Theorem ~\ref{quad}. Namely, mentioned theorem says that they are stable with respect to the right half-plane $H_{\pi/2}$, but we are going to show their hyperstablity with respect to the right half-plane $H_{\pi/2}$.

\begin{theorem}\label{half-plane} Let $R_j \in \mathbb{C}^{n, n}$ $(j=0,1,2)$  be  Hermitian positive semi-definite and let $J\in \mathbb{C}^{n, n}$ be skew-Hermitian. Assume that 
$$
\ker R_0\cap\ker R_1\cap \ker R_2\cap\ker J=\set0\text{.}
$$
 Then the polynomial $\tilde P(z_1, z_2) = z_1z_2 R_2 +z_2 (J+R_1) +R_0$ is stable with respect to the Cartesian square of the open right half-plane  $H_{\pi/2}^2$. In consequence,  $P(\lambda) = \lambda^2 R_2 + \lambda (J+R_1) + R_0$ is a regular polynomial, hyperstable with respect to $H_{\pi/2}$.
 \end{theorem}
\begin{proof}
We will make use of Theorem \ref{poly2} (b), note that $0 \not\in D = H_{\pi/2}$. 
Consider the polynomial 
$$
\tilde P(z_1,z_2)=z_1z_2R_2+z_2 (J+R_1) +R_0\text{,}
$$
and suppose it is not stable, i.e., $\tilde P(\mu_1,\mu_2)x=0$ for some  $(\mu_1,\mu_2) \in H_{\pi/2}^2$ and $x\neq 0$. Multiplying from the left by $x^*$ and taking the real part one obtains
\begin{equation}\label{eigenvalue1}
\re(\mu_1) x^*R_2 x +  x^*R_1x + \re\left(\frac1{\mu_2}\right) x^*R_0 x=0.
\end{equation}
 Since both $\re( \mu_1)$ and $\re(\frac1{\mu_2})$ are positive and $R_2,R_1,R_0$ are positive semi-definite, we obtain  $x^*R_2 x =  x^*R_1x  = x^*R_0 x=0$. Hence, $R_2x=R_1x=R_0x=0$.
 But this implies that $0=\tilde P(\mu_1,\mu_2)x=Jx$, 
 a contradiction.

	
	\end{proof}

    \begin{example}\rm
Let us take
        \begin{eqnarray*}
\;\;&&R_0 =
            \begin{bmatrix}
1 & 1 & 0 \\
1 & 2 & 1 \\
0 & 1 & 1
            \end{bmatrix}\text{,}\;
R_1 =
            \begin{bmatrix}
1 & 1 & 1 \\
1 & 1 & 1 \\
1 & 1 & 1
            \end{bmatrix}\text{,}\;
R_2 =
            \begin{bmatrix}
2 & -1 & 0 \\
-1 & 2 & -1 \\
0 & -1 & 2
            \end{bmatrix}\text{,}\;
J =
            \begin{bmatrix}
0 & 1 & 2 \\
-1 & 0 & 1 \\
-2 & -1 & 0
            \end{bmatrix}\text{.} \\
        \end{eqnarray*}
We have $\sigma(R_0) = \{0, 1, 3\}, \sigma(R_1) = \{0, 3\}, \sigma(R_2) = \{2-\sqrt{2}, 2, 2 + \sqrt{2}\}$, hence $R_0, R_1, R_2$ are positive semi-definite. Moreover $J^* = -J$, i.e. $J$ is skew-Hermitian. Note that $\ker R_2 = \{0\}$ is trivial, therefore the condition $\ker R_0\cap\ker R_1\cap\ker R_2\cap\ker J = \{0\}$ is also satisfied. According to Theorem ~\ref{half-plane} a two-variable matrix polynomial
        \begin{eqnarray*}
&&P(z_1, z_2) = z_1 z_2
            \begin{bmatrix}
2 & -1 & 0 \\
-1 & 2 & -1 \\
0 & -1 & 2
            \end{bmatrix} 
+ z_2
            \begin{bmatrix}
1 & 2 & 3 \\
0 & 1 & 2 \\
-1 & 0 & 1
            \end{bmatrix} 
+ 
            \begin{bmatrix}
1 & 1 & 0 \\
1 & 2 & 1 \\
0 & 1 & 1
            \end{bmatrix}\\
        \end{eqnarray*}
is stable with respect to $H_{\pi/2}^2$ and a one-variable matrix polynomial
        \begin{eqnarray*}
&&P(\lambda) = \lambda^2
            \begin{bmatrix}
2 & -1 & 0 \\
-1 & 2 & -1 \\
0 & -1 & 2
            \end{bmatrix} 
+ \lambda
            \begin{bmatrix}
1 & 2 & 3 \\
0 & 1 & 2 \\
-1 & 0 & 1
            \end{bmatrix} 
+ 
            \begin{bmatrix}
1 & 1 & 0 \\
1 & 2 & 1 \\
0 & 1 & 1
            \end{bmatrix}\\
        \end{eqnarray*}
is regular, hyperstable with respect to $H_{\pi/2}$.
    \end{example}
 
    \begin{remark}\rm
Let us note that the condition $\ker R_1\cap\ker R_2\cap \ker R_3\cap\ker J=\set0$ is equivalent to the regularity of $P(\lambda)$. For a proof, recall Theorem ~\ref{kernels}.
    \end{remark}
    \begin{remark}\rm
Note that although in the proof of Theorem ~\ref{half-plane}, we multiplied $P(\lambda)$ from the left by $x^*$, nevertheless we have \emph{not} shown there that the numerical range of $P(\lambda)$ is located outside the open right half-plane $H_{\pi/2}$.
    \end{remark}
	
    \begin{corollary}\label{c-hp}
Let $R\in\mathbb C^{n,n}$  be  Hermitian positive semi-definite, let $J\in \mathbb{C}^{n, n}$ be skew-Hermitian and let  $Q, A_0, A_2 \in \mathbb{C}^{n, n}$ be such that $Q^*A_2$ and $Q^*A_0$ are Hermitian positive semi-definite. Assume also that 
$$
\ker(Q^* A_0)\cap\ker(Q^* R Q)\cap\ker(Q^*JQ)\cap\ker(Q^*A_2)=\{0\}\text{.}
$$
Then the matrix polynomial  $P(\lambda) = \lambda^2 A_2 + \lambda (J+R)Q + A_0$ is a regular polynomial, hyperstable  with respect to the open right half-plane $H_{\pi/2}$.
    \end{corollary}

	
	\begin{proof}
Firstly, observe that the matrices $Q^*A_0,\; Q^*RQ,\; Q^*A_2,\; Q^*JQ$ satisfy the assumptions of Theorem ~\ref{half-plane}, i.e. $Q^*A_0,\; Q^*RQ,\; Q^*A_2 \geq 0$ and $ (Q^*JQ)^* = -Q^*JQ$. Therefore, a polynomial $\lambda^2 Q^*A_2 + \lambda(Q^*JQ + Q^*RQ) + Q^*A_0 = Q^*P(\lambda)I$ is hyperstable with respect to $H_{\pi/2}$. Secondly, we apply Lemma ~\ref{lQ} and we get hyperstability of $P(\lambda)$ with respect to $H_{\pi/2}$.
	\end{proof}
	
Let us mention that in the following results, the argument of a complex number $\lambda \in \mathbb{C}$ is taken such that $\Arg\lambda \in (-\pi; \pi]$ and, for the sake of simplicity, we set $\Arg 0 := 0$.
	    \begin{theorem}\label{ker}
Let $R_j \in \mathbb{C}^{n, n}$ $(j=1,2,3)$  be  Hermitian positive semi-definite and let $A_0\in\mathbb C^{n,n}$ be Hermitian and let $G\in \mathbb{C}^{n, n}$ be skew-Hermitian with $-\ii G$ positive semi-definite. Assume that $\ker G \cap \ker A_0\cap \ker R_1\cap\ker R_2\cap \ker R_3 = \set0$.
 Then the following holds:
        \begin{enumerate}[\rm (i)]
\item  the multivariate matrix polynomial 
$$
P_1(z_1, z_2)= (z_1^3 z_2^3 + z_1^3 + z_2^3)R_3 + (z_1^2 z_2^3 +z_1^2)R_2 + (z_1^3z_2 + z_2) R_1 + A_0+G
$$ 
is stable with respect to $D_1^2$, where $D_1 = \{\lambda \in \mathbb{C} : 0 < \Arg\lambda < \pi/6\}$;
\item\label{P2} the multivariate matrix polynomial 
$$
P_2(z_1, z_2) = z_2^3 R_3 + z_1 z_2 R_2 + z_2 R_1 + A_0+G
$$ 
is stable with respect to $D_2^2$, where $D_2 = \{\lambda \in \mathbb{C} : 0 < \Arg\lambda < \pi/3\}$;
\item the multivariate matrix polynomial 
$$
P_3(z_1, z_2) = z_1 z_2^3 R_3 + z_1 z_2^2 R_2 + z_2^2 R_1 +  z_1 (A_0 + G)
$$ 
is stable with respect to $D_3^2$, where $D_3 = \{\lambda \in \mathbb{C} : 0 < \Arg\lambda < \pi/4\}$.
        \end{enumerate}
In particular, $P(\lambda) = \lambda^3 R_3+\lambda^2 R_2+\lambda R_1 + (A_0 + G)$ is stable with respect to $D = \{\lambda \in \mathbb{C} : 0 < \Arg\lambda < \pi/3\}$.
    \end{theorem}
    \begin{proof}
We show only (i), the proofs of (ii) and (iii) are similar. Suppose that the matrix polynomial $P_1(z_1, z_2)$ has an eigenvalue $(\mu_1, \mu_2) \in D_1^2$. Thus there exist a nonzero vector $x \in \mathbb{C}^n\setminus\{0\}$ such that
        \begin{equation}\label{mmm}
(\mu_1^3\mu_2^3 + \mu_1^3 + \mu_2^3)R_3x + (\mu_1^2\mu_2^3 + \mu_1^2)R_2x + (\mu_1^3\mu_2 + \mu_2)R_1x + (A_0+G)x = 0\text{.}
        \end{equation}
Multiplying by $x^*$ and taking the imaginary part of both sides of the equation above we obtain
        \begin{eqnarray*}
&&(x^*R_3x)\;\im(\mu_1^3\mu_2^3 + \mu_1^3 + \mu_2^3)\\
&+& (x^*R_2x)\;\im(\mu_1^2\mu_2^3 + \mu_1^2) \\
&+& (x^*R_1x)\;\im(\mu_1^3\mu_2 + \mu_2)- x^*(\ii G)x \\
&=& 0\text{.}
        \end{eqnarray*}
Note that by assumption $x^*R_3x, x^*R_1x, x^*R_2x, -x^*(\ii G)x \geq 0$ and since $0< \Arg \mu_1, \Arg \mu_2 < \pi/6$, we have 
$$
\im(\mu_1^3\mu_2^3 + \mu_1^3 + \mu_2^3), \im(\mu_1^2\mu_2^3 + \mu_1^2), \im(\mu_1^3\mu_2 + \mu_2) > 0\text{.}
$$
Therefore, we have $x^*R_3x = x^*R_2x = x^*R_1x = x^*(iG)x = 0$ and consequently  $R_3x = R_2x = R_1x = G x = 0$. Due to the equation \eqref{mmm} we have $A_0x = 0$, a contradiction. 
The `In particular' part follows by substituting $\lambda$ for $z_1$ and $z_2$ in \eqref{P2}.
    \end{proof}

The authors of ~\cite{MehMW22} have shown that a regular cubic matrix polynomial with all coefficients Hermitian positive semi-definite is stable with respect to $D = \{\lambda \in \mathbb{C} : -\pi/3 < \Arg\lambda < \pi/3\}$, cf. ~\ref{cube} in the preliminaries. Below, we present a connected result as a collorary from Theorem ~\ref{poly3} and ~\ref{ker}\eqref{P2}.
Directly from these theorems, we obtain the following result.

    \begin{corollary} 
Let $R_j \in \mathbb{C}^{n, n}$ $(j=0,1,2)$  be  Hermitian positive semi-definite, then the polynomial 
$P(\lambda) = \lambda^3 R_0 + \lambda^2 R_2 + \lambda  R_1 + R_0$ is hyperstable with respect to $D = \{\lambda \in \mathbb{C} : 0 < \Arg\lambda < \pi/3\}$. Note that each scalar cubic polynomial with non-negative coefficients is stable with respect to $D$.
    \end{corollary}

\section{Hyperstability through polarisation}
    
We show now how a polarisation operator may be used to increase the degree of the polynomial. The price for which is narrowing the set $D$.
	
    \begin{corollary}\label{deg3} 
Let $R_j \in \mathbb{C}^{n, n}$ $(j=0,1,2)$ be Hermitian positive semi-definite, and let $J\in\Comp^{n,n}$ be skew-Hermitian. Then the matrix polynomial $Q(\lambda) = \lambda^3 R_2 + (\lambda^2+\lambda)(R_1+J)  + R_0$ is a regular cubic matrix polynomial, hyperstable with respect to the angle $E = \{\lambda\in\Comp:-\pi/4 < \Arg\lambda <\pi/4\}\setminus\{0\}$.
    \end{corollary}
	\begin{proof} 
By Theorem~\ref{half-plane} we obtain the polynomial $P(\lambda)= \lambda^2 R_2 + 2 \lambda  (R_1+J) + R_0$ hyperstable with respect to the open right half-plane $H_{\pi/2}$. We apply now Theorem~\ref{increasedegree} with $p_1(\lambda)=\lambda^2$, $p_2(\lambda)=\lambda$ obtaining 
$$
(T_2P)(z_1,z_2)=  z_1z_2 A_2 + (z_1+z_2) A_1 + A_0
$$
so that $(T_2P)(\lambda^2,\lambda) = Q(\lambda)$. Finally, observe that the angle $E$ is precisely the set $p_1^{-1}(H_{\pi/2})\cap p_2^{-1}(H_{\pi/2})$ from Theorem~\ref{increasedegree}.
	\end{proof}
	
The operator $T_2$ may be used to increase or decrease the degree of a polynomial, i.e., obtain higher or lower degree polynomial from the given polynomial $P(\lambda)$. In Corollary ~\ref{deg3}, we get a qubic polynomial from the quadratic hyperstable polynomial $P(\lambda)$. Let us take now an opposite action and use $T_2$ to obtain the linear polynomial (matrix pencil) from $P(\lambda)$ of degree $2$.

    \begin{corollary}
Let $R_j \in \mathbb{C}^{n, n}$ $(j=0,1)$  be  Hermitian positive semi-definite, and let $J\in\Comp^{n,n}$ be skew-Hermitian. Consider a matrix pencil $P(\lambda)=\lambda( R_1+J) + (R_0+a J)$, where $a\geq0$. Then the eigenvalues of $P(\lambda)$ are contained in the closed left half-plane $\overline{H}_{\frac{3}{2}\pi}$. 
    \end{corollary}

\begin{proof}
The case $a=0$ was considered in \cite{MehMW22}. Now take $a>0$ and define a matrix polynomial $P(\lambda)=\lambda^2\frac{R_1}a+\lambda \cdot 2 J+ R_0$, note that it clearly satisfies the assumptions of Theorem~\ref{half-plane}. Hence it is hyperstable with respect to $H_{\pi/2}$, and by Theorem~\ref{Tkappa2} the matrix polynomial 
$$
(T_2 P)(z_1,z_2)=z_1z_2 \frac{R_1}a+(z_1+z_2) J+ R_0
$$
is hyperstable with respect to $H_{\pi/2}^2$. In particular, if we set $z_2=a$ and replace $z_1$ by $\lambda$, then we obtain the original polynomial 
$\lambda a \frac{R_1}a+ (\lambda+a)J+ R_0=\lambda A_1+A_0$
being stable with respect to $H_{\pi/2}$.
\end{proof}


   \bibliographystyle{plain}

\bibliography{szw}

\begin{thebibliography}{10}

\bibitem{ang2021}
Luis~M Anguas, Froil{\'a}n~M Dopico, Richard Hollister, and D~Steven Mackey.
\newblock Quasi-triangularization of matrix polynomials over arbitrary fields.
\newblock {\em Linear Algebra and its Applications}, 665:61--106, 2023.

\bibitem{BorB09}
Julius Borcea and Petter Br{\"a}nd{\'e}n.
\newblock The {L}ee-{Y}ang and {P}{\'o}lya-{S}chur programs.{I}. {L}inear
  operators preserving stability.
\newblock {\em Inventiones mathematicae}, 177(3):541, 2009.

\bibitem{brown1999proof}
Johnny~E Brown and Guangping Xiang.
\newblock Proof of the {S}endov conjecture for polynomials of degree at most
  eight.
\newblock {\em Journal of mathematical analysis and applications},
  232(2):272--292, 1999.

\bibitem{deB61}
Louis De~Branges.
\newblock Some {H}ilbert spaces of entire functions. {II}.
\newblock {\em Transactions of the American Mathematical Society},
  99(1):118--152, 1961.

\bibitem{Gan59}
F.~R. Gantmacher.
\newblock {\em Theory of Matrices}.
\newblock Chelsea, New York, 1959.

\bibitem{gol2013}
Gene~H Golub and Charles~F Van~Loan.
\newblock {\em Matrix computations}.
\newblock JHU press, 2013.

\bibitem{grace1902}
John~H Grace.
\newblock The zeros of a polynomial.
\newblock In {\em Proc. Cambridge Philos. Soc}, volume~11, pages 352--357,
  1902.

\bibitem{HorJ85}
R.~A. Horn and C.~R. Johnson.
\newblock {\em Matrix Analysis}.
\newblock Cambridge University Press, Cambridge, 1985.

\bibitem{HorJ91}
R.~A. Horn and C.~R. Johnson.
\newblock {\em Topics in Matrix Analysis}.
\newblock Cambridge University Press, Cambridge, 1991.

\bibitem{KalN19}
B.~Kaltenbacher and V.~Nikolic.
\newblock On the {J}ordan-{M}oore-{G}ibson-{T}hompson equation: Well-posedness
  with quadratic gradient nonlinearity and singular limit for vanishing
  relaxation time.
\newblock {\em Math. Models Methods Appl. Sci.}, 29:2523--2556, 2019.

\bibitem{kanter}
Marek Kanter.
\newblock Multivariate {G}auss-{L}ucas theorems.
\newblock {\em arXiv preprint arXiv:1203.6426}, 2012.

\bibitem{Kne19}
G.~Knese.
\newblock Global bounds on stable polynomials.
\newblock {\em Complex Analysis and Operator Theory}, 13(4):1895--1915, 2019.

\bibitem{PP2023matrix}
Vadym Koval and Patryk Pagacz.
\newblock Matrix pencils with the numerical range equal to the whole complex
  plane.
\newblock {\em Linear Algebra and its Applications}, 657:274--286, 2023.

\bibitem{LiR94}
C.-K. Li and L.~Rodman.
\newblock Numerical range of matrix polynomials.
\newblock {\em SIAM Journal on Matrix Analysis and Applications},
  15(4):1256--1265, 1994.

\bibitem{MehMW18}
C.~Mehl, V.~Mehrmann, and M.~Wojtylak.
\newblock Linear algebra properties of dissipative {H}amiltonian descriptor
  systems.
\newblock {\em {SIAM} J. Matrix Anal. Appl.}, 39(3):1489--1519, 2018.

\bibitem{MehMW21}
C.~Mehl, V.~Mehrmann, and M.~Wojtylak.
\newblock Distance problems for dissipative {H}amiltonian systems and related
  matrix polynomials.
\newblock {\em Linear Algebra Appl.}, 623:335--366, 2021.

\bibitem{MehMW22}
Christian Mehl, Volker Mehrmann, and Michal Wojtylak.
\newblock Matrix pencils with coefficients that have positive semidefinite
  {H}ermitian parts.
\newblock {\em SIAM Journal on Matrix Analysis and Applications},
  43(3):1186--1212, 2022.

\bibitem{PSW}
Piotr Pikul, Oskar~Jakub Szyma{\'n}ski, and Micha{\l} Wojtylak.
\newblock The sz\'asz inequality for matrix polynomials and functional
  calculus.
\newblock {\em arXiv preprint arXiv:2406.08965}, 2024.

\bibitem{Psa03}
Panayiotis~J. Psarrakos.
\newblock Definite triples of {H}ermitian matrices and matrix polynomials.
\newblock {\em Journal of Computational and Applied Mathematics},
  151(1):39--58, 2003.

\bibitem{Psa00}
P.J. Psarrakos.
\newblock Numerical range of linear pencils.
\newblock {\em Linear Algebra and its Applications}, 317(1-3):127--141, 2000.

\bibitem{steinerberge}
Stefan Steinerberger.
\newblock A stability version of the {G}auss--{L}ucas theorem and applications.
\newblock {\em Journal of the Australian Mathematical Society},
  109(2):262--269, 2020.

\bibitem{szego1922}
G{\'a}bor Szeg{\"o}.
\newblock Bemerkungen zu einem {S}atz von {J}{H} {G}race {\"u}ber die {W}urzeln
  algebraischer {G}leichungen.
\newblock {\em Mathematische Zeitschrift}, 13:28--55, 1922.

\bibitem{szymanski2023stability}
Oskar~Jakub Szyma{\'n}ski and Micha{\l} Wojtylak.
\newblock Stability of matrix polynomials in one and several variables.
\newblock {\em Linear Algebra and its Applications}, 670:42--67, 2023.

\bibitem{tao2020sendovs}
Terence Tao.
\newblock Sendov's conjecture for sufficiently high degree polynomials.
\newblock {\em arXiv preprint arXiv:2012.04125}, 2020.

\bibitem{trian2013}
Leo Taslaman, Fran{\c{c}}oise Tisseur, and Ion Zaballa.
\newblock Triangularizing matrix polynomials.
\newblock {\em Linear Algebra and its Applications}, 439(7):1679--1699, 2013.

\bibitem{tis2013}
Fran{\c{c}}oise Tisseur and Ion Zaballa.
\newblock Triangularizing quadratic matrix polynomials.
\newblock {\em SIAM Journal on Matrix Analysis and Applications},
  34(2):312--337, 2013.

\bibitem{vinnikov}
Victor Vinnikov.
\newblock Lmi representations of convex semialgebraic sets and determinantal
  representations of algebraic hypersurfaces: past, present, and future.
\newblock {\em Mathematical Methods in Systems, Optimization, and Control:
  Festschrift in Honor of J. William Helton}, pages 325--349, 2012.

\bibitem{lecjv}
Jan Vondr\'ak.
\newblock Non-constructive methods in combinatorics. {L}ecture notes.
\newblock 2016.
\newblock https://theory.stanford.edu/$\sim$jvondrak/MATH233-2016/Math233.html.

\bibitem{walsh1922}
Joseph~L Walsh.
\newblock On the location of the roots of certain types of polynomials.
\newblock {\em Transactions of the American Mathematical Society},
  24(3):163--180, 1922.

\end{thebibliography}

\end{document}